\pdfoutput=1
\RequirePackage{ifpdf}
\ifpdf 
\documentclass[pdftex]{sigma}
\else
\documentclass{sigma}
\fi

\newcommand{\starltsa}{\star_{\ltsa}}
\usepackage{shuffle}
\usepackage{tikz}
\usepackage{tikz-cd}

\DeclareFontFamily{U}{MnSymbolC}{}
\DeclareFontShape{U}{MnSymbolC}{m}{n}{
	<-5.5> MnSymbolC5
	<5.5-6.5> MnSymbolC6
	<6.5-7.5> MnSymbolC7
	<7.5-8.5> MnSymbolC8
	<8.5-9.5> MnSymbolC9
	<9.5-11.5> MnSymbolC10
	<11.5-> MnSymbolCb12
}{}

\newcommand{\lmss}[1]{\textrm{\normalfont{{\fontfamily{lmss}\selectfont #1}}}}
\newcommand{\sbt}{\,\begin{tikzpicture}[baseline=(X.base)]%
		\node[draw, fill,black,circle, inner sep=1pt] at (0,0.1) {};
		\node[circle,inner sep=0pt,outer sep=0pt] (X){$\ $};
	\end{tikzpicture}%
	\,
}

\def\lpbt{ { \lmss{LT} } }
\def\lpbf{ { \lmss{LF} } }
\def\st{{\rm std}}

\def\gs{\sigma}

\def\cF{\mathcal{F}}

\def\cA{\mathcal{A}}

\def\cG{\mathcal{G}}

\def\C{\mathbb{C}}

\def\norm#1{\Vert #1 \Vert}

\def\1{\mathbf{1}}

\usepackage{stmaryrd}
\def\overt{\varobar}

\def\Hom{{\rm Hom}}

\numberwithin{equation}{section}

\newtheorem{Theorem}{Theorem}[section]
\newtheorem*{Theorem*}{Theorem}
\newtheorem{Corollary}[Theorem]{Corollary}
\newtheorem{Lemma}[Theorem]{Lemma}
\newtheorem{Proposition}[Theorem]{Proposition}
 { \theoremstyle{definition}
\newtheorem{Definition}[Theorem]{Definition}

\newtheorem{Example}[Theorem]{Example}
\newtheorem{Remark}[Theorem]{Remark} }

\begin{document}

\allowdisplaybreaks

\newcommand{\arXivNumber}{2102.11816}

\renewcommand{\thefootnote}{}

\renewcommand{\PaperNumber}{096}

\FirstPageHeading

\ShortArticleName{On the Signature of a Path in an Operator Algebra}

\ArticleName{On the Signature of a Path in an Operator Algebra\footnote{This paper is a~contribution to the Special Issue on Non-Commutative Algebra, Probability and Analysis in Action. The~full collection is available at \href{https://www.emis.de/journals/SIGMA/non-commutative-probability.html}{https://www.emis.de/journals/SIGMA/non-commutative-probability.html}}}

\Author{Nicolas GILLIERS~$^{\rm a}$ and Carlo BELLINGERI~$^{\rm b}$}

\AuthorNameForHeading{N.~Gilliers and C.~Bellingeri}

\Address{$^{\rm a)}$~Institut de Math\'ematiques de Toulouse, UMR5219,\\
\hphantom{$^{\rm a)}$}~Universit\'e de Toulouse, CNRS, UPS, F-31062 Toulouse, France}
\EmailD{\href{mailto:nicolas.gilliers@gmail.com}{nicolas.gilliers@gmail.com}}

\Address{$^{\rm b)}$~Technische Universit\"at Berlin, Stra{\ss}e des 17. Juni 135, 10623 Berlin, Germany}
\EmailD{\href{mailto:bellinge@math.tu-berlin.de}{bellinge@math.tu-berlin.de}}

\ArticleDates{Received January 11, 2022, in final form November 30, 2022; Published online December 09, 2022}

\Abstract{We introduce a class of operators associated with the signature of a smooth path~$X$ with values in a $C^{\star}$ algebra $\mathcal{A}$. These operators serve as the basis of Taylor expansions of solutions to controlled differential equations of interest in noncommutative probability. They are defined by fully contracting iterated integrals of $X$, seen as tensors, with the product of $\mathcal{A}$. Were it considered that partial contractions should be included, we explain how these operators yield a trajectory on a group of representations of a combinatorial Hopf monoid. To clarify the role of partial contractions, we build an alternative group-valued trajectory whose increments embody full-contractions operators alone. We obtain therefore a notion of signature, which seems more appropriate for noncommutative probability.}


\Classification{18M60; 18M80; 60L10; 46L89}

\renewcommand{\thefootnote}{\arabic{footnote}}
\setcounter{footnote}{0}

\section{Introduction}
This work intends to explore a direction suggested in \cite{deya2013rough} and aims to use paths principles for studying the following class of differential equations
\begin{equation}\label{eqn:ncequations}
	\mathrm{d}Y_{t} = a(Y_t)\cdot{\rm d}X_{t}\cdot b(Y_t),\qquad Y_0\in \cA.
\end{equation}
In the above equation, the driving path $X\colon [0,1]\to \cA$ takes values in an unital $C^{\star}$-algebra $(\mathcal{A},\cdot, \star, \|\cdot \|)$ with unity $1_{\mathcal{A}}$ and $a,b\colon \mathbb{\cA}\to\mathbb{\cA}$ are two polynomial functions or Fourier transforms of regular measures with exponential moments, see \cite{Biane01, deya2013rough}.

This paper is the first of two whose objectives are to introduce a new notion of geometric rough paths, tailored to the class of equations \eqref{eqn:ncequations}. In this work, we focus on the algebra underlying Taylor expansions of solutions to equations \eqref{eqn:ncequations}, discarding other crucial aspects (such as measurability).

\subsection{The rough paths approach}
In the nineties \cite{lyons1998}, T.J.~Lyons proposed the appropriate mathematical framework to study controlled differential equations
\begin{equation}\label{eqn:controlleddiff}
	{\rm d}Y_{t} = \sigma(Y_{t})\mathrm{d}X_{t},\qquad Y_{0} = y_{0} \in \mathbb{R}^{d}.
\end{equation}
In \eqref{eqn:controlleddiff}, the solution $Y$ is a continuous path in $\mathbb{R}^{d}$, $\sigma\colon \mathbb{R}^{d}\to \lmss{End}\big(\mathbb{R}^{n}, \mathbb{R}^{d}\big)$ is a smooth vector field and the driving path $X$ is H\"older continuous.
If $X$ is smooth, standard differential calculus provides a rigorous interpretation to \eqref{eqn:controlleddiff}.
For paths with lower regularity, Young's theory of integration~\cite{Young1936} gives sense to equation \eqref{eqn:controlleddiff} driven by an H\"older regular path $X$ with exponent greater than $\frac{1}{2}$.
Interesting stochastic driving paths are too irregular for Young integration. For instance, Brownian trajectories are only $\frac{1}{2}-\varepsilon, \,\varepsilon > 0$, H\"older continuous. Classical It\^o integration supplements limitations of Young's theory and defines integrals driven by continuous semi-martingales as limits \emph{in probability} of Riemann sums.

Rough path theory extends the standard rules of differential and integral calculus to H\"older paths $X$ and provides a \emph{pathwise} interpretation to~\eqref{eqn:controlleddiff}. Let us add more details.
Given a~smooth field $\sigma$ and $Y_0 \in \mathbb{R}^{n}$ the solution map $\Phi\colon X\mapsto Y$ to equation \eqref{eqn:controlleddiff} is continuous with respect to the Lipschitz norm on the space of smooth driving paths $X$.
A fundamental observation is the following one: by applying Picard's iterations to \eqref{eqn:controlleddiff}, one quickly reckons that the solution map $\Phi$ is a linear function of the entire \emph{signature} of $X$, that is the infinite collection of tensors,
\begin{equation}\label{equ:signature}
	\mathbb{X}_{s,t}=\bigg(1, X_t-X_s, \int_{\Delta^2_{s,t}} \mathrm{d}X_{t_1}\otimes \mathrm{d}X_{t_2}, \dots , \int_{\Delta^{n}_{s,t}} \mathrm{d}X_{t_1}\otimes \cdots \otimes \mathrm{d}X_{t_{n}}, \cdots \bigg),
\end{equation}
where $\Delta^n_{s,t}:= \{s < t_{1} < \cdots <t_{n} < t\}$ is the $n$-dimensional simplex. Signatures of smooth paths support a one-parameter family of topologies with respect to which $\Phi$ is continuous. Complete spaces for these topologies contain H\"older paths together with the additional data of an \emph{abstract signature}. These abstract signatures are called \emph{rough paths} and can alternatively be characterized by a set of algebraic and analytical properties.
Indeed, a rough path is a two parameters function $(s,t)\to \mathbb{X}_{s,t}$ with values in a \emph{group} $(G, \circ)$, included in the completed tensor space of $\mathbb{R}^n$, with the property that for each triple of times $s,u,t\in [0,1]^3$
\begin{equation}
	\label{eqn:Chen_rel}
	\mathbb{X}_{s,t}= \mathbb{X}_{s,u}\star \mathbb{X}_{u,t}.
\end{equation}
The relations \eqref{eqn:Chen_rel} are usually called \emph{Chen's relation} after Kuo-Tsai Chen \cite{Chen54} and its secular work on the homology of loop spaces. We refer the reader to the monograph \cite{Friz2020course} for a detailed exposition of rough paths theory.

\subsection{Motivation and previous works}
We choose to have an intrinsic -- coordinate-free -- approach to \eqref{eqn:ncequations} and to work consistently with the specific class of fields we consider, that is with the algebra product.
Rough paths theory on infinite-dimensional spaces is more intricate because of several notions of tensor products between two Banach algebras, see \cite{Lyons2002}. Considering the class of equations \eqref{eqn:ncequations} the \emph{projective tensor product} is the only reasonable one since the algebra product is always continuous with respect to this topology. This is not true for the spatial (or injective) topology. This limitation strikes with the results obtained in \cite{donati01, victoir2004levy}. In these works, the authors define a rough path (in fact, a L\'evy area) over the free Brownian motion in the \emph{spacial tensor product} by using free It\^o calculus. Whereas it is possible \cite{Lyons2007} to show the existence of a free L\'evy area (up to an infinitesimal loss in regularity) in the projective tensor product, an explicit procedure is missing.

 To circumvent this issue, A.~Deya and R.~Schott introduced in \cite{deya2013rough} a weaker notion of L\'evy area tailored to the class of equations~\eqref{eqn:ncequations} when the H\"older scale lies in $\big(\frac{1}{3}, \frac{1}{2}\big]$: the {\it product L\'evy area}. This object embodies the data on the small-scale behaviour of the driving path~$X$ only in the directions required to give sense to~\eqref{eqn:ncequations}. The starting point to define it is a fine analysis of~\eqref{eqn:controlleddiff} with $X$ smooth and the expansion of the solution $Y$ obtained by applying Picard iterations. Pick $A,B\in\mathcal{A}$ and consider the following example (recall that $\cdot$ denotes the product of~$\mathcal{A}$)
\begin{equation*}
	\mathrm{d}Y_{t} = (A \cdot Y_{t}) \cdot \mathrm{d}X_{t} \cdot (Y_{t} \cdot B),\qquad Y_0=1_{\mathcal{A}}.
\end{equation*}
Writing the first two steps of the Picard Iteration, we obtain
\begin{gather}
	Y_{t}=1_{\mathcal{A}}+\int_{\Delta^1_{s,t}} A\cdot\textrm{d}X_{t_{1}}\cdot B
	 +\int_{\Delta^2_{s,t}} A^2 \cdot \textrm{d}X_{t_{1}}\cdot B \cdot \textrm{d}X_{t_{2}}\cdot B\nonumber\\
\hphantom{Y_{t}=}{}
 + \int_{\Delta^2_{s,t}} A \cdot \textrm{d}X_{t_{2}}\cdot A\cdot \textrm{d}X_{t_{1}} \cdot B^2
	+ R_{s,t},\label{eqn:seriesexpansion}
\end{gather}
where $R_{s,t}$ is a remainder term satisfying $|R_{s,t}| \lesssim |t-s|^{3}$. The above equation hints at a control, at any order, of the small variations of $Y$ by the following expressions
\begin{equation}
	\label{eqn:iteratedintegrals}
	\mathbb{X}_{s,t}^{\sigma} = \int_{\Delta^n_{s,t}} A_{0} \cdot \mathrm{d}X_{t_{\sigma^{-1}(1)}}\cdot A_1 \cdots A_{n-1}\cdot \mathrm{d}X_{t_{\sigma^{-1}(n)}} \cdot A_{n},\qquad A_{0},\dots,A_{n} \in \mathcal{A}.
\end{equation}
where $\sigma$ is a permutation of $\{1,2,\dots, n\}$. The expressions in \eqref{eqn:iteratedintegrals} are values of a multilinear operator $\mathbb{X}^{\sigma}_{s,t}$, that we call full contraction operator, depending on a choice of a permutation $\sigma$. The solution of the equation \eqref{eqn:ncequations} expands over the contracted iterated integrals \eqref{eqn:iteratedintegrals} in the way alluded to above under the constraints that the Fourier transforms of $a$ and $b$ are bounded measures on the real line.
A product L\'evy area is an abstraction of the order two full contraction operators, the ones indexed by permutations of $\{1,2\}$.

We elaborate on the observation of A.~Deya and R.~Schott and extract important algebraic and analytical properties of the multilinear operators \eqref{eqn:iteratedintegrals} with the objective of developing a rough theory for the class of equations \eqref{eqn:controlleddiff} with driving noise $X$ of arbitrary low H\"older regularity. To put it shortly, the main outcome of this work is a positive answer for that and we explain it by associating to the operators \eqref{eqn:iteratedintegrals} a smooth trajectory over a group of triangular morphisms on an algebra of operators.

The main difficulties lie in writing a Chen relation for the operators \eqref{eqn:iteratedintegrals} understood as a~certain ``algebraic rule" for computing \eqref{eqn:iteratedintegrals} over an interval knowing the values of \eqref{eqn:iteratedintegrals} over a~subdivision of this interval. Consider for instance the \emph{full contraction operator}
\begin{equation*}
	\mathbb{X}^3_{s,t}(A_{0},A_{1}, A_{2}, A_3):=\int_{\Delta^3_{s,t}}A_0\cdot {\rm d}X_{t_1} \cdot A_1 \cdot {\rm d}X_{t_3} \cdot A_2 \cdot {\rm d}X_{t_2} \cdot A_3.
\end{equation*}
Then the Chasles identity implies the following deconcatenation formula:
\begin{gather*}
\mathbb{X}^3_{s,t}(A_{0},A_{1}, A_{2}, A_3) = \mathbb{X}^3_{s,u}(A_{0},A_{1}, A_{2}, A_3)+ \mathbb{X}^3_{u,t}(A_{0},A_{1}, A_{2}, A_3) \\
\hphantom{\mathbb{X}^3_{s,t}(A_{0},A_{1}, A_{2}, A_3) =}{}
+ \int_{t_{3}\in\Delta^1_{u,t}} \int_{(t_1, t_2)\in \Delta^2_{s,u}} A_0\cdot {\rm d}X_{t_1} \cdot A_1 \cdot {\rm d}X_{t_3} \cdot A_2 \cdot {\rm d}X_{t_2} \cdot A_3 \\
\hphantom{\mathbb{X}^3_{s,t}(A_{0},A_{1}, A_{2}, A_3) =}{}
+ \int_{(t_2,t_{3})\in\Delta^2_{u,t}} \int_{t_{1}\in\Delta^1_{s,u}}A_0\cdot {\rm d}X_{t_1} \cdot A_1 \cdot {\rm d}X_{t_3} \cdot A_2 \cdot {\rm d}X_{t_2} \cdot A_3.
\end{gather*}

The term on the second line above can not be expressed by composing order two full contraction operators. Instead, we can obtain it by \emph{composing} the operator,
\begin{equation*}
	(A_0,\dots,A_4)\mapsto \int_{\Delta^2_{s,u}}A_0\cdot{\rm d} X_{t_1} \cdot A_1 \otimes A_2 \cdot {\rm d} X_{t_2} \cdot A_3 \in \mathcal{A}\otimes\mathcal{A}
\end{equation*}
with the following full contraction one
\begin{equation*}
(A_0,A_1)\mapsto\int_{\Delta^1_{u,t}} A_0\cdot {\rm d} X_{t_1} \cdot A_1.
\end{equation*}
Thus a naive approach leads in fact to relations involving not only full contraction operators but also partial contractions. A remark on the terminology: we employ the term ``contraction''
to indicate that the operators reduce the degree of an input tensor, and ``full'' to indicate that it does
so maximally.

The main results of the paper are contained in the last section, Definition~\ref{def:endofacescontraction} and Theorem~\ref{cor:algebrastar}.
In this definition, we introduce the \emph{noncommutative signature of a smooth path} and in our main Theorem~\ref{cor:algebrastar}, we prove that, as for the classical theory, it yields a trajectory in a~certain group.

 \begin{Theorem} \label{thm:main_1}
 There exists a group $(\cG, \circ)$ such that for each algebra-valued smooth trajectory $X\colon [0,1]\to \cA$ there exists a map $\mathcal{X}\colon \Delta^{(2)}\to \cG$ from the two-dimensional simplex $\Delta^{2}$ with the following properties:
 	\begin{itemize}\itemsep=0pt
 		\item For any triple $s < u < t \in [0,1]^3$ one has
 		 \begin{equation}\label{eqn:nc_Chen_rel}
 			 \mathcal{X}_{s,t}=\mathcal{X}_{u,t}\circ \mathcal{X}_{s,u}.
 		 \end{equation}
 		\item
 		For any pair $ 0 < s < t < 1$, $\mathcal{X}_{s,t}$ \emph{has a set of coordinates} $\{\mathcal{X}_{s,t}(f),\,f \in \mathcal{F}\}$ where the set~$\cF$ contains all permutations $\gs$ and $\mathcal{X}_{s,t}(f)$ is a certain bounded operator acting on folded projective tensor products of $\mathcal{A}$ which coincides with \eqref{eqn:iteratedintegrals} when $f=\gs$.
 		\item Given two elements $\mathcal{X},\mathcal{Y} \in \mathcal{G}$,
 		 \begin{equation*}
 			 \mathcal{X} = \mathcal{Y} \Leftrightarrow \mathcal{X}(\sigma) = \mathcal{Y}(\sigma) \quad \text{for all permutations}~\sigma.
 		 \end{equation*}
 	\end{itemize}
 \end{Theorem}

We call the element $\mathcal{X}_{s,t}$ the \emph{noncommutative signature} of the path $X$ and the relations~\eqref{eqn:nc_Chen_rel} \emph{noncommutative Chen's relations}.

 \begin{Remark} We will define $\mathcal{G}$ as a set of representations of a certain algebra supported by trees with decorated leaves. The result that we want to prove in this work is purely algebraic and does not state any analytical property of~$\mathbb{X}$, which could be expected from the knowledgeable reader. We will in a separate work address integration theory against an irregular path drawn in~$\mathcal{A}$, and will gather at this time the relevant analytical context.
 \end{Remark}

\subsection{Outline}
Besides the introduction, this article is divided in two additional sections. In Section~\ref{sec:tensorforest}, we introduce a Hopf monoid of levelled forests, reminiscent of the Malvenuto--Reutenauer Hopf algebra of permutations.

In Section \ref{sec:contractionoperators}, we define the partial and full contraction operators we alluded to, see Definitions~\ref{def:contractionoperators} and~\ref{def:pcontractionoperators}.
In Section \ref{sec:chenrelation}, we prove a Chen relation for these operators, see Proposition~\ref{prop:chenrelation}. Next, we explain how this yields a path on a group of triangular algebra morphisms on an algebra spanned by couples of a tree and a word.
In Section~\ref{sec:geometricproperties}, we associate to the full and partial contractions operators a path of representations on the Hopf monoid of levelled forests we introduced in Section~\ref{sec:tensorforest}, see Theorem~\ref{thm:pathofrep}.

In Section \ref{sec:facescontraction}, we adopt a slightly different point of view and let the iterated integrals of a path acting on a~set of operators we call \emph{face-contractions}, see Definition \ref{def:facescontractions}. This yields a~certain triangular algebra morphism, see Definition~\ref{def:endofacescontraction} that we relate to the one introduced in Section~\ref{sec:chenrelation}. In Proposition~\ref{prop:product}, we relate partial- to full contraction operators.

In a forthcoming article, we continue to develop the theory. In particular, we introduce geometric noncommutative rough paths, geometric noncommutative controlled rough paths, and the operations of integration and composition.

\subsection{Notations}

In the following we denote by $\cA$ a generic complex $C^\star$ algebra with product $\mu$, unity $\1$, norm $\norm{\cdot}$ and involution $\star$. By definition, ($\mathcal{A},\|\cdot\|)$ is a Banach algebra, the multiplication $\mu$ and the involution $\star$ are continuous with respect to $\|\cdot\|$, and
\begin{equation*}
\|aa^{\star}\| = \|a\|^{2},\qquad a \in \mathcal{A}.
\end{equation*}

In order to deal with a topology on the algebraic tensor product $\otimes$ which behaves correctly with $\mu$, we will use the projective tensor product (see, e.g., \cite{Ryan02}). Given two Banach spaces $(E,\|\cdot\|_{E})$ and $(F,\|\cdot\|_{F})$, the projective norm of an element $x \in E\otimes F$ is defined by
\begin{equation*}
	\| x \|_{\vee} = \inf \Big\{\sum_{i} \|a_{i}\|_{E} \|b_{i} \|_{F} \colon x = \sum_{i} a_{i} \otimes b_{i}\Big\}.
\end{equation*}
We denote by $E\check{\otimes}F$ the completion of $E\otimes F$ for the projective norm. One can check the following properties
\begin{equation*}
	\|a\otimes b\|_{\vee} = \|a\|_{E}\|b\|_{F}, \qquad \|a_{\sigma(1)} \otimes\cdots\otimes a_{\sigma(n)}\|_{\vee} = \|a_{1}\otimes\cdots\otimes a_{n}\|_{\vee},
\end{equation*}
 for any permutation $\gs$ on the set $\{1,2, \dots, n\}$ and $a_1,\dots, a_n \in E$. The definition of projective norm yields immediately that the multiplication $\mu$ extends to a continuous map $\mathcal{A}\check{\otimes} \mathcal{A}\to \cA$ and, more generally, for any given pair of $C^{\star}$ algebras $\mathcal{A}$, $\mathcal{B}$, $\mathcal{A} \check{\otimes} \mathcal{B}$ is again a $C^{\star}$ algebra. From a~broader perspective, the projective tensor product makes the category of complex $C^{\star}$ algebras a symmetric monoidal category (see Appendix~\ref{appendix}). In order to lighten the notation, we will adopt the symbol $\otimes$ to denote both the projective tensor product between $C^{\star}$ algebras and the algebraic tensor product for pure tensors. Similarly, we will replace the product $\mu$ with a dot $\cdot$.

For $n\geq 1$ an integer, we denote by $\mathfrak{S}_{n}$ the set of permutations of $[n]:= \{1, \dots, n\}$. We use one-line notation for permutations,
writing $\sigma=(\sigma_1, \sigma_2, \dots, \sigma_n)$, where $\sigma_i:=\sigma(i)$. The neutral element of $\mathfrak{S}_{n}$ is also denoted by ${\rm id}_n$. Sometimes we may omit the commas and just write $\sigma=\sigma_1\sigma_2 \cdots \sigma_n$. By abuse of notation, the only permutation of $[0]:=\varnothing$ is denoted by $\mathsf{\varnothing}$, from which we defines $\mathfrak{S}_0:=\{\varnothing\}$. Given two integers $a,b$, we denote by ${\rm Sh}(a,b)$ the set of all shuffles of the two intervals $\llbracket 1,a \rrbracket$ and $\llbracket a+1,a+b\rrbracket$, that is $\sigma \in {\rm Sh}(a,b)$ if and only if $\sigma$ is non-decreasing on $\llbracket 1,a \rrbracket$ and on $\llbracket a+1,a+b \rrbracket$.

\section{Algebraic structure on levelled forests}\label{sec:tensorforest}

The objective of the present section is to introduce the main combinatorial tool that will be used in this work: the levelled trees and forests. We will review their main properties and introduce new algebraic structures to them.

\subsection{Levelled trees and forests}\label{sec:levelledforests}

In the literature, one broadly finds several equivalent representations of a permutation, such as a bijection of a finite set or a finite word without repetitions on positive integers. We will mainly use the last one and a~third~-- tree-like~-- graphical representation, presented in different variants in the literature such as \cite[pp.~23--24]{S1999}, \cite[Definition~9.9]{BW1997} or \cite[p.~478]{AS2006}. We will follow the versions used by Loday and Ronco in \cite[Section~2.4]{loday1998hopf} and Forcey, Lauve and Sottile in \cite[Section~2.2.1]{FLS2010}.

First, recall that a \emph{planar rooted tree} is a planar graph with no cycles and one distinguished vertex which we call the root. We oriented every tree from bottom to top: the target of an edge is the vertex further to the root. In this orientation, each vertex of a tree has at most one incoming edge (the root is the only vertex with no incoming edge) and at most two outcoming edges.

A \emph{leaf} of a tree is a vertex with no outcoming edges. The \emph{degree} of a tree is the number of its leaves, we denote it by $\vert\tau\vert$ if $\tau$ is a planar tree. An \emph{internal vertex} of a tree is a vertex that is not a leaf. The set of internal vertices of a tree $\tau$ is denoted by $\mathbb{V}(\tau)$ and we set $\|\tau\|:=|\mathbb{V}(\tau)|$.
The set $\mathbb{V}(\tau)$ of internal vertices of a planar tree $\tau$ is equipped with a partial order $\leq_{\tau}$: if~$u$,~$v$ are two vertices of $\tau$, we write $u \leq_{\tau} v$ if there is an oriented path of edges of $t$, moving away from the root, from $u$ to $v$. The poset $(\mathbb{V}(\tau), \leq_{\tau})$ has one minimum (the root of $\tau$) and several maxima (the leaves of $\tau$).

\begin{figure}[!ht] \centering
	\includegraphics{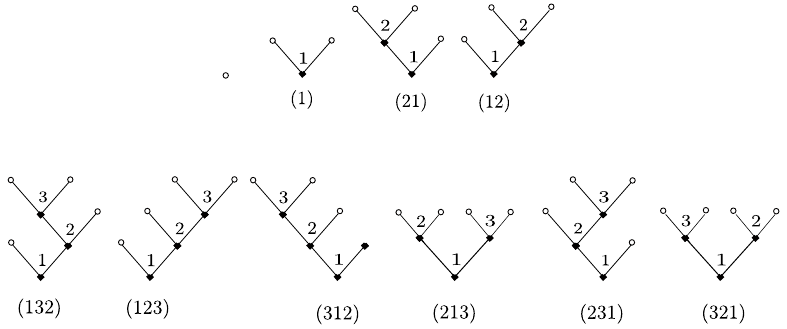}
	\caption{Examples of levelled trees $(\tau,\lambda)$ in $\lpbt$ and their associated permutations in $\mathfrak{S}_1$, $ \mathfrak{S}_2$ and $\mathfrak{S}_3$.}\label{fig:2}
\end{figure}

A \emph{planar binary tree} is a planar rooted tree for which every internal node has two children. A~\emph{levelled binary tree} (or simply \emph{levelled tree}) is a binary tree $\tau$ together with a~linear extension of the poset $(\mathbb{V}(\tau), \leq_{\tau})$. Levelled trees are also called \emph{ordered binary trees} (see~\cite{AS2006}). By definition, a levelled tree with degree one is the root tree (see Figure~\ref{fig:2}). Also, notice that the root tree has no internal vertices and corresponds to levelled tree $(\sbt,\varnothing)$ where $\varnothing$ denotes the unique function from the empty set to the empty set.

We denote by
$\lpbt(n)$ the set of levelled trees with $n$ leaves, and $\lpbt:=\cup_{n \geq 1}\,\lpbt(n)$. The complex span of $\lpbt$ is a graded vector space, and its homogeneous component of degree $n \geq 1$ is the linear span of $\lpbt(n)$.

 We justify now the terminology for levelled trees. Following \cite[p.~7]{PS2008}, a~\emph{level function} on a~tree~$t$ is a surjective increasing map
\begin{equation*}
\lambda\colon \ (\mathbb{V}(\tau), \leq_{\tau} )\to A,
\end{equation*}
where $A$ is a totally ordered set. If $\tau$ is a planar binary tree and $A=[\|\tau\|]$, then the pair $(\tau,\lambda)$ corresponds precisely to a levelled tree. If $v$ is an internal vertex of $t$, we say that $v$ has level~$\lambda(v)$.

The following result seems to be folklore. For proof of this result, see \cite[Proposition~2.3]{loday1998hopf}.
\begin{Proposition}[\cite{loday1998hopf}]\label{PermTree}
	For every integer $n\geq 0$, the set of levelled trees with $n+1$ leaves is in bijection with the set of permutations $\mathfrak{S}_{n}$.
\end{Proposition}

The bijection associates to any levelled tree $(\tau,\lambda)$ with $n+1$ leaves a permutation $\sigma=\sigma(\tau, \lambda) \in \mathfrak{S}_n$ as
follows. Label the leaves of $\tau$ with $0,1,2, \dots, n$ (in this order), from left to right. For each $1 \leq i \leq n$, let $v_i$ be the vertex which lies in between the leaves $i-1$ and $i$. Then $\sigma:=\sigma_1 \sigma_2 \cdots \sigma_n$, with $\sigma_i:=\lambda(v_i)$, see Figure~\ref{fig:2}.

\begin{figure}[!ht]	\centering
	\includegraphics{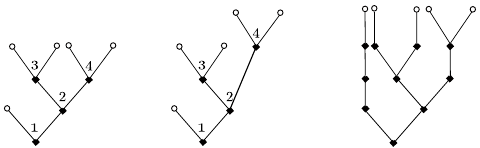}
	\caption{Example of a levelled tree and its associated sparse quasi-binary tree, before adding straight edges in the centre and after, on the right.}\label{fig:3}
\end{figure}

When illustrating a levelled tree $(\tau, \lambda)$, it will be convenient to emphasize the levelling (the map $\lambda$) of a tree without the use of labels on the vertices. To do so, we \emph{position} each vertex $v$ of the tree $\tau$ at the level $\lambda(v)$; it is represented by a dot with $y$-coordinate $\lambda(v)$. We add straight edges to $\tau$ (see Figure~\ref{fig:3}) so that the level $\lambda(v)$ is populated with $|\lambda(v)|$ vertices (in particular, on the first level, we find the root of the tree). In the resulting tree, all vertices placed on the same level $n$ have an equal distance to the root.

Notice that the resulting unlabeled tree is such that every vertex has either one or two children, and there is a unique vertex with two children.

We call such a tree a \emph{sparse quasi-binary tree}.
All operations introduced in this section have a convenient pictorial description using sparse quasi-binary trees.

In summary, we have three equivalent ways to identify the same object:
\begin{equation*}
\begin{array}{@{}c}
 \text{Permutations} \\
 \sigma\colon [n]\to [n]
\end{array} \quad
\longleftrightarrow
\quad
\begin{array}{c}
 \text{Levelled trees} \\
 \text{with } n+1 \text{ leaves}
\end{array}\quad
\longleftrightarrow
\quad
\begin{array}{c}
 \text{Sparse quasi-binary trees} \\
 \text{with } n+1 \text{ generations.}
\end{array}
\end{equation*}

We use the symbol $t$ (with super and subscript) for a levelled tree presented either as a~permutation $\sigma$ or as a~pair $(\tau,\lambda)$. The representation of $t$ as a sparse quasi-binary tree will only be used in drawings.

Levelled trees are not sufficient for our purposes. We will extend now the notion of levelled trees to \emph{forests}. A \emph{planar forest} is a word (a noncommutative monomial) on planar trees.

In the following, we denote by ${\sf nt}(\varphi)$ the number of trees in the forest $\varphi$, $|\varphi|$ the total number of leaves in the forest and we set $\|\varphi\|$ equal to the number of internal vertices of the forests. If all trees of $\varphi$ are binary trees, then $\|\varphi\|=|\varphi|-{\sf nt}(\varphi)$.
The poset $(\mathbb{V}(\varphi),\leq_{\varphi})$ of ordered vertices of $f$ is the union of the posets of vertices of the trees in $f$.

In the following, we will just consider planar forests of binary trees. The notion of level function for a tree is naturally extended to any forest. This allows considering the following analogue of levelled binary trees to binary forests.

\begin{Definition}[{levelled planar binary forests} $\lpbf$]
	A \emph{levelled planar binary forest} $f$ (or simply a \emph{levelled forest}) is a pair $(\varphi, \lambda)$ formed by a binary forest $\varphi$ and an increasing bijection
\begin{equation*}
\lambda \colon \ (\mathbb{V}(\varphi), \leq_\varphi)\to[||\varphi||].
\end{equation*}
	We denote the set of planar binary forests by $\lpbf$.
\end{Definition}

\begin{figure}[!ht]	\centering
	\includegraphics{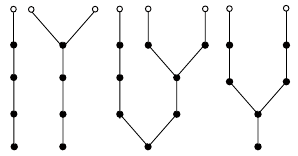}
	\caption{A sparse quasi-binary forest with four trees.}\label{fig:levelledforest}
\end{figure}

The \emph{degree} of a levelled planar forest $(\varphi, \lambda) \in \lpbf$ is the number of leaves of $\varphi$ and is denoted by $|\varphi|$. If $n \geq 1$ and $m\geq 1$, we denote by $\lpbf(m,n)$ the set of levelled forests with $n$ leaves and $m$ trees. This allows defining a bigraduation on the set $\lpbf$.

A generation of a levelled tree is a set of internal vertices on the same level, that is at the same distance from the root in the sparse quasi-binary tree representation (we thus take into account the labelling of the vertices). This notion extends to any forest.
A levelled forest $(\varphi, \lambda)$ can be pictured as a forest of quasi-binary trees, each with equal number of generations, in the same way as explained before for levelled trees, where the internal vertices are ordered vertically by adding straight edges according to $\lambda$, see Figure~\ref{fig:levelledforest}. In this representation, there is a unique internal vertex with two children among all vertices of the forest of the same generation. We call such a forest a \emph{sparse quasi-binary forest}.

We now introduce several (classical) operations on levelled trees and levelled forests.

If $n \in \mathbb{N}$, \emph{weak composition} of $n$ is a sequence $c=(c_1, c_2, \dots, c_k)$ of non-negative integers with $n=c_1 + \cdots + c_k$. If we restrict each $c_i$ to be non-zero, then $c$ is called a \emph{composition} of $n$. To a weak composition $c$, we associate the multiset $I(c):=\{c_1, c_1+c_2, \dots, c_1+c_2+\cdots+ c_{k-1}\}$. This gives a bijection between weak compositions of $n$ and multisubsets of $\{0\}\cup [n]$. If~$c$ is a~composition, $I(c)$ is a set and we obtain a bijection between compositions of $n$ and subsets of~$[n-1]$. Weak compositions of $n$ are partially ordered by \emph{refinement}. The cover relations are of the form
\begin{equation*}
(c_1, \dots, c_i+c_{i+1}, \dots, c_k) \lessdot (c_1, \dots, c_i, c_{i+1}, \dots, c_k).
\end{equation*}
We write $c \vDash n$ if $c$ is a composition of $n$, and $c \vDash_0 n$ if $c$ is a weak composition of $n$.
The length~$k$ of a (weak) composition $(c_1,\dots,c_k)$ is the number of parts of the composition.

The bijection between levelled trees and permutations used in the proposition \eqref{PermTree} extends to words $w=w_1w_2\cdots w_r \in \mathbb{N}^*$ without repetition of letters; the associated levelled tree $t$ has now a level function $\lambda\colon \left(\mathbb{V}(\tau), \leq_{t}\right)\to \{w_1, w_2, \dots, w_r\}$.

Every levelled forest $(f, \lambda) \in \lpbf(n+k,k)$ gives rise to a pair $(\sigma, c)$, where $\sigma \in \mathfrak{S}_{n}$ is obtained by concatenating the non-empty words corresponding to each tree in $f$ (from left to right) under the above-described bijection, and $c \vDash_0 n$ is the weak composition of length $k$ obtained by tracking the number of internal vertices of each tree in the forest $f$. Reciprocally any pair $(\sigma,c)$ yields a levelled planar binary forest, using the bijection between non-repeating words and levelled planar trees. We call \emph{split permutation} a pair $(\sigma,c)$ with $\sigma \in \mathfrak{S}_n$ and $c \vDash_0 n$:
\begin{equation*}
\begin{array}{@{}c}
 \text{Split permutations} \\
 \sigma:[n]\to [n]\\
 c \vDash_0 n
\end{array} \quad
\longleftrightarrow
\quad
\begin{array}{c}
 \text{Levelled forests} \\
 \text{with } n+1 \text{ leaves}
\end{array}\quad
\longleftrightarrow
\quad
\begin{array}{c}
 \text{Sparse quasi-binary forests} \\
 \text{with } n+1 \text{ generations.}
\end{array}
\end{equation*}

As for levelled trees, we use the symbol $f$ to denote a levelled forest presented either as a~pair ($\varphi,\lambda)$ or as a split permutation $(\sigma,c)$. The presentation of $f$ as a sparse quasi-binary forest will be used in the drawings only.

\begin{figure}[!ht]	\centering
 \includegraphics{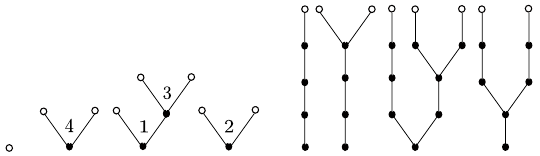}
 \caption{A split permutation $((4132),(0,1,2,1))$ drawn as a levelled forest and as a sparse quasi-binary forest.} \label{Figure04}
 \end{figure}

The next definition introduces the notion of \emph{vertical splitting} for levelled forests.
 Informally, a vertical splitting of a levelled forest $f$ consists in breaking $f$ into two forests, each bordered by a chosen path of edges in $f$, starting at a leaf of $f$ and ending at a root of a tee in $f$. The first forest (resp. the second) is on the right (resp. one the left) of this path.

\begin{Definition}[vertical splitting of levelled forests] Pick a levelled forest~$f$ presented as a~split permutation $(\sigma,c) \in \mathfrak{S}_n$, $c\vDash_0 n$ and $ d \vDash_0 n$ such that $c \lessdot d$, we define the vertical splitting $\curlyvee_{d}((\sigma,c))$ of $(\sigma,c)$ following $d$ by $(\sigma,d)$.
\end{Definition}

We introduce inverse operations to splitting. The first one takes every tree in the representation of a levelled forest as a sparse quasi-binary forest and glues all together the trees of that forest along their external paths of edges. In terms of split permutation, this operation corresponds to the projection,
\begin{equation*}
\tau^{\flat} := \sigma,\qquad f=(\sigma,c) \in \lpbf{}.
\end{equation*}

We will also need a local operation gluing two consecutive trees in the representation of a levelled tree as a sparse quasi-binary tree, once again those operations are most effectively written in terms of split permutation. We set for any $f=(\sigma,c)\in \lpbf$ and $1 \leq i \leq {\sf nt}(f)-1$
\newcommand{\glu}[1]{{|\!|_{\raisebox{-1mm}{\scalebox{0.8}{$#1$}}}}}
\begin{equation*}
\glu{i} f =\glu{i}(\sigma,c) :=(\sigma, (c_1,\dots,c_i+c_{i+1},\dots,c_{{\sf nt} (f)})).
\end{equation*}

\begin{figure}[!ht]\centering
\includegraphics{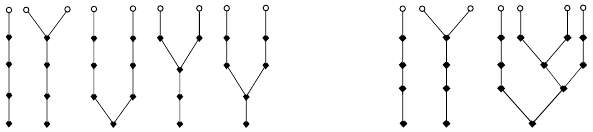}
\caption{The forests $\curlyvee_{(0,1,1,1,1)}((4132),(0,1,2,1))) $ and $\glu{3}((4132),(0,1,2,1))$.}\label{fig:my_label}
\end{figure}

We consider \emph{horizontal analogues} to the above operations of vertical splitting and gluing.
For any word $w=w_1w_2\cdots w_r $ where each letter takes value in $\mathbb{N}^*$ and $A \subseteq \text{alph}(w)$, the alphabet generated by the different letters contained in $w$, we define $w \cap A$ as the word obtained from $w$ by erasing the letters which are not in $A$. We write $w' \subseteq w$ if there exists $A \subseteq \text{alph}(w)$ such that $w'=w\cap A$. In this case, we say that $w'$ is a \emph{subword} of $w$. We use now the definition of subword to define the notion of subtree and subforest. Let $t \in \lpbt(n)$ a levelled binary tree, represented as a permutation $\sigma$. A \emph{levelled subtree} (or just \emph{subtree}) of~$t$ is a levelled binary tree~$t^{\prime}$ with associated permutation $\sigma^{\prime}$ of the form $\sigma^{\prime}=\sigma \cap [p]$, for $0 \leq p \leq n$.

In this case, we write $\tau^{\prime} \subseteq \tau$. In terms of sparse quasi-binary trees, $\tau^{\prime}$ is a subtree of $\tau$ if
there exists $0 \leq p \leq \|\tau\|$ such that $\tau^{\prime}$, seen as a quasi-binary tree, coincides with the quasi-binary tree associated with $\tau$ by erasing all vertices on generations strictly bigger than $p$.

\begin{figure}[!ht]\centering
\includegraphics{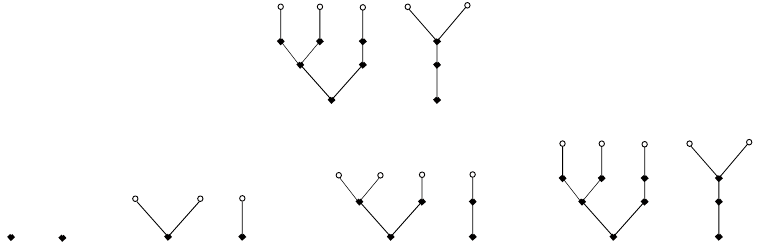}	
\caption{On the first line: levelled forest $f=((213),(2,1))$, drawn as a quasi-binary forest. On the second line: all levelled subforests included in $f$.}\label{fig:exsubforest}
\end{figure}

This notion extends to levelled forests. Pick a levelled forest $f = (\sigma,(c_1,\dots,c_k)) \in \lpbf{}$ and denote by $(\sigma_1,\dots, \sigma_k)$ the restrictions of $\sigma$ to the parts of $c$:
\begin{equation*}
\sigma_i = \sigma_{]c_1+\cdots+c_{i-1},c_1+\cdots+c_i]},\qquad 1 \leq i \leq k,
\end{equation*}
with the convention that $\sigma_{\varnothing} = \varnothing$.
A \emph{levelled forest} $f^{\prime}=(\sigma^{\prime},(c^{\prime}_1,\dots,c^{\prime}_k))$ with same number of trees as $f$ (the composition $c^{\prime}$ has the same number of parts as $c$) is a subforest of $f$ is there exists $p \leq \|f\|$
such that
\begin{equation*}
\sigma^{\prime}_i = \sigma_i \cap [p].
\end{equation*}
We write in this case $f^{\prime} \subset f$. In terms of levelled forest, writing $(f, \lambda)$ and $(f^{\prime}, \lambda^{\prime})$, then we have $(f^{\prime}, \lambda^{\prime}) \subseteq (f, \lambda)$ if $f_\alpha$, $f_\beta$ have the same number of trees and
the labelling $\lambda_{\alpha}$ restricts to $\lambda_{\beta}$ on the internal vertices of $f_\beta$. Equivalently, considering $\alpha$ and $\beta$ as quasi-binary forests, $\beta \subseteq \alpha$ if and only if $\beta$ comprise all nodes of $\alpha$ up to a certain generation of the quasi-binary forest $\alpha$.

Consider a word $w=w_1w_2\cdots w_r $ where each letter takes value in $\mathbb{N}^*$ and $A \subseteq \text{alph}(w)$.
Let $w(A) = \{i^w_A(1) < \cdots < i^w_A(q)\} $ be the increasing sequence of indices of the letters of $w$ in $A$ and define the composition ${\sf wc}(w,A)\vDash_0 r$
\begin{equation*}
{\sf wc}(w,A):=\big(|]1,i^w_A(1)]|, \dots, |]i^w_A(i), i^w_A(i+1)]|, |]i^w_A(q),r]|\big).
\end{equation*}

For example, if $r=10$ and $w=1293548851$, let $A=\{2,3,4,7,9\}$. Then $[9]\setminus A=\{1,5,6,8\}$, $w(A)=\{2 < 3 < 4 < 6\}$. Therefore, $\textsf{wc}(w,A)=(1,1,1,2,4)$. Thanks to this notion, we introduce the horizontal splitting of a levelled tree.
\begin{Definition}[horizontal splitting of trees]
Pick $\tau$ a levelled tree seen as a permutation $\sigma$ and an integer $ 0 \leq p\in \| \sigma \|$. The \emph{horizontal splitting} of $\sigma$ at $p$ is the couple
\begin{equation*}
\succ_p\!(t):=\big(\sigma \cap [p], (\sigma \cap \{p+1, \dots, n\}, \textsf{wc}(\sigma,[p]))\big),
\end{equation*}
formed by the subword of $\sigma$ containing the letters in $[p]$, and the pair $(\sigma \cap \{p+1, \dots, n\}, \textsf{wc}(\sigma,[p]))$ representing the sequence of subwords of $\sigma$ obtained from $\sigma$ after erasing the letters in $[p]$. We call the first component of $\succ_p$ the \emph{lower component} of the cut and the second component the \emph{upper component}.
\end{Definition}

For instance, $\succ_2(25143)=(21; (543),(0,1,0,2))$. Horizontal splitting acts on the sparse quasi-binary tree representation by detaching the first lower $p$ generations (we include all edges connected to the vertices of the $p^{\rm th}$ generation. The resulting levelled tree forms the lower component of the cut and the generations above it yield the upper component of the cut.

\begin{figure}[!ht] \centering
\includegraphics{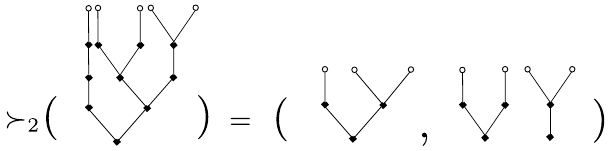}
\caption{Example of a horizontal splitting of a sparse quasi-binary tree.}\label{fig:hor_split}
\end{figure}

This operation extends to levelled forests in a straightforward manner.
\begin{Definition}[horizontal splitting of levelled forests]
Let $f=(\sigma,c)$ be a forest in $\lpbf(k,n)$, with $\sigma \in \mathfrak{S}_{n-1}$ and $(c_1,\dots,c_k) \vDash_0 n-1$ a weak composition of length $k$. Once again, denote by $\sigma_i$, $1 \leq i \leq \|k\|$ the restriction of $\sigma$ to the parts of $c$. The \emph{horizontal splitting of the levelled forest} $f$ at the level $p$ is the forest
\begin{gather*}
 {\succ}_p(f)= ({\succ}_p(f)_-,{\succ}_p(f)_+) := ({\succ}_p(\sigma_1)_-\cdots{\succ}_p(\sigma_k)_-,{\succ}_p(\sigma_1)_+\cdots {\succ}_p(\sigma_k)_+).
\end{gather*}
 \end{Definition}
For convenience, we use the shorter notations
\begin{equation*}
f_-^p := \succ_p(f)^-,\qquad f_+^p:= \succ_p(f)^+.
\end{equation*}
Consider for example the forest $\curlyvee_{\! (1,1,3,0,1)}(542631)=(5,4,263,1)$. Since $\succ_3\!(5)=(\varnothing; (5))$, $\succ_{3}\!(4)=(\varnothing; (4))$, $\succ_3\!(263)=(23; (\varnothing, 6, \varnothing))$, $\succ_3\!(\varnothing)=(\varnothing; (\varnothing))$ and $\succ_3\!(1)=(1; (\varnothing))$, we obtain
\begin{equation*}
{\succ}_3(5,4,263,1)=((\varnothing, \varnothing, 23, \varnothing, 1);(5, 4, \varnothing, 6, \varnothing, \varnothing, \varnothing)).
\end{equation*}

A horizontal cut of a levelled forest $f$ outputs a pair of forests, the lower component being a~subforest of $f$. And reciprocally, one can check that a subforest $f^{\prime} \subset f$ yields a unique horizontal cut of~$f$, the cutting level $k$ being the number of internal nodes of $f^{\prime}$ (equivalently the number of generations). We choose then to index horizontal cuts of a levelled forest by its set of subforests augmented with the root tree.

\begin{Definition}\label{def:backslash_forest}
For any given levelled forest $f$ and subforest $f^{\prime} \subset f$ we denote by $f\backslash f^{\prime}$ the upper component of the horizontal cut induced by $f'$. We say that the levelled forest $f$ is \emph{compatible} with the levelled forest $f^{\prime}$ if the number of trees of $f$ matches the number of leaves of $f^{\prime}$ in the representation of $f$ and $f^{\prime}$ as sparse quasi-binary trees or as levelled trees.
\end{Definition}

By definition of horizontal splitting one has immediately the identity
$\succ_{||f^{\prime}||} f = (f^{\prime},f\backslash f^{\prime}). $ Moreover, if both levelled forests are given as split permutations, $f=(\sigma,(c_1,\dots,c_k))$, $\sigma \in \mathfrak{S}(n)$ and $f^{\prime}=(\sigma^{\prime},(c^{\prime}_1,\dots,c^{\prime}_q))$, $\sigma^{\prime} \in \mathfrak{S}(n^{\prime})$ compatibility means that the number $k$ of parts of $c_1$ is equal to $n^{\prime}+q = {\sf nt}(f^{\prime})+\|f^{\prime}\|$.

Given two permutations $\sigma = \sigma_1\cdots\sigma_n \in\mathfrak{S}_{n}$ and $\sigma^{\prime} = \sigma^{\prime}_1\cdots\sigma^{\prime}_{n^{\prime}}\in\mathfrak{S}_{n^{\prime}}$ we denote by $\sigma\otimes\sigma^{\prime}$ the shifted concatenation of $\sigma$ and $\sigma^{\prime}$,
\begin{equation*}
\sigma \otimes \sigma^{\prime} = \sigma_1\cdots\sigma_{n}(\sigma(1)+n)\cdots (\sigma(n)+n^{\prime})
\end{equation*}
The following operation, inverse to horizontal splitting, is better understood in terms of levelled forests.
\begin{Definition}[horizontal gluing]
Pick two compatible levelled forests $f=(\tau,\lambda) \in \lpbf{}$ and $f^{\prime}=(\tau^{\prime},\lambda^{\prime})\in \lpbf{}$, we define the horizontal gluing the levelled forests $f \# f^{\prime} = (\phi\# \phi^{\prime}, \lambda^{\prime\prime})$ where
\begin{enumerate}\itemsep=0pt
\item $\phi \# \phi^{\prime}$ is the planar forest obtained by successively superimposing a root of a tree in $\phi^{\prime}$ to a leaf of $\phi$.
\item The labelling $\lambda^{\prime \prime}$ restricts to $\lambda$ on the internal vertices of $\phi$ in $\phi\# \phi^{\prime}$ and to the labelling $\lambda^{\prime}$ translated by $\|f\|$ on the internal vertices of $\phi^{\prime}$ in $\phi\# \phi^{\prime}$.
\end{enumerate}
\end{Definition}

In terms of sparse quasi-binary trees and forests, horizontal gluing corresponds to stacking the sparse quasi-binary trees representing $f^{\prime}$ above the one representing $f$. Writing this operation in the representation of levelled forests as split permutations is cumbersome and is left to the reader, see also the figure below.
\begin{figure}[!ht] \centering
\includegraphics{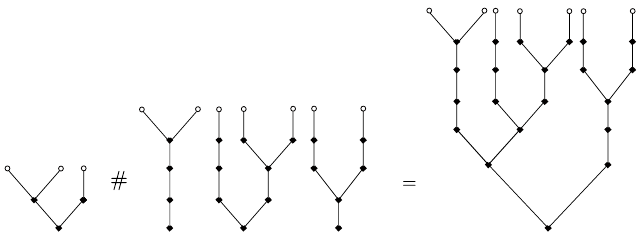}
\caption{Horizontal gluing of a tree and a forest.} \label{fig:horizontal_gluing}
\end{figure}

\subsection{A first algebraic structure on levelled forests}
We briefly recall important Hopf algebraic structures on permutations.
The vector space $\bigoplus_{n \geq 0}\mathbb{C}\mathfrak{S}_n$ has a graded Hopf algebra structure, called $FQ\textsf{Sym}$, introduced in \cite{MR1995}, and also considered in \cite{GKLLRT1995}, where it is called the algebra of free quasi-symmetric functions. This Hopf algebra is noncommutative, non-cocommutative, graded and self-dual. In the following, we introduce the product $\shuffle$ dual to the aforementioned coproduct of the Hopf algebra of noncommutative symmetric functions.
First, for two permutations $\alpha \in \mathfrak{S}_n$ and $\beta\in \mathfrak{S}_p$, we introduce
\begin{equation*}
 \textsf{Sh}(\alpha,\beta):=\{\sigma \in \mathfrak{S}_{n+p} \colon \st(\sigma \cap \{1, \dots, n\})=\alpha,\, \st(\sigma \cap \{n+1, \dots, n+p\})= \beta\},
\end{equation*}
the set of \emph{shuffles} of $\alpha$ and $\beta$. In the above equation ${\rm std}$ denotes the standardization map:
${\rm std}(\sigma \cap \{a_1 < \cdots <a_k\})$ is the image of $\sigma \cap \{a_1 < \cdots < a_k\}$ by the map $a_k \mapsto k$.

In particular, if ${\rm id}_n$ and ${\rm id}_p$ are the identity permutations of the symmetric groups~$\mathfrak{S}_n$ and~$\mathfrak{S}_p$, respectively, we let
\begin{equation*}
 \text{Sh}(n,p):=\sf{Sh}({\rm id}_n, {\rm id}_p).
\end{equation*}
Given a permutation $\tau \in \mathfrak{S}(n)$ and integer $p\geq 1$ we denote by $\tau + p$ the bijection of $\{p+1,\allowbreak \dots,p+n\}$ defined by
\begin{equation*}
(\sigma+p)(p+k) = \sigma(k)+p.
\end{equation*}
With $\sigma \in\mathfrak{S}(p)$, $\sigma\otimes(\tau +p)\in\mathfrak{S}(n+p)$ is the permutation equal to $\sigma$ on $\{1,\dots,p\}$ and $(\tau +p)$ on $\{p+1,\dots,p+n\}$,
\begin{equation*}
\Delta(G_{312})=\1^*\otimes G_{312}+ G_1 \otimes G_{21}+ G_{12}\otimes G_1 +G_{312}\otimes \1^*.
\end{equation*}

The \emph{shuffle product of two permutations} $\sigma \in \mathfrak{S}_p$ and $\tau\in\mathfrak{S}_n$ is defined by
\begin{equation*}
\sigma \shuffle \tau := \sum_{\alpha \in {\rm Sh}(p,n)} \alpha \circ [\sigma \otimes (\tau + p)],
\end{equation*}
where we have used the explicit notation $\circ$ for the composition of permutations. We extend the shuffle product from levelled trees (permutations) to levelled forests. If the permutations~$\sigma$ and~$\tau$ are presented as sparse quasi-binary trees, computing their shuffle product is done by adding straight edges at the bottom of $\tau$, then gluing this tree to the right of~$\sigma$ (we identify the outer paths of edges) and finally shuffling vertically the generations.
\begin{Definition}[shuffle product of levelled planar forests]\label{shuffle_product_forests}
	Let $f=(\sigma,(c_1,\dots,c_k))$ and $g=(\tau,(c^{\prime}_1,\dots,c^{\prime}_q))$ be two levelled forests, we define the shuffle product of $f$ and $g$ by
	\begin{equation*}
		f \shuffle g = (\sigma \shuffle \tau, (c_1,\dots,c_k+c_1^{\prime},\dots,c^{\prime}_q)).
	\end{equation*}
\end{Definition}
\begin{Example} We present in detail the product $12\shuffle (12,(1,1))$, presented as sparse quasi-binary forests in Figure~\ref{fig:shuffleproduct}.
\begin{figure}[!ht]\centering
\includegraphics{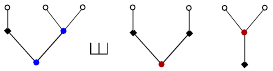}
\caption{Shuffle product of levelled forests.}\label{fig:shuffleproduct}
\end{figure}
We decorate the branching nodes with two different colors to track the two original permutations in the product. The result of the product is given by the following sum of sparse quasi-binary forests (see Figure~\ref{fig:computations}).
\begin{figure}[!ht]\centering
\includegraphics{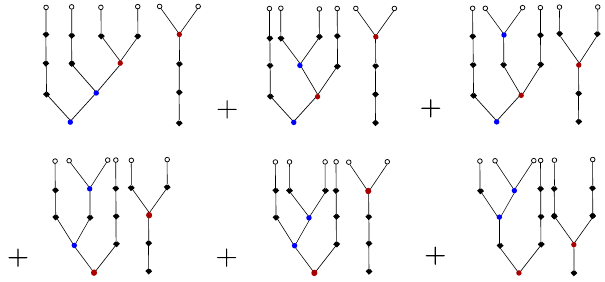}
\caption{The product in Figure~\ref{fig:shuffleproduct} expanded over sparse quasi-binary forests.}\label{fig:computations}
\end{figure}

We notice that by the construction of $\shuffle$ all terms contained in the sum the blue and red dots preserve the same ordering of generations in the initial factors.
\end{Example}

A permutation $\sigma \in \mathcal{S}_n$ acts on the left of a levelled forest $f=(\sigma,c) \in \lpbf$ with $n$ generations as follows
\begin{equation*}
\sigma \cdot (\tau,c) := (\sigma \circ \tau, c).
\end{equation*}

We denote by $c_{n}$ the permutation $(n,1)(n-1,2)\cdots \big(n-\lfloor\frac{n}{2}\rfloor,\lfloor \frac{n}{2}\rfloor\big)$. We use the right action of $c_{n}$ to define a involution on $\lpbf(n-1)$, which is the \emph{horizontal mirror symmetric of a forest}
\begin{align*}
	\theta \colon \ \mathbb{C}\lpbf & \to \mathbb{C}\lpbf,\\ f=(\sigma,c) & \mapsto \big(\sigma c_{\|f\|}^{-1}, c_n,\dots,c_1\big).
	\end{align*}

\begin{Proposition} $(\mathbb{C}\lpbf,\shuffle, \theta)$ is an involutive algebra.
\end{Proposition}
\begin{proof}
	The result follows as a direct a consequence of the following two facts: the left and right actions of $\mathfrak{S}_{n}$ on $\lpbf(n-1)$ commute and $c_{n+m} = \tau_{n,m} \circ (c_{n}\otimes (c_{m}+n))$ (we add $n$ to all letters of $c_{m}$ and concatenate the resulting word to the one representing $c_n$) where $\tau_{n,m}$ is the shuffle in ${\rm Sh}(n,m)$ determined by
	$\tau_{n,m}(1) = m+1$, $\tau_{n,m}(n)=m+n$.
\end{proof}

We will sometimes refer to $\theta$ as the \emph{horizontal involution}, for obvious reasons, to distinguish it from a second involution permuting vertically the generations of a levelled forest that we define below.

\subsection{Hopf monoid of levelled forests}\label{sec:aconvolutiongp}
\def\bclf{\mathcal{L}\mathcal{F}}
\def\bclt{\mathcal{L}\mathcal{T}}
In this section, we introduce a Hopf algebraic structure on the bicollection of spanned by levelled forests and denoted $\bclf$,
\begin{equation}\label{levelled_forests}
\bclf(n,m) = \mathbb{C} \lpbf(m,n),\qquad n,m\geq 1.
\end{equation}

In addition, we set $\bclf(0,0)=\mathbb{C}$, $\bclf(0,n)= \bclf(m,0)=0$, $n,m\geq 1$ and we denote by $\bclt$ the collection spanned by levelled binary trees
\begin{equation}\label{levelled_trees}
\bclt(n)= \mathbb{C}\lpbt(n), \qquad n\geq 1.
\end{equation}
This Hopf algebra is an object in the category of bicollections endowed with the vertical tensor product $\overt$. In general, as is briefly explained in the Appendix \ref{appendix}, the two-folded vertical tensor product $A\overt A$ of a monoid $A$ in the monoidal category (Coll$_{2}$, $\overt$) is not a monoid in the same category. Owing to the fact that the monoid generated by $\bclf$ in $(\mathrm{Coll}_{2},\overt)$ is \emph{symmetric}, in particular, $\bclf\overt \bclf$ is a monoid in a natural way, it makes sense to require compatibility between a product and a coproduct on $\bclf$. We write the unit $\mathbb{C}_{\overt}$ for the vertical tensor product as\looseness=-1
\begin{equation*}
	\mathbb{C}_{\overt} = \bigoplus_{n\geq 0} \mathbb{C}1_{n}.
\end{equation*}
Recall that we denote by $|f|$ the number of leaves of a levelled forest $f$ and ${\sf nt}(f)$ the number of trees in $f$.

We begin with the definition of the coproduct acting on the bicollection ${\sf \bclf}$ of levelled forests. Let $f$ be a levelled forest. Let $f^{\prime}$ be a levelled subforest of $f$ (recall that $f^{\prime}$ contains the roots of all trees in $f$). By definition of the forest $f \backslash f^{\prime}$, the number of outputs of the forest $f\backslash f^{\prime}$ is equal to the number of inputs of the forest $f^{\prime}$ (the number of trees of $f\backslash f^{\prime}$ matches the number of leaves of $f^{\prime}$), the following makes senses
\begin{equation}	\label{eqn:coproduct}
	\Delta(f) = \sum_{f^{\prime} \subset f} f^{\prime} \overt f\backslash f^{\prime},\qquad f \in \bclf.
\end{equation}
This operation is a genuine coproduct with respect to the vertical tensor product.
\begin{Proposition}\label{prop:coprodut}
	The morphism $\Delta\colon \bclf\rightarrow \bclf\overt\bclf$ is coassociative
	\begin{equation*}
		(\Delta \overt \mathrm{id}_{\bclf}) \circ \Delta = (\mathrm{id}_{\bclf} \overt \Delta) \circ \Delta
	\end{equation*}
	and the morphism $\varepsilon\colon \bclf \rightarrow \mathbb{C}_{\overt}$ given by
	\begin{equation*}
	\varepsilon(f)= \begin{cases}
 1_{n}& \text{if}\ f= \underbrace{\bullet \cdots \bullet}_{n \ \text{times}},\\
 0 & \text{otherwise}\end{cases}
	\end{equation*}
	is the counity for $\Delta $, i.e.,
	\begin{equation}\label{eqn:conunitpp}
		(\varepsilon \overt \mathrm{id}_{\bclf}) \circ \Delta = (\mathrm{id}_{\bclf} \overt \varepsilon) \circ \Delta = \mathrm{id}.
	\end{equation}
\end{Proposition}
\begin{proof}
	Let $f$ be a levelled forest, to show coassociativity we notice that
	\begin{equation*}
		((\Delta \overt \mathrm{id}_{\bclf}) \circ \Delta)(g) = \sum_{\substack{f^{\prime\prime},f^{\prime},f \\ f^{\prime\prime}\# f^{\prime}\# f = g}} f^{\prime\prime} \overt f^{\prime} \overt f = ((\mathrm{id}_{\bclf}\overt \Delta) \circ \Delta)(g).
	\end{equation*}
	Equation \eqref{eqn:conunitpp} is trivial.
\end{proof}

We proceed now with the definition of a vertical product on levelled forests.
\begin{Definition}[monoidal product on levelled forests]
Given two forests $f$ and $f^{\prime}$ with ${{\sf nt}(f^{\prime})=|f|}$, we define $\nabla(f \overt f^{\prime})$ as the sum of forests obtained by first stacking $f^{\prime}$ up to $f$ and then shuffling the generations of $f^{\prime}$ with the generations of $f$ (see Section \ref{sec:levelledforests} for the definition of the action of a permutation on the generations of a forest),
\begin{equation}	\label{eqn:nabla}
	\nabla(f\overt f^{\prime}) = \sum_{s\in\lmss{Sh}(\|f\|,\|f^{\prime}\|)}s\cdot(f \, \# \, f^{\prime}).
\end{equation}
\end{Definition}
The associativity of the product $\nabla$ is easily checked. The unit
$\eta\colon \mathbb{C}_{\overt} \rightarrow \bclf$ is defined by $\eta(1_{m})=\sbt^{m}$.
Let $n\geq 1$, recall that we denote by $c_{n}$ the maximal element for the Bruhat order in $\mathfrak{S}_{n}$:
\begin{equation*}
	c_{n} = n(n-1)\cdots 321.
\end{equation*}
For example, $c_{1} = 1$, $c_{2} = 21$, $c_{3} = 321$, $c_{4} = 4321$.
Given these notions, we state the main theorem of the section
\begin{Theorem}\label{thm:monoidlevelledforest}
$(\bclf, \nabla, \eta, \Delta, \varepsilon)$ is a conilpotent Hopf algebra in the category $(\mathrm{Coll}_{2}$, $\overt, \mathbb{C}_{\overt})$.
\end{Theorem}
To achieve this result we introduce an explicit antipode map.
\begin{Definition}Pick $n,m\geq 1$ two integers. Let $f\in \mathcal{L}\mathcal{F}(n,m)$ be a levelled forest and define its \emph{vertical mirror symmetric} $f^\star \in \mathcal{L}\mathcal{F}(n,m)$ by
\begin{equation*}
f^{\star} = c_{\|f\|} \cdot f.
\end{equation*}
	We extend $\star$ as a conjugate-linear morphism on the bicollection $\bclf$.
\end{Definition}
\begin{Proposition}	\label{prop:antipode}
	Let $f$ be a levelled forest. The map $S\colon \bclf \rightarrow \bclf$ defined by
	\begin{equation*}
		S(f) = (-1)^{\| f \|} f^{\star}
	\end{equation*}
	is an antipode:
	$\nabla\circ(S \overt \mathrm{id}_{\bclf}) \circ \Delta = \nabla\circ(\mathrm{id}_{\bclf} \overt S) \circ \Delta = \varepsilon \circ \eta$.
\end{Proposition}
\begin{proof}
	Let $a$, $b$ be two integers greater than one. Set $n=a+b$. The set of shuffles $\lmss{Sh}(a,b)$ is divided into two mutually disjoint subsets, the set of shuffles sending $a$ (the subset $\lmss{Sh}(a,b)_{+}$) to $n$ and the set of shuffles that do not (resp. $\lmss{Sh}(a,b)_{-}$).

	Recall that if $f$ is a forest then $f^{k}_{-}$ denotes the forest obtained by extracting the~$k$ first lowest generations of~$f$ and~$f^{k}_{+}$ denotes the forest obtained by extracting the~$k$ highest generations of~$f$.
	By definition, one has
	\begin{equation*}
		\nabla (f^{\prime} \overt f\backslash f^{\prime} ) = \sum_{s\in \lmss{Sh}(\|f^{\prime}\|,\|f\backslash f^{\prime}\|)} s \cdot ( f^{\prime} \,\#\,(f\backslash f^{\prime})^{\star} ), \qquad f^{\star} = c_{\|f\|} \cdot f.
	\end{equation*}
	The following relation is easily checked and turns to be the cornerstone of the proof:
	\begin{equation}
		\label{eqn:antipode:telescopique}
		\tilde{s} \circ (c_{n} \otimes \mathrm{id}_{m}) = s \circ (c_{n+1} \otimes \mathrm{id}_{m-1}),\qquad s \in \lmss{Sh}(m-1,n+1)_{-},
	\end{equation}
	with $\tilde{s}$ the unique shuffle in $\lmss{Sh}(m,n)_{+}$ such that $\tilde{s}(m)=n+m$, $\tilde{s}(i)=s(i)$.
	Set $\bar{S}(f) = (-1)^{\|f\|}f^{\star}$. We prove by induction that $S=\bar{S}$. Assume that $S(f)=\bar{S}(f)$ for any forest $f$ with at most $N\geq 1$ generations and pick a forest $f$ with $N+1$ generations. Then, from the induction hypothesis we get
	\begin{gather*}
		S(f) + f + \big(\mathrm{id} \overt \bar{S}\big) \circ \bar{\Delta} (f) = 0,\\
		\nabla \circ \big(\mathrm{id}\overt\bar{S}\big) \circ \bar{\Delta} (f) = \sum_{f^{\prime} \subset f} (-1)^{\|f\backslash f^{\prime}\|} \sum_{s\in\lmss{Sh}(\|f^{\prime}\|, \|f\backslash f^{\prime} \|)} s \cdot \big[f^{\prime}\,\#\,(f\backslash f^{\prime})^{\star} \big] \\
\hphantom{\nabla \circ (\mathrm{id}\overt\bar{S}) \circ \bar{\Delta} (f) }{}
= \sum_{k=1}^{\|f\| -1} (-1)^{k} \sum_{s\in\lmss{Sh}(\|f\|-k,k)}s\cdot \big[ f_{-}^{\| f\| -k } \,\#\, \big(f_{+}^{k}\big)^{\star} \big] .
	\end{gather*}
	We divide the sum over the set {\rm Sh}$(\|f\|-k,k)$ into two sums. The first sums ranges over the subset $\lmss{Sh}(\|f\|-k,k)_{+}$ and the second one ranges overs $\lmss{Sh}(\|f\|-k, k)_{-}$. Then, we gather the sums over $\lmss{Sh}(\|f\|-k,k)_{+}$ and $\lmss{Sh}(\|f\|-k+1,k-1)_{-}$:
	\begin{gather*}
		 \nabla \circ \big(\mathrm{id}\overt \bar{S}\big) \circ \bar{\Delta} =\sum_{k=2}^{\|f\|-2} (-1)^{k}\sum_{s \in \lmss{Sh}(\|f\|-k,k)_{+}} s\cdot \left[f^{\|f\|-k}_{-}\,\#\,(f_{+}^{k})^{\star} \right] \\
\hphantom{\nabla \circ \big(\mathrm{id}\overt \bar{S}\big) \circ \bar{\Delta} =}{}
-(-1)^{k} \sum_{s \in \lmss{Sh}(\|f\|-k+1,k-1)_{-}} s\cdot \big[ f^{\|f\|-k+1}_{-} \,\#\, \big(f_{+}^{k-1}\big)^{\star} \big] \\
\hphantom{\nabla \circ \big(\mathrm{id}\overt \bar{S}\big) \circ \bar{\Delta} =}{}
+(-1)\sum_{s\in\lmss{Sh}(1,\|f\|-1)_{-}} s\cdot \big[ f^{\|f\|-1}_{-}\,\#\,\big(f_{+}^{1}\big)^{\star} \big] \\
\hphantom{\nabla \circ \big(\mathrm{id}\overt \bar{S}\big) \circ \bar{\Delta} =}{}
+ (-1)^{\|f \|-1}\sum_{s\in\lmss{Sh}(\|f\|-1,1)_{+}} s\cdot \big[f^{1}_{-} \# \big(f_{+}^{\|f \|-1}\big)^{\star}\big].
	\end{gather*}
	Using equation \eqref{eqn:antipode:telescopique}, the right-hand side of the last equation is equal to
	\begin{gather*}
	 \nabla \circ \big(\bar{S} \overt \mathrm{id}\big) \circ \bar{\Delta} = 0- \sum_{s\in\lmss{Sh}(1,\|f\|-1)_{-}} s\cdot \big[ f^{1}_{-} \,\#\, \big(f_{+}^{\|f\|-1}\big)^{\star} \big] \\
\hphantom{\nabla \circ \big(\bar{S} \overt \mathrm{id}\big) \circ \bar{\Delta} =}{}
+ (-1)^{\|f\|-1}\sum_{s\in\lmss{Sh}(1,\|f\|-1)_{+}} s\cdot \big[ f^{\|f\|-1}_{-}\,\#\ \big(f_{+}^{1}\big)^{\star}\big] \\
\hphantom{\nabla \circ \big(\bar{S} \overt \mathrm{id}\big) \circ \bar{\Delta}}{}
 = -f + (-1)^{\|f\|-1}f^{\star}.\tag*{\qed}
	\end{gather*}\renewcommand{\qed}{}
\end{proof}

We defined the three structural morphisms $\nabla$, $\Delta$, $S$. To turn~$\lpbf$ into a Hopf monoid, we have to check compatibility between the coproduct $\Delta$ and the product~$\nabla$; the coproduct $\Delta$ should be a morphism of the monoid $(\bclf, \nabla$). This only makes sense provided that we can define a~product on the tensor product $\bclf\overt\bclf$.

Recall that if $f$ is a levelled forest and $0\leq k \leq \|f\|$, one denotes by $f_{-}^{k}$ the levelled subforest of $f$ corresponding to the~$k$ generations at the bottom of $f^{\prime}$: $t_{f_{-}^{k}}$ is the planar subforest of~$t_f$ with a set of internal vertices the set of internal vertices of $f$ labelled by an integer less than~$k$ and for leaves the vertices (including the leaves) of $t_{f}$ connected to one of the latter internal vertices. The levelled forest $f_{+}^{k}$ is obtained similarly by extracting the $k$ top generations of $f^{\prime}$.

\begin{figure}[!ht]\centering
\includegraphics{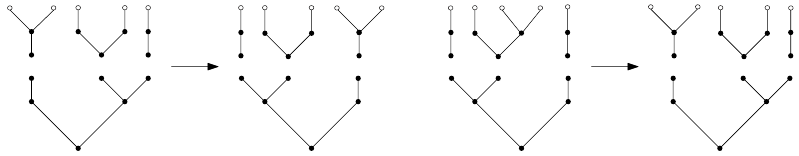}
\caption{Actions of the braiding map.}\label{fig:braiding}
\end{figure}

With $p,q \geq 1$ two integers, we denote by $\tau_{p,q}$ the shuffle in ${\rm Sh}(p,q)$ satisfying $\tau_{p,q}(1) = q+1$ and $\tau_{p,q}(p) = p+q$.
\begin{Definition}	\label{def:braidingmap}
	Define the braiding map
\begin{equation*}
{\sf K} \colon \ \bclf \overt \bclf \rightarrow \bclf \overt \bclf
\end{equation*}
by, for $g$ and $f$ levelled forests such that $f \overt g \in \bclf \overt \bclf$,
	\begin{equation*}
		{\sf K}(f \overt g) = \big(\tau_{\|f\|,\|g\|} \cdot (f \# g)\big)_{-}^{\|g\|} \overt \big(\tau_{\|f\|,\|g\|} \cdot (f \# g)\big)_{+}^{\|f\|}.
	\end{equation*}
\end{Definition}
We pictured in Figure~\ref{fig:braiding} examples of the action of the braiding map on pairs of levelled forests.

We defined the braiding map ${\sf K}$ as acting on $\bclf\overt\bclf$. We extend ${\sf K}$ as a $2$-functor on the product of the monoid generated by $\bclf$ in (Coll$_{2}$,$\overt$). This means in particular that for integers $p,q \geq 1$, we define a bicollection morphism
\begin{equation*}
{\sf K}_{p,q}\colon \ \bclf^{\overt p} \overt \bclf^{\overt q}\rightarrow \bclf^{\overt q} \overt \bclf^{\overt p}.
\end{equation*}

Pick $f_{1}\overt\cdots\overt f_{p} \in \bclf^{\overt p}$ and $g_1\overt\cdots\overt g_q \in \bclf^{\overt q}$.
We define the levelled forest $\succ_{c}h$ with $c\vDash_0 n$ to be the element in $\bclf^{\overt p}$ obtained by the following iterative application of horizontal splittings:
\begin{gather*}
\succ^{-}_c {}={} \succ_{c_1}^-( \cdots (\succ_{c_1+ \cdots + c_{n-1}}^- (h))\cdots ),
\\
{\sf K}_{p,q}(f_{1}\overt\cdots\overt f_{p} \overt g_{1} \overt \cdots \overt g_{q}) \\
\qquad{} ={} \succ^-_{\|g_{1}\|,\dots,\|g_{p}\|, \|f_{1}\|,\dots,\|f_{q}\|}\big(\tau_{\|f_{1}\|+\cdots+\|f_{p}\|, \|g_{1}\|+\cdots+\|g_{q}\|} \cdot f_{1} \# \cdots \# g_{q}\big).
\end{gather*}

The collection of morphisms ${\sf K}_{p,q}$ yields a $2$-functor on the category with objects $\bclf^{\otimes p}$ but with restricted classes of morphisms.
First, it is not difficult to see that $K = {\sf K}_{1,1}$ is an involution and therefore that ${\sf K}_{p,q}$ is an involution too, for any $p,q \geq 1$. It follows from the fact that given~$f$ and~$g$ two levelled forests,
\begin{equation*}
\big[\tau_{||f||,||g||} \cdot (f \# g)\big]^{\|g\|}_{-} \# \big[\tau_{\|f\|,\|g\|}\cdot (f \# g)\big]_+^{\|f\|} = \tau_{\|f\|,\|g\|} \cdot (f \# g),
\end{equation*}
which yields
\begin{gather*}
{\sf K}^2(f \overt g)= \big[\tau_{\|g\|,\|f\|}\cdot\tau_{\|f\|, \|g\|}(f \# g)\big]^{\|f\|}_- \overt \big[\tau_{\|g\|,\|f\|}\cdot\tau_{\|f\|, \|g\|}(f \# g)\big]_+^{\|g\|}
= f \overt g,
\end{gather*}
since $\tau_{\|g\|,\|f\|}\cdot\tau_{\|f\|, \|g\|} = \mathrm{id}_{\|f\|+\|g\|}$.
\begin{Definition}
Let $p,q\geq 0$ be integers and $\varphi\colon \bclf^{\otimes p} \to \bclf^{\otimes q}$, we say that $\varphi$ is \emph{gluing equivariant} if $\varphi$ commutes with the operations $|\!|_{\raisebox{-1mm}{$_i$}}$,namely, with $f_1\overt\cdots\overt f_m \in \bclf^{\overt m}$,
\begin{equation*}
 \varphi(|\!|_{\raisebox{-1mm}{$_i$}}[f_1\overt \cdots \overt f_p]) = |\!|_{\raisebox{-1mm}{$_i$}}[\varphi(f_1\overt\cdots\overt f_p)],\qquad 1 \leq i \leq k,
\end{equation*}
where
\begin{equation*}
|\!|_{_i}[f_1\overt \cdots \overt f_p] := \succ^-_{\|f_1\},\dots,\|f_p\|}|\!|_{\raisebox{-1mm}{$_i$}}[f_1 \sharp \cdots \sharp f_p].
\end{equation*}
\end{Definition}
We denote by ${\rm Hom}_{\rm eq}(p,q)$ the class of all gluing equivariant morphisms between $\bclf^{\overt p}$ and~$\bclf^{\overt q}$. Note that the identity morphisms are gluing equivariant and that the composition of two gluing equivariant morphisms is gluing equivariant.
Also, for each $p,q\geq 0$, $K_{p,q}$ is gluing equivariant.
\begin{Proposition}\label{proposition2.18} The monoid generated by the bicollection $\bclf$ in $(\mathrm{Coll}_{2},\overt)$ with morphisms restricted to the gluing equivariant morphisms is a symmetric monoidal category with symmetry constraints $({\sf K}_{p,q})_{p,q \geq 0}$,
	\begin{equation*}
		{\sf K}_{p,q} \circ {\sf K}_{q,p} = \mathrm{id} \qquad \text{and} \qquad (\mathrm{id}_{\bclf^{\overt q}}\overt {\sf K}_{p,r})\circ ({ \sf K}_{p,q} \overt \mathrm{id}_{\bclf^{\overt r}}) = {\sf K}_{p,q+r}.
	\end{equation*}
\end{Proposition}
\begin{proof}
	Both assertions are trivial and rely on the following relations between the permutations~$\tau_{p,q}$, $p,q\geq 0$:
\begin{equation*}
		\tau_{p,q} \circ \tau_{q,p}=\mathrm{id},\qquad (\mathrm{id}_{q} \otimes \tau_{p,r}) \circ (\tau_{p,q} \otimes \mathrm{id}_{r}) = \tau_{p,q+r},\qquad p,q,r\geq 0.\tag*{\qed}
\end{equation*}\renewcommand{\qed}{}
\end{proof}

Using the above-defined symmetry constraint $K$, we can endow the two-fold tensor product $\bclf\overt\lpbt$ with an algebra product:
\begin{equation*}
	(\nabla \overt \nabla) \circ (\mathrm{id} \overt {\sf K} \overt \mathrm{id})\colon \ \bclf^{\overt 4} \rightarrow \bclf^{\overt 2}.
\end{equation*}

\begin{Proposition}	\label{prop:verticalhopf}
	The two bicollection morphisms $\Delta \colon \bclf\to\bclf\overt\bclf$ and $\nabla\colon \bclf\overt\bclf\rightarrow\bclf$ are vertical algebra morphisms. With $\nabla^{(2)} = \nabla \circ (\nabla \overt \mathrm{id}) = \nabla \circ (\mathrm{id} \overt \nabla)$, this means that
	\begin{equation*}
		\nabla^{(2)} = \nabla^{(2)} \circ (\mathrm{id}\overt {\sf K} \overt \mathrm{id}),\qquad (\nabla \overt \nabla) \circ (\mathrm{id} \overt {\sf K} \overt \mathrm{id}) \circ (\Delta \overt \Delta) = \Delta \circ \nabla.
	\end{equation*}
\end{Proposition}
\begin{Remark}We can rephrase the fact that $\nabla$ is an algebra morphism by saying that $(\bclf,\nabla)$ is, in fact, a commutative algebra.
\end{Remark}
\begin{proof}
	We begin with the first assertion. Pick $f_{1}$, $f_{2}$, $f_{3}$, $f_{4}$ compatible levelled forests (the number of inputs of~$f_{i}$ matches the number of outputs of $f_{i+1}$, $1 \leq i\leq 3$),
	\begin{gather*}
		 \big(\nabla^{(2)} \circ {\sf K}\big)(f_{1} \overt f_{2} \overt f_{3} \overt f_{4}) \\
		 \quad{} =\sum_{s \in {\rm Sh}(\|f_{1}\|,\|f_{3}\|,\|f_{2}\|,\|f_{4}\|)} s \cdot \big[f_{1}\# \big(\tau_{\|f_{2}\|,\|f_{3}\|} \cdot (f_{2} \# f_{3})\big)_{-}^{\|f_{3}\|} \# \big(\tau_{\|f_{2}\|,\|f_{3}\|} \cdot (f_{2} \# f_{3})\big)_{+}^{\|f_{2}\|} \# f_{4} \big] \\
		\quad{} =\sum_{s \in {\rm Sh}(\|f_{1}\|,\|f_{3}\|,\|f_{2}\|,\|f_{4}\|)} \big(s(\mathrm{id}\otimes \tau_{\|f_{2}\|, \|f_{3}\|})\big) \cdot \big(f_{1} \# f_{2} \# f_{3} \# f_{4} \big) \\
		\quad{} =\sum_{s \in {\rm Sh}(\|f_{1}\|,\|f_{2}\|,\|f_{3}\|,\|f_{4}\|)} s \cdot \big(f_{1} \# f_{2} \# f_{3} \# f_{4} \big) = \nabla^{(2)}(f_{1}\overt f_{2} \overt f_{3}\overt f_{4}).
	\end{gather*}
	For the second assertion, we write first
	\begin{align*}
		(\Delta \circ \nabla)(f \overt g) = \sum_{1 \leq k \leq \|f\|+\|g\|} \sum_{s \in {\rm Sh}(k,\|f\|+\|g\|-k)} (s\cdot(f \# g))_{-}^{k} \overt (s\cdot(f\# g))_{+}^{\|f\|+\|g\|-k}.
	\end{align*}
	For each integer $1 \leq k \leq \|f\|$, we split the set of shuffles ${\rm Sh}(\|f\|,\|g\|)$ according to the cardinal~$q$ of the set $s^{-1}(\llbracket 1,k \rrbracket) \cap \llbracket \|f\|+1,\|f\|+\|g\| \rrbracket $. Then a shuffle $s \in {\rm Sh}(\|f\|,\|g\|)$ $s = (s_{1}\otimes s_{2}) \circ \tilde{\tau}_{k,q}$ with $\tilde{\tau}_{k,q}$ the unique shuffle that sends the interval $\llbracket \|f\|+1, \|f\|+q \rrbracket$ to the interval $\llbracket k-q+1,k\rrbracket$ and fixes the interval $\llbracket \|f\|+q+1,\|f\|+\|g\| \rrbracket$,
\begin{gather*}
		\sum_{\substack{1 \leq k \leq \|f\|, 1 \leq q \leq \|g\|, \\ 1 \leq q\leq k} } \sum_{\substack{s_{1} \in {\rm Sh}(k-q,q),\\ s_{2} \in {\rm Sh}(\|f\|-(k-q)\|g\|-q) }} ((s_{1} \otimes s_{2}) \circ \tilde{\tau}_{k,q}) \cdot(f \# g))_{-}^{k} \\
\qquad{} \overt ((s_{1} \otimes s_{2}) \circ \tilde{\tau}_{k,q} \cdot(f\# g))_{+}^{\|f\|+\|g\|-k}.
\end{gather*}
	Notice that $\tilde{\tau}_{k,q} = \tau_{k-q,q}$ and
\begin{equation*}
		\tilde{\tau}_{k,q} \cdot (f \, \# \, g) = f_{-}^{k-q} \# \big(\tau_{\|f\|-(k-q),q}\cdot\big(f_{+}^{\|f\|-(k-q)}\,\#\, g_{-}^{q}\big)\big) \# g_{+}^{\|g\|-q}.
\end{equation*}
	It follows that
	\begin{align*}
		(s_{1} \otimes s_{2}) \circ \tilde{\tau}_{k,q} \cdot(f \# g))_{-}^{k} & = \big((s_{1} \otimes \mathrm{id}) \cdot f_{-}^{k-q} \# \big(\tau_{\|f\|-(k-q),q}\cdot\big(f_{+}^{\|f\|-(k-q)}\,\#\, g_{-}^{q}\big)\big) \# g_{+}^{\|g\|-q}\big)^{k}_{-} \\
&= s_{1} \cdot f_{-}^{k-q} \# \big(\tau_{\|f\|-(k-q),q} \cdot f_{+}^{\|f\|-(k-q)} \# g_{-}^{q}\big)_{-}^{q}.
	\end{align*}
	Similar computations show that
\begin{gather*}
		((s_{1} \otimes s_{2}) \circ \tilde{\tau}_{k,q} \cdot(f\# g))_{+}^{\|f\|+\|g\|-k} \\
\qquad{} = s_{2}\cdot \big(\big(\tau_{\|f\|-(k-q),q}\cdot\big(f_{+}^{\|f\|-(k-q)} \,\#\, g_{-}^{q}\big)\big)_{+}^{\|f\|-(k-q)} \# g_{+}^{\|g\|-q}\big).
\end{gather*}
	The case $\|f\|+1\leq k\leq \|f\|+\|g\|$ is similar, we split the set of shuffles {\rm Sh}$(\|f\|,\|g\|)$ according to the cardinal of the set $s^{-1}(\llbracket k+1, \|f\|+\|g\|\rrbracket) \cap \llbracket 1,\|f\| \rrbracket $) and omitted for brevity. Finally, we obtain for $\Delta \circ \nabla (f \overt g)$ the expression:
	\begin{gather*}
\sum_{\substack{1 \leq k \leq \|f\|, \\ 1 \leq q \leq \|g\|}} \sum_{ \substack{ s_{1} \in {\rm Sh}(k,q) \\ s_{2} \in {\rm Sh}(\|f\|-k,\|g\|-q)}} s_{1} \cdot \big(f_{-}^{k} \# \big(\tau_{\|f\|-k,q}\cdot\big(f_{+}^{\|f\|-k}\,\#\, g_{-}^{q}\big)\big)_{-}^{q}\big) \\
\qquad{} \overt s_{2}\cdot\big(\big(\tau_{k,q}\cdot\big(f_{+}^{\|f\|-k} \,\#\, g_{-}^{q}\big)\big)_{+}^{\|f\|-k} \# g_{+}^{\|g\|-q}\big),
\end{gather*}
which is easily seen to be equal to
$(\nabla \overt \nabla) \circ (\mathrm{id} \overt {\sf K} \overt \mathrm{id}) \circ (\Delta \overt \Delta)(f \overt g)$.
\end{proof}

By collecting altogether the statements of Proposition~\ref{prop:verticalhopf} (proving compatibility between the product~$\nabla$ and the coproduct~$\Delta$) and Proposition~\ref{prop:antipode} proves Theorem~\ref{thm:monoidlevelledforest}. Notice that compatibility between the coproduct and coproduct makes sense because the monoid generated by $\mathcal{L}\mathcal{F}$ is symmetric (for the symmetry constraint ${\sf K}$) as stated in Proposition~\ref{proposition2.18}.

\section{Iterated integrals of a path as operators}\label{sec:iteratedintegrals}
Let us fix a smooth path, $X \colon [0,1] \rightarrow \mathcal{A}$. In this section, we use the algebraic tools developed previously and introduce \emph{partial-} and \emph{full contraction} operators. These operators are indexed by levelled forests and provide a different perspective on the iterated integrals of $X$, as a representation of the monoid of levelled forest introduced in the previous section (see Theorem~\ref{thm:monoidlevelledforest}), rather than as a sequence of tensors.

\subsection{Full and partial contraction operators}\label{sec:contractionoperators}
In what follows, we will intensively use the identifications in the previous sections between levelled trees, and levelled forests and their corresponding permutation and split permutations. In what follows, for any couple of Bananch spaces $A$, $B$ we use the notation ${\rm Hom}(A, B)$ to denote the set of linear continuous maps between $A$ and $B$.
\begin{Definition}\label{def:contractionoperators}
	For any integer $n \geq 1$ and levelled tree $\sigma$ in $\lpbt(n)$, we define the \emph{full contraction of $X$ along $\sigma$} as the map $\mathbb{X}^{\gs}\colon [0,1]^2\to {\rm Hom} \big(\cA^{\otimes n}, \cA\big)$ given for any $A_1, \dots , A_n\in \cA$ by
	\begin{equation}\label{defn_full_contraction}
		\mathbb{X}_{s,t}^{\sigma}(A_{1}\otimes \cdots \otimes A_{n})= \int_{\Delta_{s,t}^{n-1}} A_{1}\cdot\mathrm{d}X_{t_{\sigma(1)}} \cdot A_2 \cdots \mathrm{d}X_{t_{\sigma(n-1)}}\cdot A_{n},
	\end{equation}
	where $\sigma$ is identified with a permutation in $\mathfrak{S}_{n-1}$ when $n\geq 2$ and $\mathbb{X}_{s,t}^{\bullet}= \mathrm{id}_{\mathcal{A}}$.
\end{Definition}
\begin{Remark}The above definition may be misleading since the identity \eqref{defn_full_contraction} defines a linear map on the algebraic tensor product, whereas we used $\otimes$ to denote the projective tensor product. However, the algebraic tensor product is a dense subspace of $\cA^{\otimes n}$ and we interpret $\mathbb{X}^{\sigma}_{s,t}$ as the unique continuous operator extending the values in \eqref{defn_full_contraction}. Similar considerations apply throughout the paper.
\end{Remark}
If linearly extended to the vector space spanned by all levelled trees (or equally permutations), the map $\sigma\mapsto\mathbb{X}^{\sigma}_{s,t}$ yields naturally a morphism between the collection $\bclt$ in \eqref{levelled_trees} and the endomorphism collection ${\rm End}_{\mathcal{A}}$ given by
\begin{equation*}
	 {\rm End}_{\mathcal{A}}(n)={\rm Hom}\big(\mathcal{A}^{\otimes n},\mathcal{A}\big),\qquad n\geq 1,
\end{equation*}
see Appendix \ref{appendix}. The partial contraction operators, that we now introduce, extend $\sigma\mapsto\mathbb{X}^{\sigma}_{s,t}$ to a morphism between $\bclf$ to many-to-many operators, i.e., elements of ${\rm Hom}_{\mathrm{ Vect}_{\mathbb{C}}}\big(\mathcal{A}^{\otimes n},\mathcal{A}^{\otimes m}\big)$, $m<n$.

To properly define them, we denote by ${\rm End}^{2}_{\mathcal{A}}$ the bicollection of \emph{noncommutative polynomials on multilinear maps on~$\mathcal{A}$ with values in} $\mathcal{A}$. That is using the notation in the Appendix~\ref{appendix}
\def\endd{{{\rm End}_{\mathcal{A}}^{2}}}
\begin{equation*}
{\rm End}_{\mathcal{A}}^{2}(m,n) := T({\rm End}_{\mathcal{A}})(m,n)=\bigoplus_{\substack{k_1+\cdots + k_m = n}} {\rm End}_{\mathcal{A}}(k_1)\otimes\cdots\otimes {\rm End}_{\mathcal{A}}(k_m),
\end{equation*}
when $n\geq 1$ and $m \geq 1$ and the condition $k_1+\cdots + k_m = n$ is satisfied for some integer $k_1, \dots, k_m\geq 1$. Moreover, we set ${\rm End}^2_{\mathcal{A}}(0,0) = \mathbb{C}$ and ${\rm End}^2_{\mathcal{A}}(m,n)=0$ otherwise. The bicollection $\endd$ is endowed with a monoidal product $\nabla_{{\rm End}^{2}_{\mathcal{A}}}$ associated to the vertical tensor product~$\overt$. This operation extends the usual canonical operadic structure $\circ$ on ${\rm End}_{\mathcal{A}}$ as a monoidal morphism, see Appendix~\ref{appendix}. For example, given two non-trivial elements $v\in {\rm End}_{\mathcal{A}}^{2}(k,n)$, $v= v_1 \otimes \cdots \otimes v_k$ with and $u\in {\rm End}_{\mathcal{A}}^{2}(m,k)$, $u=u_1 \otimes \cdots \otimes u_m $ one has
\begin{equation*}
	\nabla_{{\rm End}^{2}_{\mathcal{A}}}(u \overt v) = (u_1\circ v_{1}\otimes\cdots\otimes v_{k_1})\otimes \cdots \otimes (u_m\circ(v_{k-k_m+1}\otimes\cdots\otimes v_{k}).
\end{equation*}

We introduce also the double tensor algebra of $\mathcal{A}$. We start from the unital tensor algebra
\begin{equation*}
\bar{T}(\mathcal{A}) = \mathbb{C}\varnothing \oplus \bigoplus_{n\geq 1} \mathcal{A}^{\otimes n}.
\end{equation*}
Elements of $ \bar{T}(\mathcal{A})$ are linear combinations of words $w=a_1\cdots a_k:=a_1\otimes \cdots \otimes a_k$ and $\varnothing$ is the unity for concatenation of words. The double tensor algebra is given by
\begin{equation*}
\bar{T}^2(\mathcal{A}):=\bar{T}(\bar{T}(\mathcal{A}))=\mathbb{C} 1 \oplus\bigoplus_{n\geq 1}\bar{T}(\mathcal{A})^{\otimes n}.
\end{equation*}
Elements of $\bar{T}^2(\mathcal{A})$ are represented as words of words. To distinguish the internal concatenation of $\bar{T}(\mathcal{A})$ and the second order concatenation, we use the symbol $|$ for the concatenation product on $\bar{T}\big(\bar{T}(\mathcal{A})\big)$ and the symbol $1$ stands now for the unit of $|$. In the following, for any integers $m \geq 1$ and $n\geq 1$ we use the notation
\begin{equation*}
\bar{T}^{2}(\mathcal{A})(m,n) := \bigoplus_{\substack{n_1+\cdots + n_m = n \\ n_i \geq 1}} \mathcal{A}^{\otimes (n_1-1)}\otimes\cdots\otimes\mathcal{A}^{\otimes (n_m-1)}
\end{equation*}
with the convention that $\mathcal{A}^{0}=\mathbb{C}\varnothing$. For both words on words in $\mathcal{A}$ and words on endomorphisms in $\mathcal{A}$, we freely identify the sequence of vector spaces $\endd(m,n)$ (resp.~$\bar{T}^2(\mathcal{A})(m,n)$) with their direct sum. We will however make clear this distinction for other collections and bicollections. We relate $\bar{T}^{2}(\mathcal{A})$ and ${\rm End}_{\mathcal{A}}^{2}$ via an explicit representation.
\begin{Definition}[representation of the algebra $\bar{T}^{2}(\mathcal{A})$] We define a representation $\mathrm{Op}\colon \bar{T}^2(\mathcal{A}) \rightarrow \mathrm{End}^2_{\mathcal{A}}$ of the algebra $\big(\bar{T}^{2}(\mathcal{A}),\, |\,\big)$ extending the following values, for $A_{1}\otimes\cdots\otimes A_{n} \in T(\mathcal{A})$, $X_{i} \in \mathcal{A}$,
\begin{gather*}
{\mathrm {Op}}(A_{1}\otimes\cdots\otimes A_{n})(X_{0}\otimes\cdots\otimes X_{n}) = X_{0}\cdot A_{1} \cdots A_{n}\cdot X_{n}, \\
{\rm Op}(\varnothing)(X_1)={\rm id}_{\mathcal{A}}.
\end{gather*}
\end{Definition}

The representation Op has one crucial property. By definition, ${\rm Op}$ is compatible with the concatenation product $|$ on $\bar{T}^2(\mathcal{A})$. As explained, ${\rm End^{2}_{\mathcal{A}}}$ is endowed with a vertical monoidal structure $\nabla_{{\rm End}^{2}_{\mathcal{A}}}$. The same kind of structure exists on $\bar{T}^2(\mathcal{A})$. Indeed, $\bar{T}(\mathcal{A})$ can be endowed with an operadic structure $\circ$, that we call \emph{words insertions}. Given a word $a_1\otimes \cdots \otimes a_n \in \bar{T}(\mathcal{A})$ and $w_0,\dots,w_n \in \bar{T}(\mathcal{A})$, one defines
\begin{equation*}
	a_1 \cdots a_n \circ (w_0\otimes\cdots\otimes w_n) := w_0 \otimes a_1 \otimes w_1\cdots w_{n-1}\otimes a_n \otimes w_n.
\end{equation*}
One can check that $\circ$ satisfies the associativity and unitality constraints of an operadic composition. We then extend this operadic composition as a horizontal monoidal morphism and define this way an associative product
\begin{equation*}
\nabla_{\bar{T}^2(\mathcal{A})}\colon \ \bar{T}^2(\mathcal{A}) \overt \bar{T}^2(\mathcal{A}) \to \bar{T}^2(\mathcal{A}).
\end{equation*}

Then ${\rm Op}$ is compatible with respect to the products $\nabla_{{\rm End}_{\mathcal{A}}^{2}}$ and $\nabla_{\bar{T}^2(\mathcal{A})}$. That is
\begin{equation}\label{compatibility_OP}
\big(\nabla_{\rm End^{2}_{\mathcal{A}}} \overt \nabla_{\rm End^{2}_{\mathcal{A}}}\big) \circ ({\rm Op} \overt {\rm Op}) = {\rm Op} \circ \nabla_{\bar{T}^2(\mathcal{A})}.
\end{equation}

\begin{Example} Pick $a_1,a_2 \in \mathcal{A}$ and $b_1,b_2,b_3 \in \mathcal{A}$ and consider
\begin{equation*}
w_1=a_1|a_2,\qquad w_2 = b_1 b_2 | \varnothing | b_3 |\varnothing.
\end{equation*}
The two words $w_1$ and $w_2$ are compatible, since $w_1$ has four inputs and $w_2$ has four outputs, we compose them together
\begin{equation*}
 \nabla_{\bar{T}^{2}(\mathcal{A})}(w_1\overt w_2)=b_1b_2a_1 | b_3a_2,
\end{equation*}
and apply Op to the result,
\begin{equation}\label{eqn:opopnabla}
 {\rm Op}(b_1b_2a_1 | b_3a_2)(X_0,\dots,X_6) = X_0\cdot b_1\cdot X_1\cdot b_2\cdot X_2\cdot a_1\cdot X_3\otimes X_4\cdot b_3 \cdot X_5 \cdot a_2 \cdot X_6.
\end{equation}
We can apply first ${\rm Op}$ to $w_1$ and $w_2$ and compose together the resulting operators,
\begin{gather*}
 {\rm Op}(w_1)(Y_0,Y_1,Y_2, Y_3) = Y_0\cdot a_1\cdot Y_1\otimes Y_2\cdot a_2\cdot Y_3,\\
 {\rm Op}(w_2)(X_0,\dots,X_6) = X_0\cdot b_1\cdot X_1\cdot b_2\cdot X_2\otimes X_3 \otimes X_4\cdot b_3\cdot X_5 \otimes X_6.
\end{gather*}
Substituting to $Y_0\otimes Y_1 \otimes Y_2\otimes Y_3$ the right-hand side of the last equation, we recover \eqref{eqn:opopnabla},
\begin{equation*}
\nabla_{\bullet{\rm End}^{2}_{\mathcal{A}}}(v \overt u) = (v_1\circ (u_1\otimes\cdots\otimes u_{n_1}))\cdots (v_p\circ(u_{n_1+\cdots+n_{i-1}}\otimes\cdots\otimes u_{n_1+\cdots+n_{i}})).
\end{equation*}
\end{Example}

The vector space $\bar{T}(\mathcal{A})$ is the natural space wherein the signature of a smooth path $X$ takes values, see \eqref{equ:signature}. To define partial contractions we need to implement the freedom of permutations and the double tensor algebra inside the usual signature. Let $w$ be a word in $T(\mathcal{A})$ with length $n \geq 1$. Let $c=(c_{1},\dots,c_{k})$ be a composition of $n$.
The composition $(c_{1},\dots,c_{k})$ yields a splitting of $w$: we define the element $\left[w\right]_{(c_{1},\dots,c_{k})}\in \bar{T}^{2}(\mathcal{A})$ by
\begin{equation*}
	\left[w\right]_{(c_{1},\dots,c_{k})}= w_{1}\cdots w_{c_{1}}| w_{c_{1}+1}\cdots w_{c_{1}+c_{2}}| \cdots |w_{c_{1}+\cdots+c_{k-1}+1} \cdots w_{c_{1}+\cdots+c_{k}},
\end{equation*}
with the convention $w_{c_{1}+\cdots+c_{i-1}+1}\cdots w_{c_{1}+\cdots+c_{i-1}+c_{i}}=\varnothing$ if $c_{i}=0$.
\begin{Definition}
For any integer $n \geq 1$ and levelled tree $\sigma$ in $\lpbt(n)$, we denote by $X^{\gs}$ the map $ X^{\gs}\colon [0,1]^2\to \bar{T}(A)$ given by
\begin{equation*}
	X_{s,t}^{\sigma}= \int_{\Delta_{s,t}^{n-1}} \mathrm{d}X_{t_{\sigma(1)}} \otimes \cdots \otimes \mathrm{d}X_{t_{\sigma(n-1)}},
\end{equation*}
where $\sigma$ is identified with a permutation in $\mathfrak{S}_{n-1}$ when $n\geq 2$ and $\mathbb{X}_{s,t}^{\bullet}=\varnothing$. For any levelled forest $f=(\sigma,c)\in \lpbf(m,n)$ we denote by $X^f$ the application $X^f_{s,t}=\big[{{X}}_{s,t}^{\sigma} \big]_{c}$.
\end{Definition}

\begin{Definition}\label{def:pcontractionoperators}
For any integers $n \geq 1$, $n\geq m\geq 1$ and any levelled forest $f=(\sigma,c)\in \lpbf(m,n)$, we define the \emph{partial-contraction of $X$ along the forest $f$} as a map $\mathbb{X}^{f}\colon [0,1]^2\to \mathrm{Hom}\big(\cA^{\otimes n}, \cA^{\otimes m}\big)$ given by
\begin{equation*}
 \mathbb{X}_{s,t}^{f} = \mathrm{Op}\big(X^f_{s,t}\big) .
\end{equation*}
\end{Definition}
\begin{Example}Let us calculate the partial contraction associated with the levelled forest $f$ in Figure~\ref{fig:levelledforest}. In this case, $n=5$ and the word on words representing $f$ is
\begin{equation*}
f=(13245, (2,0,2,1,0)).
\end{equation*}
We associate it with the formal expression
\begin{equation*}
{\rm d}X_{t_1} {\rm d}X_{t_3} \otimes \varnothing \otimes {\rm d}X_{t_2}{\rm d}X_{t_4} \otimes {\rm d}X_{t_5} \otimes \varnothing.
\end{equation*}
The first term on the left of the above expression corresponds to the first tree in $f$, with two vertices labelled $1$ and $3$. It is followed on its right by a straight tree, yielding the first $\varnothing$.
We then interleave a $10$-tuple $(A_0, \dots, A_9)$ of elements in $\mathcal{A}$ between each ${\rm d}X_{t_i}$, replacing the empty letter $\varnothing$ by one of the $A$'s,
\begin{equation*}
A_0 \cdot {\rm d}X_{t_1} \cdot A_1 \cdot {\rm d}X_{t_3} \cdot A_2\otimes A_3 \otimes A_4 \cdot {\rm d}X_{t_2}\cdot A_5\cdot {\rm d}X_{t_4}\cdot A_6 \otimes A_7\cdot {\rm d}X_{t_5} \cdot A_8 \otimes A_9.
\end{equation*}
Finally, we integrate over $\Delta^n_{s,t}$ and obtain the following formula for $\mathbb{X}_{s,t}^{f}(A_0, \dots,A_9)$,
\begin{gather*}
\int_{\Delta_{s,t}^5}\! (A_0 \cdot {\rm d}X_{t_1} \cdot A_1 \cdot {\rm d}X_{t_3} \cdot A_2 )\otimes A_3 \otimes ( A_4 \cdot {\rm d}X_{t_2} \cdot A_5 \cdot {\rm d}X_{t_4} \cdot A_6)\otimes( A_7 \cdot {\rm d}X_{t_5} \cdot A_8) \otimes A_9.
\end{gather*}
\end{Example}

\subsection{Chen relation}\label{sec:chenrelation}
In this section, we use from time to time the symbol $\circ$ in place of $\nabla_{\bar{T}^2(\mathcal{A})}$ or $\nabla_{\rm End^{2}_{\mathcal{A}}}$ to improve readability. We describe how the concatenation of paths lifts to the full and partial contractions operators, that is we write an extension of Chen identity over iterated integral, see \cite{Chen54} for these operators.

\begin{Proposition}[Chen relation]\label{prop:chenrelation}
For any forest $f\in\lpbf$ and any three times {$(s,u,t)\in [0,1]^3$} one has
	\begin{equation}\label{eq:chen_eq}
		\mathbb{X}^{f}_{s,t} = \sum_{f^{\prime} \subset f} \nabla_{\rm End_{\mathcal{A}}^{2}} \big[\mathbb{X}_{u,t}^{f^{\prime}} \overt \mathbb{X}^{f\backslash f^{\prime}}_{s,u} \big].
	\end{equation}
Written in term of the notations introduced in \eqref{eqn:coproduct} and the map $\mathbb{X}_{s,t}\colon \mathcal{L}\mathcal{F} \rightarrow \mathrm{End}^2_{\mathcal{A}}$ defined by $\mathbb{X}_{s,t}(f) = \mathbb{X}_{s,t}^{f}$, the equation \eqref{eq:chen_eq} becomes
 \begin{equation*}
		{\mathbb{X}}_{s,t} = \nabla_{\mathrm{End}^2(\mathcal{A})}\circ (\mathbb{X}_{u,t} \overt \mathbb{X}_{s,u}) \circ \Delta.
	\end{equation*}
\end{Proposition}
\begin{Example}
Before writing the proof, we check equation \eqref{eq:chen_eq} on an explicit example given by the levelled forest $f=(213,(2,1))$ in Figure~\ref{fig:exsubforest}, to see how the operations combine themselves. In that case, the operator $\mathbb{X}^{f}_{s,t}$ is given by
\begin{equation*}
\mathbb{X}^{f}_{s,t}(A_0, \dots, A_5 )= \int_{\Delta_{s,t}^3} ( A_0 \cdot {\rm d}X_{t_2} \cdot A_1 \cdot {\rm d}X_{t_1} \cdot A_2 ) \otimes (A_3 \cdot {\rm d}X_{t_3} \cdot A_4 ).
\end{equation*}
Using the standard properties of Lebesgue integration, we can easily write
\begin{gather}
 \mathbb{X}^{f}_{s,t}= \mathbb{X}^{f}_{s,u}+\mathbb{X}^{f}_{u,t}+ \int_{t_1\in \Delta_{u,t}^1}\int_{(t_2,t_3)\in\Delta^2_{s,u}} ( A_0 \cdot {\rm d}X_{t_2} \cdot A_1 \cdot {\rm d}X_{t_1} \cdot A_2 ) \otimes (A_3 \cdot {\rm d}X_{t_3} \cdot A_4 )\nonumber\\
 \hphantom{\mathbb{X}^{f}_{s,t}=}{}
 + \int_{(t_1,t_2)\in\Delta_{u,t}^2}\int_{t_3\in \Delta^1_{s,u}} ( A_0 \cdot {\rm d}X_{t_2} \cdot A_1 \cdot {\rm d}X_{t_1} \cdot A_2 ) \otimes (A_3 \cdot {\rm d}X_{t_3} \cdot A_4 ).\label{eq:ex_2}
\end{gather}
At the same time, we list all subforests in $f'\subset f$ in Figure \ref{fig:exsubforest} together with $f\setminus f'$.
\begin{figure}[!ht]
$$
\begin{array}{ |c|c| }
 \hline
 f & f\setminus f' \\
 \hline
 (\varnothing, (0,0)) & (213,(2,1)) \\
 (1, (1,0)) & (23, (1,1)) \\
 (21, (2,0)) & (3, (0,1)) \\
 (213,(2,1)) & (\varnothing, (0,0)) \\
 \hline
\end{array}
$$
\caption{The subforests in Figure \ref{fig:exsubforest} presented as words on words in the first column. In the second column, the cut of $f$ by each of these forests.}\label{fig:tableforests}
\end{figure}

Proposition \ref{prop:chenrelation} implies the equality
\begin{equation*}
 \mathbb{X}^{f}_{s,t}= \mathbb{X}^{f}_{s,u}+\mathbb{X}^{f}_{u,t}+ \mathbb{X}^{(1, (1,0))}_{u,t} \circ \mathbb{X}^{(23, (1,1))}_{s,u} + \mathbb{X}^{(21, (2,0))}_{u,t} \circ \mathbb{X}^{(3, (0,1))}_{s,u},
\end{equation*}
which is exactly \eqref{eq:ex_2}.
\end{Example}

\begin{proof} It is sufficient to show the identity when $s<u<t$. The statement of the proposition is implied by the same statement but for the iterated integrals $X^{f}_{s,t}$, $f \in \lpbf$ since Op is a~representation of the word-insertions operad (see equatio~\eqref{compatibility_OP}). We prove the identity
\begin{equation*}
		X^{f}_{s,t} = \sum_{f^{\prime} \subset f} X_{u,t}^{f^{\prime}} \circ X^{f\backslash f^{\prime}}_{s,u}
\end{equation*}
by induction on the generation of $f$ and with $\circ$ the operation $\nabla_{\rm \bar{T}^{2}(\mathcal{A})}$. The initialization is done for forests with $0$ generations. Assume that the results as been proved for forests having at most $N$ generations and let $f$ be a forest with $N+1$ generations. Splitting the simplex $\Delta_{s,t}^{n+1}$ according to $s<u<t$ one has
\begin{gather}	\label{eqn:recursivestepchenrelation}
	X_{s,t}^{f} = \int_{s}^{t}{\rm d}{X}_{t_{1}} \circ X_{s,t_{1}}^{f\backslash f_{1}}+ \int_{u}^{t}{\rm d}{X}_{t_{1}} \circ X_{s,t_{1}}^{f\backslash f_{1}}= X_{s,u}^{f} + \int_{u}^{t}{\rm d}{X}_{t_{1}} \circ X_{s,t_{1}}^{f\backslash f_{1}},
\end{gather}
where $f_{1}={\succ}_1(f)^-$ and $\circ$ is the word insertion of the element
\begin{equation*}
	\mathrm{d}X_{t_{1}} = \varnothing^{\otimes i-1} \otimes \mathrm{d}X_{t_{1}} \otimes \varnothing^{|f|-i},
\end{equation*}
where $i$ is the order of the $i^{\rm th}$ tree in the forest $f$ whose root is decorated by $1$. By construction of~$f\backslash f_{1}$, this forest has only~$N$ generations and the recursive hypothesis to the forest $f \backslash f_{1}$ implies
\begin{equation*}
	X_{s,t_{1}}^{f\backslash f_{1}} = \sum_{f'' \subset f\backslash f_{1}} X_{u,t_{1}}^{f''} \circ X^{(f\backslash f_{1}) \backslash f''}_{s,u}.
\end{equation*}
We insert this last relation into equation~\eqref{eqn:recursivestepchenrelation} to get the identity
	\begin{equation*}
	\int_{u}^{t} \mathrm{d}X_{t_{1}} \circ X_{u,t_{1}}^{f^{\prime}} = \sum_{f'' \subset f\backslash f_{1}} \int_{u}^{t} \mathrm{d}X_{t_{1}} \circ X_{u,t_{1}}^{f''} \circ \big[X^{(f\backslash f_{1}) \backslash f''}_{s,u}\big]= \sum_{f^{\prime} \subset f\, f'\neq \varnothing} X_{u,t}^{f^{\prime}} \circ X^{(f\backslash f')}_{s,u}.\tag*{\qed}
	\end{equation*}\renewcommand{\qed}{}
\end{proof}

\begin{Remark}We apply the above formula to the levelled tree $f$ which is a right comb tree obtained by grafting corollas with two leaves with each other, always on their rightmost node. By cutting such a tree at a certain level, we obtain on one hand a smaller comb tree and on the other hand, we obtain a levelled forest with only straight trees, except for the last one, the rightmost, which is also a comb tree.
 By denoting $\textrm{comb}_n$ the comb tree with $n$ internal nodes, we thus get for a tuple $A_0,\dots,A_n \in \mathcal{A}$,
\begin{equation*}
 \mathbb{X}_{s,t}^{c_n} (A_0,\dots,A_n) = \sum_{k=0}^n \mathbb{X}_{u,t}^{c_k}\big(A_0,\dots,A_{k-1},\mathbb{X}^{c_{n-k}}_{s,u}(A_k,\dots,A_n)\big),
\end{equation*}
 more explicitly
 \begin{gather*}
 \int_{\Delta^{n}_{s,t}}A_0\cdot \mathrm{d}X_{t_1}\cdots \mathrm{X}_{t_n} \cdot A_n
 = \sum_{k=0}^{n} \int_{\Delta^{k}_{u,t}} A_1\cdot \textrm{d}X_{t_1} \cdot A_2 \cdots \mathrm{d}X_{t_{k-1}}\cdot A_k\cdots\\
 \hphantom{\int_{\Delta^{n}(_{s,t}}A_0\cdot \mathrm{d}X_{t_1}\cdots \mathrm{X}_{t_n} \cdot A_n =}{}
 \times \int_{\Delta^{n-k}_{s,u}}A_k \cdot \mathrm{d}X_{u_1} \cdots \mathrm{d}X_{u_{n-k}}\cdot A_n.
 \end{gather*}
 The above relation is implied by contracting the famous Chen relation satisfied by the tensor iterated integrals, with elements
$A_0,\dots,A_n$. This is the only one, among~\eqref{eqn:recursivestepchenrelation} when $f$ ranges levelled forests, implied by a linear transformation of the Chen relation for tensor iterated integral (contraction). When $f$ is not a comb tree, $\mathbb{X}_{s,t}^f$ will expand following~\eqref{eqn:recursivestepchenrelation} over iterated integrals contracted with elements of~$\mathcal{A}'s$ but with \emph{permuted noises}~$\mathrm{d}X_t$.
\end{Remark}

To the family of operators $\big\{\mathbb{X}^{f}_{s,t}, \, f \in \lpbt \big\}$, we now associate a family of
endomorphisms on
\begin{equation*}
	\lpbt(\mathcal{A}) := \bigoplus_{n\geq 0} \mathcal{A}^{\otimes n} \otimes \bigoplus_{\substack{\tau \in \lpbt \\ |\tau|=n}}\mathbb{C}[\tau] .
\end{equation*}
For the remaining part of the article, we use the lighter notations
\begin{equation*}
a\otimes \tau = a \tau \in \mathcal{A}^{\otimes |\tau|}\otimes \mathbb{C}\tau, \qquad \lpbt(\mathcal{A})(\tau):= \mathcal{A}^{\otimes |\tau|} \tau \subset \lpbt(\mathcal{A}),
\qquad
\mathbb{X}_{\cdot} = \mathrm{id}\otimes \mathbb{X}_{\cdot} \circ \Delta.
\end{equation*}

Although it is not yet clear if it is possible to associate to the full and partial contractions operators a path on a certain convolution group of representations, our statement of the Chen relation makes clear that any prospective deconcatenation product $\Delta$ should act on a tree by cutting it in all possible ways, generations after generations. In Section~\ref{sec:aconvolutiongp} we prove this cutting operation yields a comonoid structure on $(\mathcal{L}\mathcal{F}, \overt)$.

From Proposition \ref{prop:chenrelation} we immediately deduce the following properties.
\begin{Proposition}\label{prop:model}
The family of maps $\bar{\mathbb{X}}_{s,t}$ defined by
\begin{align}
		\bar{\mathbb{X}}_{s,t}\colon \ \lpbt(\mathcal{A}) & \rightarrow \lpbt(\mathcal{A}),\nonumber \\
 a \tau & \mapsto \sum_{\tau^{\prime}\subset \tau}\bar{\mathbb{X}}_{s,t}^{\tau \backslash \tau^{\prime}}(a) \tau^{\prime}\label{eqn:barncrp}
\end{align}
have the following properties:
\begin{enumerate}\itemsep=0pt
\item[$(1)$] for every levelled tree $\tau\in \lpbt$
\begin{equation*}
(\bar{\mathbb{X}}_{s,t}-\mathrm{id})(\lpbt(\mathcal{A})(\tau)) \subset \bigoplus_{\tau^{\prime} \subsetneq \tau} \lpbt(\mathcal{A})(\tau^{\prime}),
\end{equation*}
\item[$(2)$] 
for any $(s,u,t)\in [0,1]^3$ the so-called \emph{noncommutative Chen's relations} hold
\begin{equation*}
\bar{\mathbb{X}}_{s,t} = \bar{\mathbb{X}}_{u,t} \circ \bar{\mathbb{X}}_{s,u}.
\end{equation*}
\end{enumerate}
\end{Proposition}

\begin{proof}It is sufficient to show only point (2) for any $s < u <t$, since point~(1) is trivial.
For any given and $A \tau \in \lpbt(\mathcal{A})$, to the Chen's relation in Proposition~\ref{prop:chenrelation} implies
	\begin{align*}
		\bar{\mathbb{X}}_{s,t}(A \tau) = \sum_{\tau^{\prime} \subset \tau} \bar{\mathbb{X}}_{s,t}^{\tau\backslash \tau^{\prime}}(A)\cdot \tau^{\prime} & = \sum_{\tau^{\prime} \subset \tau} \sum_{\tau^{\prime\prime} \subset \tau\backslash \tau^{\prime}} \bar{\mathbb{X}}_{u,t}^{\tau^{\prime\prime}} \big( \bar{\mathbb{X}}_{s,u}^{(\tau \backslash \tau^{\prime}) \backslash \tau^{\prime\prime}}(A_{1}\otimes\cdots\otimes A_{|\tau|})\big) \tau^{\prime} \\
	 & =\sum_{\tau^{\prime} \subset \tau} \sum_{\tau^{\prime\prime} \subset \tau\backslash \tau^{\prime}} \bar{\mathbb{X}}_{u,t}^{\tau^{\prime\prime}} \big( \bar{\mathbb{X}}_{s,u}^{\tau \backslash (\tau^{\prime\prime} \sharp \tau^{\prime}) }(A_{1}\otimes\cdots\otimes A_{|\tau|})\big) \tau^{\prime}.
	\end{align*}
	By performing the change of variable $g = f^{\prime\prime} \sharp f^{\prime}$, $g^{\prime} = f^{\prime}$, we obtain
	\begin{align*}
		\bar{\mathbb{X}}_{s,t}(A_{1}\otimes\cdots\otimes A_{|f|}\, f)&=\sum_{g \subset f} \sum_{g^{\prime} \subset g} \bar{\mathbb{X}}_{u,t}^{g \backslash g^{\prime} } \big( \bar{\mathbb{X}}_{s,u}^{f \backslash g}(A_{1}\otimes\cdots\otimes A_{|f|})\big) g^{\prime} \\
&= \big(\bar{\mathbb{X}}_{u,t} \circ \bar{\mathbb{X}}_{s,t}\big)(A_{1} \otimes \cdots \otimes A_{|f|} f).\tag*{\qed}
	\end{align*}\renewcommand{\qed}{}
\end{proof}

\subsection{Integration by part properties}\label{sec:geometricproperties}

Another important property of iterated integral are integration by part formulae, e.g., the identity
\begin{equation*}
	\int_{\Delta^2_{s,t}} \mathrm{d}X_{t_{1}}\otimes \mathrm{d}X_{t_{2}} + \int_{\Delta^2_{s,t}} \mathrm{d}X_{t_{2}}\otimes \mathrm{d}X_{t_{1}} = (X_{t}-X_{s})\otimes (X_{t}-X_{s}).
\end{equation*}

We investigate the consequences of such identities at the level of full-partial contractions and their related associated endomorphisms $\bar{\mathbb{X}}_{s,t}$. These relations imply a specific ``compatibility condition'' with respect to a certain product on $\lpbt(\mathcal{A})$ defined from the shuffle operations on trees in Definition~\ref{shuffle_product_forests}. We start with the consequences of partial contractions.
\begin{Proposition}
For any given couple of forests $f$ and $g$ with ${\sf nt}(g)=|f|$ and any $(s,t)\in [0,1]^2$ one has
\begin{equation}\label{eq:product_eq}
	\nabla_{\rm End_{\mathcal{A}}^{2}}\big(\mathbb{X}^{f}_{s,t} \overt \mathbb{X}^{g}_{s,t}\big)= \sum_{s\in\lmss{Sh}(\|f\|,\|g\|)}\mathbb{X}_{s,t}^{s\cdot(f \shuffle g)} .
\end{equation}
Written in term of the notations introduced in \eqref{eqn:nabla} and the map $\mathbb{X}_{s,t}\colon \mathcal{L}\mathcal{F} \rightarrow \mathrm{End}^2_{\mathcal{A}}$ defined by $\mathbb{X}_{s,t}(f) = \mathbb{X}_{s,t}^{f}$, the equation~\eqref{eq:product_eq} becomes
 \begin{equation*}
		 \nabla_{\mathrm{End}^2(\mathcal{A})} (\mathbb{X}_{s,t} \overt \mathbb{X}_{s,t}) = \mathbb{X}_{s,t} \nabla.
	\end{equation*}
\end{Proposition}

\begin{proof}
Let us fix $s<t$. Writing the forests $f$ and $g$ as the split permutations $f=(\sigma_1, c_1)$ and $g=(\sigma_2, c_2)$, it follows from the shuffle identity for iterated integrals of $X$ that one has the identity
	\begin{gather*}
	 X^{\sigma_1}_{s,t}\otimes X^{\sigma_2}_{s,t} =\sum_{s\in{\rm Sh}(|f|-1,|g|-1)}\int_{\Delta^{|f|+|g|-2}_{s,t}}\mathrm{d}X_{t_{(s\circ\sigma)(1)}}\otimes\cdots\otimes\mathrm{d}X_{t_{(s\circ\sigma)(|f|+|g|-2)}},
	\end{gather*}
where $\sigma =\sigma_{1}\otimes\sigma_{2}$. By applying the split permutation of $f\# g$ to both sides we deduce the identity
 \begin{gather*}
 X^{f\shuffle g}_{s,t} =\sum_{s\in{\sf Sh}(|f|-1,|g|-1)}\int_{\Delta^{|f|+|g|-2}_{s,t}} X^{s\cdot (f\shuffle g)}_{s,t}.
 \end{gather*}
Composing with $\mathrm{Op}$ we conclude.
\end{proof}

We restate this identity at the level of $\bar{\mathbb{X}}_{s,t}$. This task can be done by introducing an operadic composition $L$ on a collection of words with entries in $\mathcal{A}$ different from before. Together with the shuffle product on levelled trees, this operadic composition yields a \emph{structural} map ${\sf L}$ on~$\lpbt^{\sharp}(\mathcal{A})$. Further properties of~${\sf L}$ will turn central in better understanding the composition of the Taylor series for the fields~$a$,~$b$ in equation~\eqref{eqn:controlleddiff}.

\def\fs{\mathcal{F}\mathcal{S}}
\def\ltsa{{{\sf LT}^{\#}(\mathcal{A})}}
\begin{Definition}[faces substitution]
 We define the collection of vector spaces $\fs$ by
 \begin{equation*}
 \fs(n) = \mathcal{A}^{\otimes n}, \qquad n\geq 1.
 \end{equation*}
 Next, define $L \colon \fs\circ\fs \rightarrow \fs$ as follows. Pick a word $U \in \mathcal{A}^{\otimes p+1}$ and words $A^{i} \in \mathcal{A}^{\otimes m_{i}}$, $1 \leq i \leq p+1$, $A^i=\big( A^i_{(1)} \otimes \cdots \otimes A^i_{(m_i)}\big) $ and set
 \begin{gather*}
	 L\big(U\otimes A^{1} \otimes \cdots \otimes A^{p}\big) :=\big( U_{(1)}\cdot A^{1}_{(1)}\big) \otimes A^{1}_{(2)} \otimes \cdots\\
\hphantom{L\big(U\otimes A^{1} \otimes \cdots \otimes A^{p}\big) :=}{}
 \otimes \big(A_{(m_{1})}^{1} \cdot U_{(2)}\cdot A^{2}_{(1)}\big) \otimes \cdots \otimes \big(A^{p}_{(m_{p})}\cdot U_{p}\big).
\end{gather*}
The word $\1\otimes \1$ acts as the unit for $L$.
\end{Definition}

We denote by ${\sf FS}$ the graded vector space equal to the direct sum of all vector spaces in the collection $\fs$. Notice that elements of $\mathcal{A}$ are $0$-ary operators in the collection $\fs$ and, for example, the above formula for $L$ gives $L(U_{1}\otimes U_{(2)} \circ A) = U_{1}\cdot A\cdot U_{2} \in \mathcal{A}$, with $U_{(1)}\otimes U_{(2)}\in\mathcal{A}^{\otimes 2}$.

\begin{Proposition} ${\rm FS}=(\fs, L, \1\otimes \1)$ is an operad.
\end{Proposition}
\begin{proof}
The following proposition holds and rests on the associativity of the product on $\mathcal{A}$.
\end{proof}

\begin{figure}[!ht] \centering
	\includegraphics{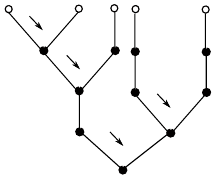}
	\caption{The faces of a levelled tree are indicated with arrows.}\label{fig:treewithfacesasentries}
\end{figure}

In the collection $\fs$, a word with length $n$ is an operator with $n-1$ entries, the inner gaps between the letters. So far, a levelled tree was considered as an operator with as many inputs as it has of leaves. However, there is an alternative way to see such a tree as an operator: by considering the \emph{faces} of the tree as inputs. A face is a region enclosed between two consecutive leaves and delimited by two paths of edges meeting at the least common ancestor, see Figure~\ref{fig:treewithfacesasentries}.

We denote by $\lpbt^{\#}$ the set of levelled trees graded by the numbers of faces, $\lpbt^{\#}(n)$ the set of levelled trees with $n$ faces, and ${\sf LT}^{\sharp}(\mathcal{A})$ the space ${\sf LT}(\mathcal{A})$ seen as a graded vector space with ${\sf LT}^{\#}(\mathcal{A})(n) = \mathbb{C}{\sf LT}^{\#}(n)\otimes \fs(n)$. Notice that the endomorphism $\bar{\mathbb{X}}_{s,t}$ we defined in the previous section satisfies:
\begin{equation*}
	\big(\bar{\mathbb{X}}_{s,t}-\mathrm{id}\big)\big(\lpbt^{\#}(\mathcal{A})(n)\big) \subset \bigoplus_{k<n} \lpbt^{\#}(\mathcal{A})(k).
\end{equation*}
We also set for any levelled tree $\tau \in \lpbt$,
\begin{equation*}
 \ltsa(\tau) = \mathbb{C}\tau \otimes\mathcal{FS}(\|\tau\|).
\end{equation*}
The space $\ltsa$ is equipped with an involution $\starltsa$, defined by
\begin{equation}\label{starlta}
\starltsa (A_1\otimes \cdots \otimes A_p \tau) = A_p^{\star}\otimes\cdots\otimes A_1^{\star} \theta (\tau).
\end{equation}
The graded vector space $\lpbt^{\#}(\mathcal{A})$ yields a collection $\mathcal{L}\mathcal{T}^{\#}(\mathcal{A})$ by setting the space $n$-ary operators $\mathcal{L}\mathcal{T}^{\#}(\mathcal{A})$ equal to $\lpbt^{\#}(\mathcal{A})(n)$. We set abusively
\begin{equation*}
	\lpbt^{\#}(\mathcal{A}) \circ \lpbt^{\#}(\mathcal{A}) := \bigoplus_{n\geq 0}\big(\mathcal{L}\mathcal{T}^{\sharp}(\mathcal{A})\circ \mathcal{L}\mathcal{T}^{\sharp}(\mathcal{A})\big)(n).
\end{equation*}
\def\lts{\lpbt^{\#}}
For any $U,A \in \lpbt^{\#}(\mathcal{A})$, we introduce the notation
\begin{equation*}
	U \circ A := \sum_{\substack{\tau,\tau_{1},\dots,\tau_{\|\tau\|}}} U^{\tau}\, \tau \otimes A^{\tau_{1}} \, \tau_{1} \otimes \cdots \otimes A^{\tau_{\|f\|}} \, \tau_{\| \tau \|} \in \lpbt^{\#}(\mathcal{A}) \circ \lpbt^{\#}(\mathcal{A}) .
\end{equation*}
Observe that $U\circ A$ is only linear on $U$, not on $A$.

\begin{Definition}\label{def:operatorL}
 Define the morphism of graded vector spaces
\begin{gather}
{\sf L} \colon \ \lpbt^{\#}(\mathcal{A}) \circ \lpbt^{\#}(\mathcal{A}) \rightarrow \lpbt^{\#}(\mathcal{A}),\nonumber\\
\label{eq:operatorL}
	 {\sf L}\big(U \otimes \big(A_{1} \otimes\cdots\otimes A_{\|\alpha\|}\big)\big) = \sum_{\substack{\alpha \\\tau_{1},\dots,\tau_{\|\alpha\|}}} L\big(U^{\alpha}\otimes A^{\tau_{1}}_{1} \otimes \cdots \otimes A^{\tau_{\|\alpha\|}}_{{\|\alpha\|}}\big) \cdot \tau_{1}\shuffle \cdots \shuffle \tau_{\|\alpha\|},
\end{gather}
where
\begin{equation*}
U = \sum_{\alpha \in \lpbt^{\#}}U^{\alpha} \alpha, \qquad A_{i}=\sum_{\tau_{i} \in \lpbt^{\#}} A_{i}^{\tau_{i}}\,\tau_i\, \in \lpbt^{\#}(\mathcal{A}).
\end{equation*}
\end{Definition}
\begin{Lemma}\label{lemma:shufflexbar}
	Let $\alpha$ and $\beta$ be two levelled trees in $\lpbt$. For any $A \in \mathcal{A}^{\otimes |\alpha|+|\beta|-1}$ one has
	\begin{equation}\label{lemma:shuffle}
		\bar{\mathbb{X}}_{s,t}(A \alpha \shuffle \beta) = \sum_{\tau_{\alpha} \subset \alpha, \tau_{\beta} \subset \beta } \bar{\mathbb{X}}_{s,t}^{(\alpha \backslash \tau_{\alpha}) \shuffle (\beta \backslash \tau_{\beta})}(A) \tau_{\alpha} \shuffle \tau_{\beta}.
	\end{equation}
\end{Lemma}
\begin{proof}The proof consists essentially of a re-summation. It stems from the definition of the map $\bar{\mathbb{X}}_{s,t}$ that
\begin{equation}\label{eqn1:lem:operadxbarsimple}
		\bar{\mathbb{X}}_{s,t}(A \cdot \alpha \shuffle \beta) = \sum_{\substack{\tau \in \alpha\shuffle \beta \\ \tau^{\prime} \subset \tau}} \bar{\mathbb{X}}_{s,t}^{\tau \backslash \tau^{\prime}}(A) \tau^{\prime}.
	\end{equation}
	Let $\tau \in \alpha \shuffle \beta$ a tree obtained by shuffling vertically the generations of $\alpha$ and $\beta$ and pick $\tau^{\prime} \subset \tau$ a subtree. Let $s$ be the shuffle in ${\rm Sh}(\sharp \alpha, \sharp \beta)$ such that $\tau^{-1} = (\alpha \otimes \beta) \circ s^{-1}$. We associate to the pair $(\tau,\tau^{\prime})$ a triple which consists in the tree $\tau$, and two others trees $\tau^{\prime}_{\alpha} \subset \alpha$ and $\tau^{\prime}_{\beta} \subset \beta$ satisfying
	\begin{equation*}
		\tau^{\prime} = (\tau^{\prime}_{\alpha} \otimes \tau^{\prime}_{\beta}) \circ \tilde{s}^{-1},
	\end{equation*}
where $\tilde{s}$ is a shuffle in $\lmss{Sh}\big(\|\tau^{\prime}_{\alpha}\|, \|\tau^{\prime}_{\beta}\|\big)$. Such a permutation $\tilde{s}$ is unique, in fact it is obtained from~$s^{1}$ by extracting the first $\|\tau^{\prime}\|$ letters of the word representing~$s^{-1}$, followed by standardization and finally inversion. Recall that \emph{standardization} means that we translate the first $\|\tau^{\prime}\|$ letters representing~$s^{-1}$, while maintaining their relative order to obtain a word on integers in the interval $\llbracket 1, \sharp \tau^{\prime} \rrbracket$. It is clear that the map $\phi\colon (\tau,\tau^{\prime}) \mapsto (\tau, \tau^{\prime}_{\alpha}, \tau^{\prime}_{\beta})$ is injective. Now, given $\tau_{\alpha} \subset \alpha$, $\tau_{\beta} \subset \beta$, and two shuffles $s_{-} \in \lmss{Sh} (\|\tau^{\prime}_{\alpha}\|, \| \tau^{\prime}_{\beta}\| )$, $s_{+} \in \lmss{Sh}(\|\alpha \backslash \tau^{\prime}_{\alpha}\|, \|\beta \backslash \tau^{\prime}_{\beta}\|)$,
	we define a~third shuffle $s_{-+}$ in $\lmss{Sh}(\|\alpha\|, \|\beta\|)$ by requiring
\begin{gather*}
		s_{-+}(i) = s_{-}(i), \qquad 1 \leq i \leq \|\tau^{\prime}_{\alpha}\|,\\ s_{-+}(\|\tau^{\prime}_{\alpha}\| + i) = s_{+}(i) + s_{-}(\| \tau^{\prime}_{\alpha}\|),\qquad 1 \leq i \leq \|\tau_{\alpha} \backslash \tau^{\prime}_{\alpha}\|.
\end{gather*}
	The map $\delta\colon (\tau^{\prime}_{\alpha},\tau^{\prime}_{\beta}, s_{+},s_{-}) \mapsto (\tau, \tau^{\prime}_{\alpha},\tau^{\prime}_{\beta})$ with $\tau^{-1} = \alpha \otimes \beta \circ s_{-+}^{-1}$ is a bijection between the image of $\phi$ and
\begin{equation*}
		\mathcal{S}=\{(\tau_{\alpha},\tau_{\beta},s_{+},s_{-}),\,\tau_{\alpha} \subset \alpha, \,\tau_{\beta} \subset \beta,\, s_{-} \in \lmss{Sh}(\|\tau_{\alpha} \|,\|\tau_{\beta} \|),\, s_{+} \in \lmss{Sh}(\|\alpha \backslash \tau_{\alpha} \|, \|\beta \backslash \tau_{\beta} \|)\}.
\end{equation*}

	We can thus rewrite the sum on the right-hand side of \eqref{eqn1:lem:operadxbarsimple} as follows:
	\begin{gather*}
		\sum_{\substack{\tau \in \alpha \shuffle \beta \\ \tau^{\prime} \subset \tau}} \bar{\mathbb{X}}_{s,t}^{\tau \backslash \tau^{\prime}}(A) \tau^{\prime} = \sum_{\tau_{\alpha},\tau_{\beta},s_{+},s_{-} \in \mathcal{S}} \bar{\mathbb{X}}_{s,t}^{(\alpha \otimes \beta) \circ s_{-+} \backslash (\tau^{\prime}_{\alpha} \otimes \tau^{\prime}_{\beta}) \circ s_{-}^{-1}} (\tau^{\prime}_{\alpha} \otimes \tau^{\prime}_{\beta}) \circ s_{-}^{-1}.
	\end{gather*}
	Now, we observe that the forest $(\alpha \otimes \beta) \circ s_{-+} \backslash (\tau^{\prime}_{\alpha} \otimes \tau^{\prime}_{\beta}) \circ s_{-}^{-1}$ does only depend on the trees~$\tau_{\alpha}$,~$\tau_{\beta}$ and the shuffle $s_{+}$. Summing over all shuffles~$s_{+}$, we get $\alpha \backslash \tau_{\alpha} \shuffle \beta \backslash \tau_{\beta}$. The statement of the lemma follows by computing the sum over~$s_{-}$.
\end{proof}

\begin{Definition}[product on the face-substitution collection]
For any $A\in \mathcal{FS}(n)$ and $B\in \mathcal{FS}(m)$, we define their product $A\cdot B$ by
\begin{equation*}
 A\cdot B = A_{(1)}\otimes\cdots\otimes (A_{(n+1)}\cdot B_{(1)})\otimes\cdots\otimes B_{(m+1)}.
\end{equation*}
The product $\cdot$ is a product on the collection $\mathcal{FS}$ with unit $1\in \mathcal{FS}(0)$,
$A \cdot B \in \mathcal{FS}(n+m)$. Note that we use the same symbol $\cdot$ for the above-defined product and the product on the algebra $\mathcal{A}$ for the reason that the former restricts to the latter on $\mathcal{F}\mathcal{S}(0)$.
\end{Definition}
\begin{Remark}The product $\cdot$ has a very special form, namely,
	\begin{equation*}
		A \cdot B = (\1\otimes \1\otimes \1)\circ(A\otimes B)=L((\1\otimes\1\otimes\1)\otimes (A \otimes B)),
	\end{equation*}
	and the relation $L(m \otimes ({\rm id}_{\sf FS} \otimes m)) = L(m \otimes (m\otimes {\rm id}_{\sf FS})$ with $ m=\1\otimes\1 \otimes \1$ entails associativity of the product $\cdot$. We say that $m \in \fs(2)$ is a \emph{multiplication} in the operad~$(\fs, L)$. In addition, associativity of the operadic composition $L$ results in the following distributivity law
	\begin{equation*}
		(A \cdot B)\circ C = (A\circ B) \cdot (B \circ C),\qquad A,B,C \in {\sf FS}.
	\end{equation*}
\end{Remark}
Conjointly with the shuffle product on levelled trees, the product $\cdot$ brings in a graded algebra product $\shuffle\colon \lpbt^{\#}(\mathcal{A}) \otimes \lpbt^{\#}(\mathcal{A}) \rightarrow \lpbt^{\#}(\mathcal{A})$, namely
\begin{equation}\label{eqn:shuffleproductdeux}
	(A \alpha) \shuffle (B \beta) = (A \cdot B) \alpha \shuffle \beta,
\end{equation}
with unit $\1 \cdot \sbt$. The above relation on full contraction operators yields compatibility of the endomorphism $\bar{\mathbb{X}}_{s,t}$ with the product $\shuffle$ defined in~\eqref{eqn:shuffleproductdeux}.
\begin{Proposition}\label{prop:geometricityun}Let $\alpha$ and $\beta$ be two levelled trees and pick $A \in \mathcal{A}^{|\alpha|}$, $B\in\mathcal{A}^{\otimes |\beta|} $,
	\begin{equation*}
		\bar{\mathbb{X}}_{s,t}((A \alpha) \shuffle (B \beta))=\bar{\mathbb{X}}_{s,t}(A \alpha)\shuffle\bar{\mathbb{X}}_{s,t}(B \beta).
	\end{equation*}
\end{Proposition}
\begin{proof}The result is a simple consequence of the previous Proposition \ref{lemma:shufflexbar} and the shuffle relation for the partial contraction operators~\eqref{lemma:shuffle}. In fact, one has the trivial identities
	\begin{align*}
	\bar{\mathbb{X}}_{s,t}((A\cdot B)\cdot \alpha \shuffle \beta) & = \sum_{\substack{\tau_{\alpha} \subset\alpha, \tau_{\beta} \subset \beta}} {\bar{\mathbb{X}}}^{\alpha\backslash\tau_{\alpha} \shuffle \beta\backslash\tau_{\beta}}_{s,t}(A\cdot B)\, \tau_{\alpha} \shuffle \tau_{\beta} \\& = \sum_{\substack{\tau_{\alpha} \subset\alpha, \tau_{\beta} \subset \beta}} {\bar{\mathbb{X}}}^{\alpha\backslash\tau_{\alpha}}(A) \cdot{\bar{\mathbb{X}}}_{s,t}^{\beta\backslash\tau_{\beta}}(B)\, \tau_{\alpha} \shuffle \tau_{\beta}= \bar{\mathbb{X}}_{s,t}(A)\cdot\bar{\mathbb{X}}_{s,t}(B).
	\end{align*}
Thereby obtaining the desired identity.
\end{proof}

\begin{Corollary}\label{cor:geometricity}
	For all times $0 < s < t < 1$, it holds that
	\begin{equation*}
		{\sf L} \circ \big(\mathrm{id} \circ \bar{\mathbb{X}}_{s,t}\big) = \bar{\mathbb{X}}_{s,t} \circ {\sf L}.
	\end{equation*}
\end{Corollary}
\begin{proof}According to Proposition \ref{prop:geometricityun}, one has the identity
	\begin{gather*}
		 {\sf L}\big(U^{\alpha} \otimes \bar{\mathbb{X}}_{s,t}\big(A^{\beta_{1}} \beta_{1}\big) \otimes \cdots \otimes \bar{\mathbb{X}}_{s,t}\big(A^{\beta_{\sharp \alpha}} \beta_{\sharp \alpha}\big)\big)\\
	\qquad{}	 =\bar{\mathbb{X}}_{s,t}\big(U^{\alpha}_{(1)}\sbt\big) \shuffle \bar{\mathbb{X}}_{s,t}\big(A^{\beta_{1}} \beta_{1}\big) \shuffle \bar{\mathbb{X}}_{s,t}\big(U^{\alpha}_{(2)}\sbt\big) \cdots \bar{\mathbb{X}}_{s,t}\big(A^{\beta_{\sharp \alpha}} \cdot \beta_{\sharp \alpha}\big) \shuffle \bar{\mathbb{X}}_{s,t}\big(U^{\alpha}_{(|\alpha|)} \sbt\big) \\
	\qquad{}	 =\bar{\mathbb{X}}_{s,t}\big(\big(U^{\alpha}_{(1)}\sbt\big) \shuffle \big(A^{\beta_{1}}\,\beta_{1}\big) \shuffle \big(U^{\alpha}_{(2)}\sbt\big) \cdots \big(A^{\beta_{\sharp \alpha}} \cdot \beta_{\sharp \alpha}\big) \shuffle \big(U^{\alpha}_{(|\alpha|)} \sbt\big)\big)\\
		\qquad{} =\bar{\mathbb{X}}_{s,t}\big({\sf L}\big(U^{\alpha} \alpha \otimes A^{\beta_{1}} \beta_{1}\otimes \cdots \otimes A^{\beta_{\sharp \alpha}}\beta_{\sharp \alpha}\big)\big).
	\end{gather*}
From which we deduce the announced equality.
\end{proof}

Summing up the result in this section and the previous one, we can actually state the properties of $\bar{\mathbb{X}}_{s,t}$ by introducing a suitable group of endomorphisms.
\begin{Theorem}\label{thm:pathofrep}Denoting by $G(\mathcal{A})$ the group
\begin{gather}
G(\mathcal{A})=\bigg\{ \bar{\mathbb{X}} \in \mathrm{Hom}_{\mathrm{Alg}}\big(\ltsa, \ltsa\big)\colon \nonumber\\
\hphantom{G(\mathcal{A})=\bigg\{}{}
 \big(\bar{\mathbb{X}} - \mathrm{id}\big)\big(\ltsa(\tau)\big)\subset \bigoplus_{\tau^{\prime}\subset \tau} \ltsa(\tau^{\prime}) \bigg\},\label{eqn:modelspaceun}
\end{gather}
the endomorphisms $\bar{\mathbb{X}}_{s,t}$ satisfy the following properties
\begin{enumerate}\itemsep=0pt
 \item[$(1)$] $\bar{\mathbb{X}}_{s,t} \in G(\mathcal{A})$,
 \item[$(2)$] for every time $(s,u,t) \in [0,1]^3$, $\bar{\mathbb{X}}_{s,t} = \bar{\mathbb{X}}_{u,t}\circ\bar{\mathbb{X}}_{s,u}$.
\end{enumerate}
\end{Theorem}

\section{Group of signatures}
Looking again at the formal Peano--Picard expansion sketched in \eqref{eqn:seriesexpansion}, we see a sum of full contraction operators $\{\mathbb{X}^{\tau}_{s,t}\}_{\tau\in \lpbt} $ applied to generic elements of the algebra $\mathcal{A}$, those operators appear thus as the fundamental objects to generalise in a rough path setting. As explained in the previous section, we were forced to consider partial contraction operators indexed by levelled binary forests to write the Chen relation for these operators. These operators appear as coefficients of an endomorphism $\mathbb{X}$ acting on $\lpbt({\sf \mathcal{A}})$. These ``coefficients'' associated with forests can not be related to the ``coefficients'' associated with trees if $\mathcal{A}$ is truly infinite-dimensional. A bit more formally, the application corestricting an endomorphism $\mathbb{X}_{s,t}$, for any pair of times $s<t$ to $\mathbb{C} [\sbt ]$
\begin{equation*}
	\bar{\mathbb{X}} (\cdot)\mapsto \sum_{\tau \in \lpbt} \bar{\mathbb{X}}^{\tau}_{s,t}(\cdot) \sbt
\end{equation*}
is not injective, \emph{we are lacking relations between partial and full contraction operators}. Worth is the fact, that data of $\mathbb{X}$ is in fact equivalent to the data of all iterated integrals of~$X$, so that the previous section is in fact a mere, though much more intricated reformulation of the classical theory. Yet, we explain in this section how to get rid of these partial contraction operators while maintaining a Chen relation for an object comprising only full contraction operators. These partial contraction operators are turned into ``technical proxies'', that can be constructed from the operators associated with levelled trees and bear no additional information on the small-scale behaviour of the paths but allows for an efficient formulation of the Chen relation.

\subsection{Algebra of face-contractions}\label{sec:facescontraction}
\def\lpbtsfc{\lpbt^{\#}({\sf FC})}

To define a proper group where full contraction take value, we introduce a new collection of operators, that we call \emph{face-contractions}. These operators will replace words with entries in $\mathcal{A}$ of the previous section.%

\begin{Definition}[face-contractions]\label{def:facescontractions}
For any $\tau \in \lpbt$, $\tau\neq \sbt$ and $A_{1}\otimes\cdots\otimes A_{(|\tau|)} \in \mathcal{A}^{|\tau|}$ we associate the \emph{global face-contraction map}
\begin{equation*}
\sharp((A_{1}\otimes\cdots\otimes A_{|\tau|}) \tau)\colon \ \mathcal{A}^{\otimes \|\tau\|}\to \mathcal{A},
\end{equation*}
which evaluates on a tuple $X_{1},\dots,X_{\|\tau\|} \in \mathcal{A}$ as
\begin{equation*}
\sharp(A_{1}\otimes\cdots\otimes A_{|\tau|} \tau)(X_{1},\dots,X_{\|\tau\|}) = A_{1} \cdot X_{\tau^{-1}(1)}\cdot A_2 \cdots \cdot X_{\tau^{-1}}(\|\tau\|) \cdot A_{|\tau|}.
\end{equation*}
We denote by ${\sf FC}(\tau)$ the closure with respect to the operator norm the span of the global face-contraction maps; i.e., for any fixed $\tau \in \lpbt$, $\tau\neq \sbt$ we set \begin{equation*}
{\sf FC}(\tau) := \mathrm{Cl}\big(\big\{\sharp((A_{1}\otimes\cdots\otimes A_{|\tau|}) \, \tau)\colon A_{1}\otimes\cdots\otimes A_{|\tau|} \in \mathcal{A}^{\otimes|\tau|}\big\}\big),
\end{equation*}
and we call the elements of ${\sf FC}(\tau)$ the \emph{$\tau$-face-contractions}. Moreover, we define
\[\lpbtsfc :=\bigoplus_{n\geq 0}\bigoplus_{\substack{\tau \in \lpbt \\ \|\tau\|=n}}{\sf FC}(\tau), \]
 and set ${\sf FC}(\sbt)=\mathcal{A}$.
\end{Definition}
 We introduce for any face-contractions operators $m = \sum_{\tau\in \lpbt} m^{\tau}$
the \emph{face-contraction norm}
\begin{equation*}
	\pmb{|}m\pmb{|} := \sum_{\tau \in \lpbt}\frac{1}{\|\tau\|!} \|m^{\tau}\|,
\end{equation*}
where $\|\cdot \|$ is the usual operator norm induced by $\mathcal{A}$. From the notation $\sharp((A_{1}\otimes\cdots\otimes A_{|\tau|}) \cdot \tau)$ we also denote by $\sharp$ the morphism of graded vector spaces
\begin{align*}
	\sharp \colon \ \lpbt^{\#}(\mathcal{A}) & \to \lpbtsfc, \\
A_{1}\otimes \cdots \otimes A_{(|\tau|)} \tau & \mapsto \sharp (A_{1}\otimes\cdots\otimes A_{(|\tau|)} \tau).
	\end{align*}

Notice that the operator $\sharp((A_{1}\otimes\cdots\otimes A_{|\tau|}) \cdot \tau)$ has $\|\tau\|=|\tau|-1$ inputs. Its output can be computed by drawing a sparse quasi-binary tree $\tau$ and placing $A_1, \dots, A_{\tau}$ up to the leaves of $\tau$ and the ${\rm d}X_{t_i}$ on the unique vertex with two children on the~$i^{\rm th}$ generation of $\tau$. Whereas in the previous section the arguments of the multilinear operators we considered were located on the leaves, in this section they are located on the faces. Some operations we defined on trees can be push-forward via $\sharp$ to define a proper unital Banach algebra with involution. We denote these operations with similar notation as the operations defined over~$\lpbt^{\#}(\mathcal{A})$.

\newcommand{\starltfc}{\star_{\lpbtsfc}}

\begin{Definition}[shuffle product on face-contractions operators]For any $m \in {\sf FC}(\tau)$ and $m^{\prime}\in {\sf FC}(\tau^{\prime})$ we define $m \shuffle m^{{\prime}} \in \lpbtsfc$ and $\star_{{\sf LT}^{\#}({\sf FC})}(m)$ on every tuple $X_1,\dots,X_{\|\tau\|+\|\tau^{\prime}\|} \in \mathcal{A}$ by
\begin{gather*}
(m \shuffle m^{{\prime}}) \big(X_{1}\otimes\cdots\otimes X_{\|\tau\|+\|\tau^{\prime}\|}\big)\nonumber \\
\qquad {} := \sum_{s\in{\rm Sh}(\|\tau\|,\|\tau^{\prime}\|)} m\big(X_{s(1)}\otimes\cdots\otimes X_{s(\|\tau\|)}\big)\cdot m'\big(X_{s(\|\tau\|+1)}\otimes\cdots\otimes X_{s(\|\tau\|+\|\tau'\|)}\big), \\ 
 \star_{\lpbtsfc} (m)\big(X_{1}\otimes\cdots\otimes X_{\|\tau\|}\big) :=\star \big(m\big(\star (X_{\|\tau\|})\otimes\cdots\otimes \star (X_{1})\big)\big).
\end{gather*}
We call $\shuffle$ and $\star_{{\sf LT}^{\#}({\sf FC})} (m)$ the \emph{shuffle product on face-contractions} and the \emph{involution on face-contrac\-tions}.
\end{Definition}

\begin{Proposition}\label{prop:product}The triple $\big(\lpbtsfc, \shuffle, \star, \pmb{|}\cdot \pmb{|}\big)$ is a unital Banach algebra with involution, i.e., for any pairs of levelled trees $\tau$ and $\tau^{\prime}$ and operators $m \in {\sf FC}(\tau)$, $m'\in{\sf FC}(\tau^{\prime})$ one has the properties
\begin{equation*}
	\starltfc(m \shuffle m^{{\prime}})= \starltfc(m^{{\prime}})\shuffle \starltfc(m), \qquad	\pmb{|} m\shuffle m^{\prime}\pmb{|} \leq \pmb{|}m\pmb{|} \pmb{|} m^{\prime} \pmb{|}.
\end{equation*}
Moreover, $\sharp$ is a morphism of unital Banach algebra, i.e., one has the identities
\begin{gather*}
\sharp\big(A_{1}\otimes\cdots\otimes A_{|\tau|} \tau\shuffle B_{1}\otimes\cdots\otimes B_{|\tau'|} \tau'\big) = \sharp\big(A_{1}\otimes\cdots\otimes A_{|\tau|} \tau\big)\shuffle \sharp\big(B_{1}\otimes\cdots\otimes B_{|\tau'|} \tau'\big), \\
\sharp\big(\starltsa(A_{1}\otimes\cdots\otimes A_{|\tau|} \tau)\big) =\starltfc \sharp\big(A_{1}\otimes\cdots\otimes A_{|\tau|} \tau\big)
\end{gather*}
for any $A_1, \dots, A_{|\tau|}$ and $B_1, \dots, B_{|\tau'|}$ of elements in $\mathcal{A}$ and $\tau, \tau'$ in $\lpbt$.
\end{Proposition}
\begin{proof}We first prove the morphism property for $\sharp$. We fix $\tau$, $\tau'$ in $\lpbt$ and adopt the notations $A= A_{1}\otimes\cdots\otimes A_{|\tau|}$ and $B=B_{1}\otimes\cdots\otimes B_{|\tau'|}$. Using the identity~\eqref{eqn:shuffleproductdeux}
and the explicit form of $\tau \shuffle \tau'$ one has
\begin{equation*}
A \tau\shuffle B \tau'= (A\cdot B) \tau\shuffle\tau'= \sum_{\substack{\st( \sigma_1\cdots\sigma_{|\tau|} )=\tau \\ \st(\sigma_{|\tau|+1}\cdots\sigma_{|\tau|+ |\tau'|})= \tau'}} (A\cdot B)\cdot \sigma.
\end{equation*}
By evaluating the right-hand side on a generic tuple $X=(X_{1},\dots,X_{\|\tau\|+\|\tau^{\prime}\|})$, we have
\begin{gather*}
 \sum_{\substack{\st( \sigma_1\cdots\sigma_{|\tau|} )=\tau \\ \st(\sigma_{|\tau|+1}\cdots\sigma_{|\tau|+ |\tau'|})= \tau'}} (A\cdot B)\cdot \sigma(X )
 \\
 \quad{}=\sum_{\substack{\st( \sigma_1\cdots\sigma_{|\tau|} )=\tau \\ \st(\sigma_{|\tau|+1}\cdots\sigma_{|\tau|+ |\tau'|})= \tau'}} A_{1} \cdot X_{\sigma(1)} \cdots X_{\sigma(\|\tau\|)}\cdot A_{|\tau|} \cdot B_{1}\cdot X_{\sigma(\|\tau\|+1)} \cdots X_{\sigma(\|\tau\|+\|\tau^{\prime}\|)}\cdot B_{|\tau'|}\\
 \quad{}= \sum_{s \in \scalebox{0.6}{${\rm Sh}(\|\tau\|, \| \tau^{\prime}\|)$}}\sharp (A \tau)(X_{s(1)},\dots,X_{s(\|\tau\|)}) \cdot \sharp (B \tau^{\prime})\big(X_{s(\|\tau\|+1)},\dots,X_{s(\|\tau\| + \|\tau^{\prime}\|)}\big)\\
\quad{} = \big(\sharp(A \tau)\shuffle \sharp(B \tau') \big)(X).
\end{gather*}
Thereby obtaining the algebra morphism property for $\sharp$. Moreover, from the previous identity, we deduce also the following estimate in terms of the operator norm
	\begin{gather*}
		\| \sharp (A \tau \shuffle B \tau^{\prime}) \| \leq | {\rm Sh}(\|\tau\|,\|\tau^{\prime}\|)|\| \sharp (A \tau )\|\|\sharp (B \tau^{\prime})\|=\frac{(\|\tau\| + \| \tau^{\prime}\|)!}{\|\tau\|! \|\tau^{\prime}\| !} \| \sharp (A \tau )\|\|\sharp (B \tau^{\prime})\|.
	\end{gather*}

From this, we deduce by density the Banach algebra property with respect to the norm $\pmb{|} \cdot \pmb{|}$. Compatibility of $\sharp$ with respect to the involution $\starltsa$ and $\starltfc$ follows from the definition of $\starltsa$ in equation \eqref{starlta} and $\starltfc$ right above.

In addition, for any $m \in {\sf FC}(\tau)$, $m'\in{\sf FC}(\tau^{\prime})$ one has
\begin{gather*}
\starltfc(m\shuffle m') \big(X_1\otimes\cdots\otimes X_{\|\tau\|+\|\tau^{\prime}\|}\big) \\
\qquad{} =\sum_{s}\star_{\mathcal{A}} m'\big(\star_{\mathcal{A}}\big(X_{s(\|\tau\|+\|\tau^{\prime}\|)}\big)\otimes \cdots \otimes \star_{\mathcal{A}}\big(X_{s(\|\tau\|+1)}\big)\big)\\
\qquad\quad{}
\times \star_{\mathcal{A}} m\big(\star_{\mathcal{A}}(X_{s(\|\tau\|)})\otimes \cdots \otimes \star_{\mathcal{A}}(X_{1})\big) \\
\qquad{} =\sum_{s} \starltfc(m')\big(X_{s(\|\tau\|+1)},\dots,X_{s(\|\tau\|+\|\tau^{\prime}\|)}\big)\cdot \starltfc(m)\big(X_{s(1)},\dots,X_{\|\tau\|}\big)\\
\qquad{} =\starltfc(m')\shuffle \starltfc(m),
\end{gather*}
where the sums above are taken over $s$ in ${\rm Sh}(\|\tau\|,\|\tau^{\prime}\|)$.
\end{proof}

\begin{Remark}
The global face-contraction map $\sharp$ yields a morphism of operads. In the same way as formula \eqref{eq:operatorL} introduces an operadic composition on $\mathcal{F}\mathcal{S}$. It is also possible to define an operadic composition, that we denote by the symbol $\tilde{L}$, on the collection $\mathcal{F}\mathcal{C}$ of face-contractions operators,
\begin{equation*}
\mathcal{F}\mathcal{C}(n) := \bigoplus_{n\geq 0}\bigoplus_{\substack{\tau \in \lpbt \\ \|\tau\|=n}}{\sf FC}(\tau)\,
\end{equation*}
induced by the canonical operadic structure on $\mathrm{End}_{\mathcal{A}}$, that is
	\begin{equation*}
		{\tilde{L}}(V \circ (W_{1} \otimes\cdots\otimes W_{p})) = V \circ ( W_{1}\otimes\cdots\otimes W_{p}),
	\end{equation*}
where $V \in {\sf FC}(p)$, $W_{i} \in {\sf FC}(n_{i})$, $1\leq i \leq p$ and the symbol $\circ$ in the right-hand side of the above equation stands for functional composition in $\mathrm{End}_{\mathcal{A}}$,
\begin{equation*}
{\sf \tilde{L}\colon \ \lpbtsfc \circ \lpbtsfc \rightarrow \lpbtsfc}.
\end{equation*}
We set FC$=\big(\mathcal{F}\mathcal{C},\tilde{L},{\rm id}_{\mathcal{A}}\big)$. Notice that with this definition, the map $\sharp$ is a morphism between the operads {\rm FS} and {\rm FC}, namely, for $A_{1},\dots,A_{p+1} \in \mathcal{A}$ and $W_1,\dots,W_p \in \mathcal{F}\mathcal{S}$
	\begin{equation*}
		\tilde{L}(\sharp (A_{1}\otimes\cdots\otimes A_{p+1}) \circ \sharp W_{1} \otimes\cdots\otimes \sharp W_{p}) = \sharp {L}(A_{1}\otimes \cdots \otimes A_{p} \circ (W_{1} \otimes\cdots\otimes W_{p})).
	\end{equation*}
\end{Remark}

\subsection{Group acting on faces-contractions}

From the algebra structure defined on $\lpbtsfc$ and the properties the morphism $\sharp \colon \lpbt^{\#}(\mathcal{A}) \to \lpbtsfc $, we will also introduce a group which plays the same role of $G(\mathcal{A})$ in~\eqref{eqn:modelspaceun} (the group where the maps $\bar{\mathbb{X}}_{s,t}$ takes value) at the level of face-contractions.

To achieve this, we first rewrite the vector space $\lpbtsfc$ in an equivalent way so that we can speak of components. Using the identification between permutations and levelled trees from Proposition~\ref{PermTree} and the intrinsic product of $\mathfrak{S}_n$, for any levelled tree $\tau$ with $\|\tau\| =n$ and $\sigma\in\mathfrak{S}_{ \|\tau\|}$ the map which sends $\sharp (A_{1} \otimes\cdots\otimes A_{n} \tau)$ to $\sharp (A_{1} \otimes\cdots\otimes A_{n} \sigma\tau)$ extends continuously to a linear map
\begin{equation*}
\phi_{\sigma}\colon \ {\sf FC}(\tau) \to {\sf FC}(\sigma \tau),
\end{equation*}
which evaluates on $m \in {\sf FC}$ as
\begin{equation*}
\phi_\sigma(m)(Y_1\otimes \cdots \otimes Y_{n}) = m(Y_{\sigma^{-1}(1)}\otimes\cdots\otimes Y_{\sigma^{-1}(n)}) = m(\sigma \cdot (Y_1\otimes\cdots\otimes Y_n)).
\end{equation*}
Each map $\phi_{\sigma}$ is continuous and has inverse given by $\phi_{\sigma^{-1}}$. Combining the action of the maps $\{\phi_{\sigma}\colon \sigma \in \mathfrak{S}_{ \|\tau\|}\}$, we introduce the map
\begin{equation*}
\phi\colon \bigoplus_{\substack{\tau \in \lpbt \\ \|\tau\|=n}} {\sf FC}(\tau) \rightarrow {\sf FC}(n) \otimes \mathbb{C}(\lpbt_n) , \qquad \phi:= \sum_{\substack{\tau \in \lpbt \\ \|\tau\|=n}}\phi_{\tau^{-1}}\otimes \tau,
\end{equation*}
and by extension of $\phi$ to $\lpbtsfc$ we obtain a continuous isomorphism
\begin{gather*}
\phi\colon \ \lpbtsfc \rightarrow \bigoplus_{n\geq 0}{\sf FC}(c_n) \otimes \mathbb{C}\lpbt_n ,
\end{gather*}
where for each $n\geq 0$, $c_n$ is the levelled tree represented by the permutation ${\rm id}_n$. For brevity, we use the notation ${\sf \lpbtsfc}(n) := {\sf \lpbtsfc}(c_n)$. Also, a generic element of the tensor product ${\sf FC}(n) \otimes \mathbb{C}\lpbt_n$ will be denoted $m \tau$ (we omit the symbol $\otimes)$ where $m \in {\sf FC}(n)$ and $\tau \in \lpbt_n$.

\begin{Definition}Let $\mathcal{X}\colon \lpbtsfc\to\lpbtsfc$ be an endomorphism of $\lpbtsfc$ and $\tau,\tau^{\prime}\in \lpbtsfc$ a couple of levelled trees. We define the \emph{components of} $\mathcal{X}$ as the set of continuous linear maps
\begin{equation*}
\big\{\mathcal{X}(\tau^{\prime},\tau),\,\tau,\tau^{\prime} \in \lpbtsfc \big\}, \qquad {\mathcal{X}}(\tau^{\prime},\tau)\colon \ {\sf FC}({\|\tau\|}) \to {\sf FC}({\|\tau^{\prime}\|}),
\end{equation*}
defined by the relation
\begin{equation}\label{operator-components}
	\mathcal{X}(m) = \sum_{\tau^{\prime},\,\tau} \phi^{-1}\big(\mathcal{X}(\tau^{\prime},\tau)\big(\phi(m)\big)\,\tau^{\prime}\big),\qquad m\in {\sf FC}(\tau).
\end{equation}
\end{Definition}

We remark that for any given family $\{\mathcal{X}(\tau^{\prime},\tau),\,\tau,\tau^{\prime} \in \lpbtsfc \}$ like above formula \eqref{operator-components} defines actually a graded endomorphism of $\lpbtsfc$. Moreover, for any given graded endomorphism $\mathcal{X}\colon \lpbtsfc\to\lpbtsfc$ the component $\mathcal{X}(\tau^{\prime},\tau)$ can be computed on any $m \in {\sf FC}({\|\tau\|}) $ as\looseness=-1
\begin{equation*}
\mathcal{X}(\tau^{\prime},\tau)(m)= (\phi( \mathcal{X}(\phi_{\tau} (m)))_{\tau'},
\end{equation*}
where $ (\cdot)_{\tau'}$ is the natural projection on the component associated with $\tau'$ as the right factor.

With the notations introduced so far, and omitting conjugation by $\phi$, fora an operator $\mathcal{X} \in {\rm End}\big(\lpbtsfc\big)$, one writes
\begin{equation*}
\mathcal{X}(m \tau) = \sum_{\tau^{\prime} \in \lpbt} \mathcal{X}(\tau^{\prime},\tau) \tau^{\prime}.
\end{equation*}
\begin{Remark}
The involution $\starltfc$ defined in the previous section induced through $\phi$ an involution on $\bigoplus_{n\geq 0}{\sf FC}(c_n) \otimes \mathbb{C}\lpbt_n$, denoted by the same symbol, one has
\begin{equation*}
\starltfc(m \tau) = \starltfc(m) \theta(\tau).
\end{equation*}
\end{Remark}
In the next definition, we introduce a specific class of operators.
\begin{Definition}
For any integer $k\geq 1$ and $Y_{1},\dots,Y_{k} \in \mathcal{A}$ we introduce the operator
\begin{equation*}
C_{Y_{1},\dots,Y_{k}} \colon \ \lpbtsfc \rightarrow \lpbtsfc
\end{equation*}
defined by the components
\begin{equation*}
C_{Y_{1},\dots,Y_{k}}(\tau^{\prime},\tau)(m) =\begin{cases} m(X_{1},\dots,X_{\|\tau^{\prime}\|}, Y_{1},\dots,Y_{k}), & {\rm if}~\|\tau\|=\|\tau^{\prime}\|+k, \ \tau'\subset \tau, \\ 0 & {\rm otherwise},
\end{cases}
\end{equation*}
where $m\in{\sf FC}(\tau)$ and $X_1,\dots,X_{\|\tau\|}$ are elements of $\mathcal{A}$.
\end{Definition}
\begin{Example}For example, by taking $m = \sharp (A_1\otimes A_2 \otimes A_3 \otimes A_4 \otimes A_5 \, \,2413)\in {\sf FC}(2413)$ and $Y_{1},Y_{2},Y_{3} \in \mathcal{A}$, from formula \eqref{operator-components} we deduce
\begin{gather*}
C_{Y_1}(m) = \sharp (A_1 \otimes A_2 \cdot Y_1\cdot A_3 \otimes A_4 \otimes A_5\,\, 213), \\
C_{Y_1,Y_2}(m) = \sharp (A_1 \otimes A_2 \cdot Y_1\cdot A_3 \otimes A_4\cdot Y_2\cdot A_5 \,\, 21), \\
C_{Y_1,Y_2,Y_3}(m) = \sharp (A_1\cdot Y_1 \cdot A_2 \cdot Y_2 \cdot A_3 \otimes A_4 \cdot Y_2 \cdot A_5 \,\, 1).
\end{gather*}
\end{Example}
The coefficient $C_{Y_{1},\dots,Y_{k}}(\tau',\tau)$ depends only the forest $\tau\backslash \tau^{\prime}$, for any $Y_1,\dots,Y_{k}$ in $\mathcal{A}$ and any integer $k\geq 1$, one has
\begin{equation}\label{eqn:cstdiagonal}
	C_{Y_1,\dots,Y_k}(\alpha',\alpha) = C_{Y_1,\dots,Y_k}(\beta',\beta),\qquad \text{if} \quad \alpha\backslash \alpha^{\prime}=\beta\backslash \beta^{\prime}, \ \alpha^{\prime}\subset\alpha, \ \beta^{\prime} \subset \beta.
\end{equation}

\begin{Definition}We denote by $U({\sf FC})$ \emph{the group of triangular algebra morphisms} of $\lpbtsfc$ with the identity on the diagonal, i.e.,
\begin{equation*}
	U({\sf FC}) = \bigg\{ \mathcal{X} \in \mathrm{End}_{\mathrm{Alg}}\big(\lpbtsfc\big)\colon (\mathcal{X}-\mathrm{id})({\sf FC}(\tau)) \subset \bigoplus_{\tau^{\prime} \subsetneq \tau} {\sf FC}(\tau^{\prime}),\, \tau \in \lpbt\bigg\}.
\end{equation*}
We denote by $\mathcal{C}$ the closure of the linear span of $C_{Y_1,\dots,Y_k}$, $Y_1,\dots,Y_k \in \mathcal{A}$ augmented with $\mathrm{id}_{\lpbtsfc}$ and we set $U_{\mathcal{C}}({\sf FC}) := U({\sf FC}) \cap \mathcal{C}$.
\end{Definition}
\begin{Remark}By construction of $U({\sf FC})$ an endomorphism of $\lpbtsfc$ belongs to $U({\sf FC})$ if and only if
\begin{gather*}
	\mathcal{X}(\tau',\tau) = 0 \quad \text{if} \ \tau^{\prime} \text{ is not a subtree of } \tau, \qquad \mathcal{X}(\tau,\tau)={\rm id}_{{\sf FC}({\|\tau\|})},
\\
\mathcal{X}(m\shuffle m')= \mathcal{X}(m)\circ \mathcal{X}( m'),\qquad \star \mathcal{X}(m)= \mathcal{X}(\star m),
\end{gather*}
for any $m,m^{\prime} \in \lpbtsfc$.
\end{Remark}
The triangular property in the definition of $U({\sf FC})$ implies also that $U({\sf FC})$ is a group. Besides, we note that $\mathcal{C}$ is a subalgebra of $\mathrm{End}\big(\lpbtsfc\big)$, because of the identity
\begin{equation*}
C_{A_{1},\dots,A_{k}} \circ C_{B_{1},\dots, B_{q}} = C_{B_{1},\dots, B_q, A_1,\dots,A_{k}}.
\end{equation*}
Therefore we obtain immediately that $U_{\mathcal{C}}({\sf FC})$ is a group,
\begin{equation*}
\mathcal{X}(f) := \mathcal{X}(\tau',\tau) \colon \ {\sf FC}({|f|}-1) \to {\sf FC}({{\sf nt}(f)}-1),
\end{equation*}
where $\tau^{\prime} \subset \tau \in \lpbt$ is any pair of levelled trees such that $f=\tau\backslash \tau^{\prime}$.

A \emph{diagonal} of an operator $\mathcal{X} \in U({\sf FC})$ is the set of entries $\mathcal{X}(\tau^{\prime},\tau)$ where $\tau^{\prime}\subset \tau$ are pairs of levelled trees with fixed forest $\tau\backslash \tau^{\prime}$. Therefore, a forest corresponds to a unique diagonal of $\mathcal{X}$. From~\eqref{eqn:cstdiagonal}, an operator $\mathcal{X}\in U_{\mathcal{C}}({\sf FC})$ may be viewed as a matrix with operators for coefficients and constant diagonals. We denote by~$\mathcal{X}(f)$ the common value of the entries of $\mathcal{X}$ on the diagonal corresponding to the forest~$f$,
\begin{equation*}
\mathcal{X}(f) := \mathcal{X}(\tau',\tau) \colon \ {\sf FC}({|f|}-1) \to {\sf FC}({{\sf nt}(f)}-1),
\end{equation*}
where $\tau^{\prime} \subset \tau \in \lpbt$ is any pair of levelled trees such that $f=\tau\backslash \tau^{\prime}$. By writing the components in terms of forests, we can exchange relations between operators in the group $U_{\mathcal{C}}({\sf FC})$.
\begin{Proposition}\label{lemma:exchange}
	For any couple ${\mathcal{X}}$, ${\mathcal{Y}}$ in $U_{\mathcal{C}}({\sf FC})$ and any pair of compatible levelled forests $f,f^{\prime} \in \lpbf$, $f\overt f^{\prime} \in \lpbf\overt \lpbf$, with the notation $\langle {\mathcal{X}}\overt {\mathcal{Y}}, f \overt f^{\prime} \rangle := {\mathcal{X}}(f) \circ {\mathcal{Y}}(f^{\prime}),
 $ one has
	\begin{equation*}
	\langle {\mathcal{X}}\overt {\mathcal{Y}}, f \overt f^{\prime} \rangle = \langle {\mathcal{Y}}\overt {\mathcal{X}}, {\sf K}(f\overt f^{\prime}) \rangle, \qquad f\overt f^{\prime} \in \lpbf\overt\lpbf.
	\end{equation*}
\end{Proposition}
\begin{proof}Let $A_{1},{\dots},A_{|f^{\prime}|} \!\in\! \mathcal{A}$ and call $\sigma$ (resp.~$\sigma^{\prime}$) the permutation associated with $f^{\flat}$ (resp.~$(f^{\prime})^{\flat}$). We use the notation $cb_{n}$ for the right-comb tree associated with the identity permutation $\mathrm{id}_{n}$. Next, define $s$ the permutation in $\mathfrak{S}_{\|f\|+\|f^{\prime}\|}$ by
	\begin{itemize}\itemsep=0pt
		\item $s_{f\overt f^{\prime}}(k) = i$, if the $k^{\rm th}$ face of $cb_{{\sf nt}(f)})\# f \# f^{\prime}$ (reading the faces from left to right) is the~$i^{\rm th}$ face of $f$,
		\item $s_{f\overt f^{\prime}}(k)=\|f\|+i$ if the $k^{\rm th}$ face of $f \# f^{\prime}$ is the $i^{\rm th}$ face of $f^{\prime}$,
		\item $s_{f\overt f^{\prime}}(k) = \|f\|+\|f^{\prime}\|+i$ if the $k^{\rm th}$ face is the $i^{\rm th}$ face of $cb_{{\sf nt}(f)}$.
	\end{itemize}
	With ${\sf{K}}(f\overt f^{\prime}) = f_{(1)}^{\prime} \overt f_{(1)}$, with $\|f^{\prime}_{1}\|=\|f^{\prime}\}$ and $\|f\|=\|f_{(1)}\|$ notice that
\begin{equation*}
		f^{\prime}_{(1)}{}^{\flat} = f^{\prime}{}^{\flat}, \qquad f_{(1)}{}^{\flat} = f^{\flat}, \qquad s_{{\sf{K}}(f \overt f^{\prime})} = s_{f\overt f^{\prime}}.
\end{equation*}
	Note that $(({\mathcal{X}}\overt {\mathcal{Y}}), f \overt f^{\prime})$ is non-zero only of pair of forests $f\overt f^{\prime}$ with $\|f\|=p$ and $\|f^{\prime}\|=q$. Pick two such forests $f$, $f^{\prime}$.
	Pick $U_{1},\dots,U_{{\sf nt}(f)-1} \in \mathcal{A}$. Pick $\mathcal{X}=V_{X_{1},\dots,X_{p}}$ and $\mathcal{Y}=V_{Y_1,\dots,Y_q}$ two operators in $\mathcal{C}$ and put $Z=(U_{1},\dots, U_{{\sf nt}(f)},X_1,\dots,X_p,Y_1,\dots,Y_q)$. Therefore one has
	\begin{gather*}
		 \langle {\mathcal{X}}\otimes {\mathcal{Y}}, f \overt f^{\prime} \rangle(\sharp(A_{1}\otimes\cdots\otimes A_{|f^{\prime}|}))(U_{1},\dots,U_{{\sf nt}(f)-1})\\
\qquad{}= A_{1}\cdot Z_{{s_{f\overt f^{\prime}}^{-1}(1)}}\otimes\cdots\otimes Z_{{s_{f\overt f^{\prime}}^{-1}}(\|f\|+\|f^{\prime}\|)}\cdot A_{|f^{\prime}|}\\
\qquad{} =\langle {\mathcal{Y}}_{s,t}\otimes \mathcal{X}_{s,t}, {\sf K}(f\overt f^{\prime}) \rangle(\sharp(A_{1}\otimes\cdots\otimes A_{|f^{\prime}|}))(U_{1},\dots,U_{{\sf nt}(f)-1}).\tag*{\qed}
	\end{gather*}\renewcommand{\qed}{}
\end{proof}

\begin{Proposition}	\label{prop:involutionU}The involution $\Theta \colon {\rm End}\big(\lpbtsfc\big) \to {\rm End}\big(\lpbtsfc\big)$ defined on any $\mathcal{X} \in {\rm End}(\lpbtsfc)$ by
	\begin{equation*}
		\Theta(\mathcal{X}) := \starltfc \circ \mathcal{X} \circ \starltfc, \qquad \mathcal{X}\in {\rm End}\big(\lpbtsfc\big),
	\end{equation*}
	restricts to an algebra morphism on $U_{\mathcal{C}}({\sf FC})$.
\end{Proposition}
\begin{proof}Pick $\mathcal{X},\mathcal{Y}\in U_{\mathcal{C}}({\sf FC})$, we have to show that
\begin{equation*}
	\Theta(\mathcal{X})\circ\Theta(\mathcal{Y}) = \Theta(\mathcal{X}\circ \mathcal{Y}).
\end{equation*}
	Let $\tau \in \lpbtsfc$ and $m \in {\sf FC}(||\tau||)$, one has first
	\begin{align*}
		(\mathcal{X}\circ \mathcal{Y}) (\starltfc(m) \theta(\tau)) & = \sum_{\tau^{\prime}\subset \raisebox{-1pt}{$\theta$}(\tau)}(\mathcal{X}\circ \mathcal{Y})(\theta(\tau)\backslash \tau^{\prime})(\star(m)) \tau^{\prime} \\
&= \sum_{\tau^{\prime}\subset \tau}(\mathcal{X}\circ \mathcal{Y},\theta(\tau\backslash \tau^{\prime}))(\starltfc(m)) {\theta}(\tau^{\prime}) \\
		&=\sum_{\substack{\tau^{\prime}\subset \tau
	 \\ f^{\prime} \subset \raisebox{-1pt}{$\theta$}(\tau\backslash \tau^{\prime})}} \langle\mathcal{X}\overt\mathcal{Y}, f^{\prime}\overt \theta(\tau\backslash \tau^{\prime}) \backslash f^{\prime}\rangle(\starltfc(m)) \theta(\tau^{\prime}) \\
	 & = \sum_{\substack{\tau^{\prime}\subset \tau \\ f^{\prime} \subset \tau\backslash \tau^{\prime}}} \langle \mathcal{X}\overt\mathcal{Y}, \theta(f^{\prime})\overt \theta(\tau\backslash \tau^{\prime} \backslash f^{\prime})\rangle (\starltfc(m)) \theta(\tau^{\prime}),
	\end{align*}
	where the last equality follows from the simple observation that $\theta(f\backslash f^{\prime}) = \theta (f) \backslash \theta (f^{\prime})$ for any pair of forests $f^{\prime} \subset f$. We deduce that
	\begin{equation*}
		\Theta(\mathcal{X}\circ \mathcal{Y})(m \tau) = \sum_{\substack{\tau^{\prime}\subset \tau \\ f^{\prime} \subset \tau\backslash \tau^{\prime}}} \starltfc\big(\langle\mathcal{X}\overt\mathcal{Y}, \theta(f^{\prime})\overt \theta(\tau\backslash \tau^{\prime} \backslash f^{\prime})\rangle(\star(m))\big) (\tau^{\prime}).
	\end{equation*}
	For any pair of forests $f^{\prime}\subset f$,
	\begin{align*}
		\starltfc(\langle\mathcal{X}\overt\mathcal{Y}, \theta(f^{\prime})\overt \theta( f\backslash f^{\prime}))(\star(m))\rangle & = \starltfc\big(\big(\mathcal{X}(\theta(f))\circ \mathcal{Y}(\theta(f\backslash f^{\prime}))\big)(\star(m))\big) \\
	 & =\langle \Theta(\mathcal{X}) \overt \Theta(\mathcal{Y}), f\overt f^{\prime})(m)\rangle.
	\end{align*}
	Hence,
	\begin{align*}
	 \Theta(\mathcal{X}\circ \mathcal{Y})(m \tau) = \sum_{\substack{\tau^{\prime} \subset \tau \\ f^{\prime} \subset \tau \backslash \tau^{\prime}}} \langle \Theta(\mathcal{X}) \overt \Theta(\mathcal{Y}), f\overt f^{\prime}\rangle(m) \tau^{\prime} = \Theta(\mathcal{X}) \circ \Theta(\mathcal{Y}) (m \tau).\tag*{\qed}
	\end{align*}\renewcommand{\qed}{}
\end{proof}

We now further restrict the group that will support signatures of smooth valued paths. It is defined by a set of equations on components of an operator in $U_{\mathcal{C}}$ indexed by levelled forests. To specify these equations, we resort to operations on levelled forests defined in Section~\ref{sec:levelledforests}.

Pick a integer $n\geq 1$. Let $f = (\sigma,c)$, $\sigma\in\mathfrak{S}_n$, $c\vDash_0 n$ be a levelled forest with $n$ generations and ${\sf nt}(f)$ trees.
For any subset $ I \subset [{\sf nt}(f)-1]$ of integers we denote by $f^I$ the levelled forest obtained by gluing the trees $f_i$ in $f$ at positions $i\in I$ along their external edges,
\begin{equation*}
f^{I} = (\sigma, (c_1,\dots,c_{i_1} + c_{i_1+1},\dots, c_{i_k} + c_{i_{k+1}},\dots, c_{{\sf nt}(f)+1})).
\end{equation*}
We set also
\begin{equation*}
\ell^{f}_{j} = \sum_{i=1}^{j} (c_{i}+1),\qquad 1 \leq j \leq {\sf nt}(f)-1.
\end{equation*}
The integers $\ell^{f}_j$, $1 \leq j \leq {\sf nt}(f)-1$ index spaces between consecutive sparse quasi-binary trees in the levelled forest~$f$, which can also be considered as faces of $f$ (in addition to the ones delimited by two consecutive leaves of a tree in~$f$).

Recall also that elements of $\mathcal{A}$ are considered as face-contractions operators with $0$ input. Pick $m {\sf FC}(p)$ a face-contractions operator with arity $p$, $m_{1},\dots, m_{q}$, $ 1 \leq g \leq p$ other face-contractions operators and a sequence of integers $1 \leq i_{1}<\cdots<i_{q} \leq p$. In the definition below, we denote by
\begin{equation*}
	m\circ_{i_{1},\dots,i_{q}} m_{1}\otimes\cdots\otimes m_{q}
\end{equation*}
the operator obtained by composing $m_{j}$ with the $i_{j}^{\rm th}$ input of $m$, the remaining inputs are filled with the identity.

\begin{Definition}\label{def:grouproughpaths}We define $G({\sf FC})\subset U_{\mathcal{C}}({\sf FC})$ as the subset comprising all operators $\mathcal{X}\in U_{\mathcal{C}}({\sf FC})$ satisfying for any quadruple
	\begin{itemize}\itemsep=0pt
	\item $f$ a levelled forest,
	\item an integer $1\leq q \leq {\sf nt}(f)-1$,
	\item a face-contractions operator $m \in \lpbtsfc$ with $|f|-1$ inputs,
	\item a sequence of integers $1\leq i_1 <\cdots < i_q < {\sf nt}(f)$,
	\end{itemize}
	the relation
	\begin{equation}\label{eqn:simplicial}
		{\mathcal{X}}(f)(m)\circ_{i_1,\dots,i_q}(X_{1},\dots,X_{q}) = {\mathcal{X}}\big(f^{I}\big)(m \circ_{\ell^{f}_{i_{1}},\dots,\ell^{f}_{i_{q}}}X_{1},\dots,X_q),
	\end{equation}
 where $X_1,\dots,X_q \in \mathcal{A}$.
	In addition, we denote by $G_{\star}({\sf FC})$ the subset of self-adjoint operators in $G({\sf FC})$ for the involution $\Theta$ defined in Proposition~\ref{prop:involutionU},
\begin{equation*}
	G_{\star}({\sf FC}) := \{ \mathcal{X} \in G({\sf FC}) \colon \Theta(\mathcal{X}) = \mathcal{X}\}.
\end{equation*}
\end{Definition}
\begin{Example} We give a simple example of the condition $\eqref{eqn:simplicial}$. Consider $f$ the forest represented by $(132, (2,1))$, a faces contraction operator $m$ with four inputs. Choose $i_1=1$ (this is the only possibility). Equation \eqref{eqn:simplicial} is equivalent to
\begin{equation*}
{\mathcal{X}}(f)(m)\circ_{1}(X_{1}) = {\mathcal{X}}(132)(m \circ_{3}X_{1}).
\end{equation*}
\end{Example}
\begin{Remark}\label{rk:simplicialop}On can restricts $m$ in equation \eqref{eqn:simplicial} to be an operator of
the form $\sharp (A_0\otimes\cdots\otimes A_p \tau )$ for a certain word $A_0,\dots,A_p \in \mathcal{A}$ and levelled tree $\tau$. If there exists $\bar{\mathbb{X}}$ such that $\sharp \circ \bar{\mathbb{X}} = \mathcal{X} \circ \sharp$ and
\begin{equation*}
\bar{\mathbb{X}}(f)(A_1 \otimes \cdots \otimes A_{|f|}) = {\rm Op}( [X_{s,t}^{\sigma}]_c)(A_0,\dots,A_f)
\end{equation*}
then \eqref{eqn:simplicial} is automatically satisfied since by evaluation of both sides of \eqref{eqn:simplicial} of $Y_1,\dots,Y_{{\sf nt}(f)-1-q}$ one obtains
\begin{equation*}
{\rm Op}\big(X^{\sigma}_{s,t}\big)(Z_1,\dots,Z_{\|f\|}),
\end{equation*}
where
\begin{enumerate}\itemsep=0pt
\item[(1)] $Z_{\ell^f_{i_j}} = A_{\ell^{f}_{i_j}}\cdot X_j \cdot A_{\ell^{f}_{i_j}+1}$,
\item[(2)] $Z_{\ell^f_{i_j+s}} = A_{\ell^{f}_{i_j+s}}\cdot Y_{i_j+s-j} \cdot A_{\ell^{f}_{i_j+s}+1}$, $1 < s < i_{j+1}-i_j$,
\item[(3)] $Z_{\ell^f_j+t} = A_{\ell^f_j+t-j}$, $1 < t < \ell^f_j -\ell^f_{j+1}$.
\end{enumerate}
\end{Remark}
\begin{Proposition}\label{prop:groupfacescontraction}The sets $G({\sf FC})$ and $G_{\star}({\sf FC})$ are sub-groups of $U_{\mathcal{C}}({\sf FC})$.
\end{Proposition}
\begin{proof}Pick two endomorphisms $\mathcal{X}$ and $\mathcal{Y}$ in $G({\sf FC})$. Pick $f$ a levelled forest and $I =\{i_1<\cdots < i_q \} \subset \llbracket 1, {\sf nt}(f)-1\rrbracket$. We prove that $\mathcal{X}\circ\mathcal{Y} \in G({\sf FC})$. We already know that $\mathcal{X}\circ\mathcal{Y} \in U_{\mathcal{C}}({\sf FC})$, hence it will be sufficient to prove $\eqref{eqn:simplicial}$ for $\mathcal{X}\circ\mathcal{Y}$. Pick $m\in \lpbtsfc$ as in Definition \ref{def:grouproughpaths}, one has
	\begin{align*}
	((\mathcal{X}\circ\mathcal{Y})(f))(m)\circ_{i_1,\dots,i_p}(A_{I}) & = \sum_{f^{\prime}\subset f } \mathcal{X}(f^{\prime})\big(\mathcal{Y}(f\setminus f^{\prime})(m)\big)\circ_{i_1,\dots,i_p}(A_I)\\
	& =\sum_{f^{\prime}\subset f } \mathcal{X}\big(f^{\prime}{}^{I}\big)\big(\mathcal{Y}(f\setminus f^{\prime})(m)\circ_{\ell^{f^{\prime}}_{i_1},\dots,\ell^{f^{\prime}}_{i_p}} A_{I}\big) \\& =\sum_{f^{\prime} \subset f} \mathcal{X}\big(f^{\prime}{}^{I}\big)
\big(\mathcal{Y}\big((f\setminus f^{\prime})^{\ell^{f^{\prime}}_{i_1},\dots,\ell^{f^{\prime}}_{i_p}}\big)\big)\big(m\circ_{\ell^{f}_{\ell_{i_1}^{f^{\prime}}},\dots,\ell^{f}_{\ell^{f^{\prime}}_{i_p}}}(A_{I})\big).
	\end{align*}
	Owing to associativity of $\circ$, we have $ \ell^{f}_{\ell_{i_1}^{f^{\prime}}},\dots,\ell^{f}_{\ell^{f^{\prime}}_{i_p}}= \ell_{i_1}^{f^{\prime}\overt f},\dots,\ell_{i_p}^{f^{\prime}\overt f}$. The statement follows by noticing that
	\begin{equation*}
		\big\{\big(f^{\prime}{}^{I},(f\setminus f^{\prime})^{\ell^{f^{\prime}}_{i_1},\dots,\ell^{f^{\prime}}_{i_p}}\big), f^{\prime}\subset f\big\} = \big\{ (f^{\prime},f\setminus f^{\prime}), f^{\prime}\subset f^I\big\}.\tag*{\qed}
	\end{equation*}\renewcommand{\qed}{}
\end{proof}

\subsection{The noncommutative signature of a path}
\begin{figure}[ht]	\centering
	\includegraphics{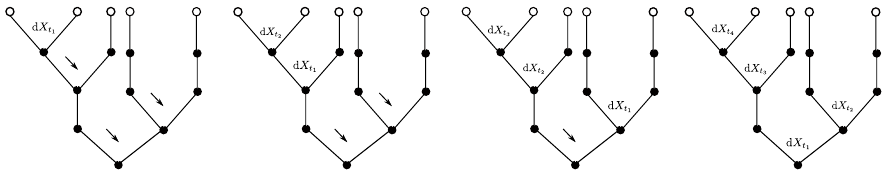}
	\caption{face-contractions of a levelled forest.}	\label{fig:faces_contraction}
\end{figure}

We are now ready to state the main Definition of the notion of the signature of a smooth path that is adapted to the class of equations \eqref{eqn:ncequations}.
\begin{Definition}\label{def:endofacescontraction}
	Pick $X\colon [0,1]\rightarrow \mathcal{A}$ a smooth path. Let $0<s < t<1$ be two times. We define a triangular endomorphism,
\begin{equation*}
{\mathcal{X}}_{s,t} \colon \ \lpbtsfc \rightarrow \lpbtsfc
\end{equation*}
defined, for $m \in {\sf FC}(\tau)$, a pair of trees $\tau^{\prime}\subset \tau \in \lpbt$, $Y_1,\dots,Y_{\|f\|} \in \mathcal{A}$, by
\begin{gather*}
{\mathcal{X}}_{s,t}(f)(m)(Y_{1},\dots,Y_{\|f\|}) = \int_{\Delta^{\|f\|}_{s,t}}\!\! m\big(\sigma^{-1}\!\cdot\! \big(Y_{1}\otimes \cdots\otimes Y_{\|{\sf nt}(f)-1\|}\otimes \mathrm{d}X_{t_{1}}\otimes \cdots \otimes \mathrm{d}X_{t_{\|f\|}}\big)\big)
\end{gather*}
for a levelled forest $f=(\sigma,c)$.
See Figure~\ref{fig:faces_contraction} for a picture representing the action of ${\mathcal{X}}_{s,t}$.
\end{Definition}

Recall that $\bar{\mathbb{X}}_{s,t}$ is defined in equation \eqref{eqn:barncrp}.
\begin{Theorem}\label{cor:algebrastar}
	For any smooth path $X$ the operators ${\mathcal{X}}_{s,t}$ satisfy the following properties:
	\begin{enumerate}\itemsep=0pt
	\item[$(1)$] ${\mathcal{X}}_{s,t}\in G({\sf FC})$,
	\item[$(2)$] ${\mathcal{X}}_{s,t} = {\mathcal{X}}_{u,t} \circ {\mathcal{X}}_{s,u}$,
	\item[$(3)$] $\sharp \circ \bar{\mathbb{X}}_{s,t} = {\mathcal{X}}_{s,t} \circ \sharp$.
	\end{enumerate}
Besides, if $X$ is a self-adjoint trajectory then $\mathcal{X}_{s,t}$ is a morphism of $\star$ algebras.
\end{Theorem}
\begin{proof}Item (2) is implied by item~(3) and the Chen relation for $\mathcal{X}$. By Remark~\ref{rk:simplicialop}, (1) is implied by item~(3). Item~(3) follows by direct computations.
\end{proof}

The above proposition implies that the endomorphism $\mathcal{X}_{s,t}$ is characterized by its values on the levelled trees. Partial contractions appear thus as technical proxies required to write in a compact form the Chen relation for the operators $\mathcal{X}_{s,t}$. However, this gives no additional information on the small-scale behaviour of $X$. This is compliant with the simple observation that expansion of a solution of an equation in the class~\eqref{eqn:controlleddiff} does only involve full contraction,
\begin{equation*}
L^{k}\colon \ \mathcal{A}^{\otimes k}\ni (A_{1}\otimes\cdots\otimes A_{k}) \mapsto L_{A_{1},\dots,A_{k}}
\end{equation*}
is well-defined and continuous. We call $\mathcal{C}$ the closure for the operator norm of the direct sum of the ranges of the operators $L^{k}$.

\appendix

\section{Appendix}\label{appendix}

We recall some definitions from the theory of operads and more generally, we underline here the categorical notions we use in this work. The reader will find below, among other things, definitons of collections, operads, bi-collections and PROs. All of concepts are standard in the algebraic literature, see, e.g., the monographs \cite{aguiar2010monoidal, loday2012algebraic}, but not very known among non-algebraists. Hence the need of this small appendix. For further details we refer to \cite{bultel2016combinatorial,vallette2007koszul}.

At the base of these definitions above lies the concept of \emph{monoidal category}. In loose words, it is a category $ \mathcal{C}= (\text{Ob}(\mathcal{C}), \text{Mor}(\mathcal{C})))$ equipped with an operation $\sbt$ and a unity element $I\in \text{Ob}(\mathcal{C})$. The operation $\sbt$ associates to any couple of objects $A,B\in \text{Ob}(\mathcal{C})$ an object $ A\sbt B\in \text{Ob}(\mathcal{C})$ and to any couple of morphisms $f\colon A\to A'$, $g\colon B\to B'$ a morphism $f \sbt g\colon A \sbt B\to A' \sbt B'$ in a~functorial way. In order that $(\mathcal{C}, \sbt, I) $ is a monoidal category, the operation $\sbt$ must satisfy two main properties, which emulate the tensor product operation on finite-dimensional vector spaces:
\begin{enumerate}\itemsep=0pt
\item (Associativity constraints) for any triple of objects $A,B,C\in \text{Ob}(\mathcal{C})$ one has that the object $(A \sbt B) \sbt C$ is isomorphic to $A \sbt (B \sbt C)$ in a \emph{functorial way}, that is there exists a~natural isomorphism between the two functors $\sbt \circ ({\rm id} \times \sbt)$ and $\sbt \circ (\sbt \times {\rm id} )$.
\item (Unitality constraints) for any object $ A \in \text{Ob}(\mathcal{C})$ the objects $A \sbt I$ and $I \sbt A$ are (naturaly) isomorphic to $A$.
\end{enumerate}
The prototypical example is the category of finite-dimensional vector spaces with monoidal product given by the tensor product of vector spaces $\otimes$. Another example is the category Set, the category of all sets with functions between sets as morphisms, with monoidal product given by the cartesian product of sets.\footnote{This monoidal category is particular in the sense that the monoidal product coincides with the categorical product. Such categories are called cartesian monoidal.}
Of interest in the present work is the $2$-monoidal category of collections and bicollections that we now define.

A monoid in a monoidal category is a categorical abstraction of a binary product on a set.
\begin{Definition}[monoid]A \emph{monoid} in a monoidal category $(\mathcal{C},\bullet, I)$ is a triple $(C,\rho,\eta)$ with $C\in \mathrm{Ob}(\mathcal{C})$, $\rho\colon C\bullet C \to C$, $\eta\colon I\to C$ meeting the constraints\samepage
\begin{enumerate}\itemsep=0pt
 \item[(1)] $\rho \circ (\rho \bullet \mathrm{id}) = \rho \circ (\mathrm{id} \bullet \rho)$,
 \item[(2)] $\rho \circ (\eta \bullet \mathrm{id}) = \mathrm{id} $.
\end{enumerate}
\end{Definition}
\begin{Definition}[comonoid]
A \emph{comonoid} in a monoidal category $(\mathcal{C}, \bullet, I)$ is a triple $(C,\Delta,\varepsilon)$ with $C\in \mathrm{Ob}(\mathcal{C})$, $\Delta\colon C \rightarrow C\bullet C$, $\varepsilon\colon C\to I$ meeting the constraints:
\begin{enumerate}\itemsep=0pt
 \item[(1)] $(\Delta \bullet \mathrm{id}) \circ \Delta = (\mathrm{id} \bullet \Delta )\circ \Delta$,
 \item[(2)] $(\varepsilon \bullet \mathrm{id}) \circ \Delta = (\mathrm{id}\bullet \varepsilon )\circ \Delta$.
\end{enumerate}
\end{Definition}
\newcommand\Coll{\mathsf{Coll}}
\begin{Definition}\label{def:collection}
We call a (reduced) \emph{collection} $P$ a sequence of complex vector spaces\footnote{The original definition involves vector spaces over a generic field but we consider only complex vector spaces, in accordance with the structures presented so far.} \linebreak $\{P(n)\}_{n\geq 1}$. A~morphism between two collections $P$, $Q$ is a sequence of linear maps $\{\phi(n)\}_{n\geq 1}$ with $\phi(n)\colon P(n) \allowbreak \rightarrow Q(n)$, $n\geq 0$. For any couple of morphisms between collections we define the composition of morphisms by composing each component. We denote the category of collections by $\Coll$.
\end{Definition}
The category $\Coll$ has a natural monoidal structure $\odot$ over it: for any couple of collections $P$ and $Q$ and morphisms $f$, $g$ we define
\begin{gather*}
	 (P \odot Q)(n): = \bigoplus_{\displaystyle{\substack{k \geq 1 \\n_{1}+\cdots+n_{k}=n } }}
	P(k) \otimes Q(n_{1}) \otimes \cdots \otimes Q(n_{k}),\\
	 (f \odot g)(n):= \bigoplus_{\displaystyle{\substack{k \geq 1 \\n_{1}+\cdots+n_{k}=n } }}
	f(k) \otimes g(n_{1}) \otimes \cdots \otimes g(n_{k}) .
\end{gather*}
Denoting by $\mathbb{C}_{\odot}$ the collection
\[\mathbb{C}_{\odot}= \begin{cases}
\mathbb{C} & \text{if}\ n = 1,\\
 0 & \text{otherwise},\end{cases}\]
it is straightforward to check that the triple $(\Coll,\odot, \mathbb{C}_{\odot})$ is a monoidal category. If the vectors spaces of the collections $P$ and $Q$ above are Banach algebras, then we might use in place of the algebraic tensor product $\otimes$ the projective one $\hat{\otimes}$. An operad is a monoid in the monoidal category $(\Coll,\odot, \mathbb{C}_{\odot})$:
\begin{Definition}\label{def:operad}
A non-symmetric \emph{operad} (or simply an operad) is a monoid in the monoidal category $(\Coll,\odot, \mathbb{C}_{\odot})$, i.e., a triple $(P, \rho,\eta_{P})$ of the following objects
\[
P \in \mathrm{Ob}(\Coll), \qquad \rho\colon \ P \odot P \rightarrow P, \qquad \eta_{P}\colon \ \mathbb{C}_{\odot} \rightarrow {P},
\]
satisfying the properties $(\rho \odot \textrm{id}_{P}) \circ \rho = (\textrm{id}_{P} \odot \rho) \circ \rho$ and $(\eta_{P} \odot \textrm{id}_{P}) \circ \rho = (\textrm{id}_{P} \odot \eta_{P}) \circ \rho = \textrm{id}_{P}$.
\end{Definition}

We keep the notation $\odot$ for the monoidal operation. It is common in the literature to denote the morphism $\rho$ by $\circ$, i.e., for every $k\geq 1$, $p\in P(k)$ and $q_{i}\in Q(n_{i})$ for $i=1, \dots, k$,
\[
p \circ (q_{1}\otimes \cdots \otimes q_{n} ):= \rho(n_1+ \cdots +n_k)(p \otimes q_{1}\cdots \otimes q_{k}).
\]
Moreover, for any $ 1 \leq i \leq k $ and $q_{i}\in Q(n_{i})$ we use also the notation $\circ_{i}$ to denote partial composition
\[
p\circ_{i}q := p \circ \big(\eta_P(1)(1)^{\otimes i-1} \otimes q \otimes \eta_P(1)(1)^{\otimes k-i}\big),
\]
where $\eta_P(1)\colon \mathbb{C}\to P(1)$. Since the maps $ \rho(n)_{n\geq 1}$ carry multiple inputs and give back one output, it is common in the literature to call them many-to-one operators. A classical example to understand this definition is given by the endomorphism operad $ {\rm End}_{V}$ of a complex vector space $V$ with elements
\[	{\rm End}_{V}(n)={\rm Hom}_{{\rm Vect}_{\mathbb{C}}}\big(V^{\otimes n},V\big),\quad n\geq 1, \qquad \eta_{{\rm End}_{V}}(1)(1)= \mathrm{id}_{V},\]
and operadic composition
\[f \circ (f_{1}\otimes \cdots \otimes f_{k})(v_1\otimes \cdots \otimes v_{n}):=f\big(f_1(v_1\otimes \cdots \otimes v_{n_1}), \dots, f_k(v_{n-n_k+1}\otimes \cdots \otimes v_{n})\big),\]
where $n=n_1+\cdots+ n_k$. In fact, it is possible to generalise the notion of an operad to model composition between many-to-many operators, that is operators with multiple in- and outputs. This leads us to define the category of bicollections.
\begin{Definition}\label{def:bigradedcollection}
We call a \emph{bicollection} a two parameters family of complex vector spaces
\[
P = \{P(n,m)\}_{n,m \geq 0}.
\]
A morphism between two bicollections $P$, $Q$ is a sequence of linear maps $\{\phi(n,m)\}_{n,m\geq 0}$ with $\phi(n,m)\colon P(n,m) \rightarrow Q(n,m)$. For any couple of morphisms between bicollections we define the composition of morphisms by composing each component. We denote the category of bicollections by~$\Coll_2$.
\end{Definition}
The category of bicollections is endowed with two compatible monoidal structures.
\begin{Definition}
For any couple of bicollections $P$ and $Q$ and morphisms $f$, $g$ we define the \emph{horizontal tensor product} $\ominus$ as follows:
\begin{gather*}
 (P\ominus Q)(n,m) := \displaystyle\bigoplus_{\substack{n_{1}+n_{2}=n \\ m_{1}+m_{2}=m}} P(n_{1},m_{1}) \otimes Q(n_{2},m_{2}),\nonumber\\
(f \ominus g)(n, m) := \displaystyle\bigoplus_{\substack{n_{1}+n_{2}=n \\ m_{1}+m_{2}=m}} f(n_{1},m_{1}) \otimes g(n_{2},m_{2})
\end{gather*}
\newcommand\overtical{\mathbin{\rotatebox[origin=c]{90}{$\ominus$}}}
together with the horizontal unity
\[\pmb{\mathbb{C}}_{\ominus}=\pmb{\mathbb{C}}_{\ominus}(m,n)= \begin{cases}
\mathbb{C} & \text{if} \ n=m=0,\\
 0 & \text{otherwise}.\end{cases}\]
We define also the \emph{vertical tensor product} $\overt$
\begin{gather*}
(P\overt Q)(m,n) := \bigoplus_{k= 0}^{+ \infty} P(m,k) \otimes Q(k,n),\nonumber \\
(f \overt g)(n,m):= \bigoplus_{k= 0}^{+ \infty} f(m,k) \otimes g(k,n)
\end{gather*}
together with the vertical unity
\[\pmb{\mathbb{C}}_{\overt}=\pmb{\mathbb{C}}_{\overt}(m,n)= \begin{cases}
\mathbb{C} & \text{if}\ n=m,\\
 0 & \text{otherwise}.\end{cases}\]
We refer to the triple $(\Coll_2, \ominus ,\pmb{\mathbb{C}}_{\ominus})$ and $(\Coll_2, \overt ,\pmb{\mathbb{C}}_{\overt})$ respectively as the category of \emph{horizontal bicollections} and the \emph{vertical bicollections}.
\end{Definition}

\begin{Lemma}
$(\Coll_2, \ominus ,\pmb{\mathbb{C}}_{\ominus})$ and $(\Coll_2, \overt ,\pmb{\mathbb{C}}_{\overt})$ are monoidal categories.
\end{Lemma}
\begin{proof}This is simple computations, based on the fact that $({\rm Vect}_{\mathbb{C}},\otimes,\C)$ is monoidal.
\end{proof}
\begin{Remark}We point at some core differences and similarities between the tensor product of vector spaces, and the two tensor products we defined on bicollections.

If $V$ and $W$ are two vector spaces, there exists an isomorphism of vector spaces $S_{V,W}\colon V\otimes W \to W\otimes V$. The set $\big\{S^{\otimes}_{V,W},V,W\in {\rm Vect}_{\C} \big\}$ defines a natural transformation, called a symmetry constraint. for any $V$ and $W$ are bicollections The horizontal tensor product is symmetric with symmetry constraint given by $S^{\otimes}_{V,W}$. Nevertheless,
the vertical tensor product $\overt$ does not have such symmetry constraints, though we constructed such one but for the monoid generated by the bicollection $\mathcal{L}\mathcal{F}$ in~\eqref{levelled_forests}.
Another important property related to the vertical and horizontal tensor product is the closedness. A category $\mathcal C$ is said \emph{closed} if for all objects $A,B \in \mathrm{Ob}(\mathcal{C})$ the set of morphisms $\Hom_{\mathcal C}(A,B)$
is an object of~$\mathcal C$.
A~monoidal category $(\mathcal{C}, \bullet, I)$ is said a \emph{closed monoidal category} if it is closed and the following compatibility holds:
for all objects $A$, $B$, $C$ in $\mathrm{Ob}(\mathcal{C})$
\begin{align*}
 \Hom_{\mathcal C}(A,\hom_{\mathcal C}(B, C))
 \cong
 \Hom_{\mathcal C}(A \bullet B,C),
\end{align*}
with the isomorphism being natural in all three arguments. The category of finite-dimensional vector spaces with the usual tensor product is closed monoidal, owing to the fact that the set of linear maps between vector spaces is again a vector space and then using usual identification of bilinear maps with linear maps on the tensor product. Now, neither $(\Coll_2, \ominus ,\pmb{\mathbb{C}}_{\ominus})$ nor $(\Coll_2, \overt,\pmb{\mathbb{C}}_{\overt})$ are closed monoidal. Indeed, they are not even closed, since there is no canonical bigrading on the set of morphisms.
\end{Remark}

There exists a functor from the category of collections to the category of bicollections, that is the free horizontal monoid functor $T\colon \Coll \rightarrow \Coll_{2}$, adjoint to the forgetful functor associating to a monoid $(P,\gamma, \eta)$ for the horizontal tensor product $\ominus$ the collection $(P(1,n))_{n\geq 1}$.
\begin{Definition}
Let $P = (P_{n})_{n\geq 1}$ be a collection, we define the bicollection $T(P)$ by
\begin{equation*}
T(P)(m,n) =\bigoplus_{\substack{k_1+\cdots + k_m = n}} P_{k_{1}}\otimes\cdots\otimes P_{k_{m}} ,
\end{equation*}
when $n\geq 1$ and $m \geq 1$ and the condition $k_1+\cdots + k_m = n$ is satisfied for some integer $k_1, \dots, k_m\geq 1$. Moreover, we set $T(P)(0,0)=\mathbb{C} $ and $T(P)(m,n)=0$ otherwise.
\end{Definition}

\begin{Proposition}\label{lem:lax}
Let $C_{i}$, $1 \leq i \leq 4$ be four bicollections, then there exists an explicit morphism
\[
R_{C_{1},C_{2},C_{3},C_{4}}\colon \ (C_{1} \overt C_{2} ) \ominus (C_{3}\overt C_{4}) \rightarrow (C_{1} \ominus C_{3} ) \overt (C_{2}\ominus C_{4} ).
\]
We call $R_{C_{1},C_{2},C_{3},C_{4}} $ the \emph{exchange law}. Besides, if the bicollections $C_{2}$ and $C_4$ are equal and in the image of $\mathcal{F}$, one has
\begin{equation*}
(C_{1} \overt T(C)) \ominus (C_{2}\overt T(C) ) \simeq (C_{1} \ominus C_{2} ) \overt T(C).
\end{equation*}
\end{Proposition}

The family of morphisms $\{R_{C_{1},C_{2},C_{3},C_{4}},\, C_{i} \in \Coll_{2}\}$ defines a \emph{natural transformation} (which is, in general, not an isomorphism) between the functors $\ominus \circ \overt \times \overt$ and $\overt \circ \ominus \times \ominus$. In particular, for any quadruplet of morphisms $f_{i}\colon C_{i}\rightarrow D_{i}$, $1 \leq i \leq 4$, one has the commutative diagram (see Figure~\ref{fig:naturalR}).
We denote by $\mathrm{Alg}_{\ominus}$ (resp.\ $\mathrm{CoAlg}_{\ominus}$) the category of all monoids (resp.\ comonoids) in $(\Coll_2, \ominus ,\pmb{\mathbb{C}}_{\ominus})$, Alg$_{\overt}$ (resp.~CoAlg$_{\overt}$ the category of monoids (resp.\ comonoids) in $(\Coll_2, \ominus ,\pmb{\mathbb{C}}_{\overt})$.

\begin{figure}[!ht]	\centering
	\begin{tikzcd}
		(C_{1} \overt C_{2}) \ominus (C_{3} \overt C_{4}) \arrow{r}[yshift=2ex]{(f_{1}\overt f_{2})\ominus (f_{3}\overt f_{4})}\arrow{d}{R_{C_{1},C_{2},C_{3},C_{4}}} & (D_{1} \overt D_{2}) \ominus (D_{3} \overt D_{4}) \arrow{d}{R_{D_{1},D_{2},D_{3},D_{4}}}\\
		(C_{1} \ominus C_{3}) \overt (C_{2} \ominus C_{4}) \arrow{r}[yshift=-2ex,swap]{(f_{1}\ominus f_{3}) \overt (f_{2}\ominus f_{4})} & (D_{1}\ominus D_{3}) \overt (D_{2}\ominus D_{4})
	\end{tikzcd}
	\caption{$R$ is a natural transformation.}\label{fig:naturalR}
\end{figure}

\begin{Proposition}[{\cite[Proposition~6.35]{aguiar2010monoidal}}]	\label{prop:monoids}
The category $(\mathrm{Alg}_{\ominus},\overt,\pmb{\mathbb{C}}_{\overt})$ is a monoidal category. Indeed for any couple of horizontal algebra $\big(A,m_{\ominus}^{A},\eta_{A}\big)$ and $\big(B,m_{\ominus}^{B},\eta_{B}\big)$, the product $m^{A\overt B}\colon A\overt B \to A\overt B$ is defined
\begin{equation*}
m_{\ominus}^{A\overt B} := \big(m_{\ominus}^{A} \overt m_{\ominus}^{B}\big) \circ R_{A,B,A,B},\qquad \eta_{A\overt B}= \eta_{A}\overt\eta_{B}.
\end{equation*}
Moreover, the bicollection $\pmb{\mathbb{C}}_{\overt}$ is a an horizontal monoid
\begin{equation*}
m^{\overt}_{\ominus}\colon \ \pmb{\mathbb{C}}_{\overt} \ominus \pmb{\mathbb{C}}_{\overt} \to \pmb{\mathbb{C}}_{\overt}, \qquad \eta^{\overt}_{\ominus}\colon \ \pmb{\mathbb{C}}_{\ominus}\to \pmb{\mathbb{C}}_{\overt} ,
\end{equation*}
which are respectively a horizontal algebra and a horizontal unity.

The category $(\mathrm{CoAlg}_{\overt},\ominus,\pmb{\mathbb{C}}_{\ominus})$ is a monoidal category. Indeed for any couple of vertical comonoid $\big(A,m_{\overt}^{A},\eta_{A}\big)$ and $\big(B,m_{\overt}^{B},\eta_{B}\big)$, the product $\Delta^{A\ominus B}\colon A \to A\ominus B$ is defined
\begin{equation*}
\Delta_{\ominus}^{A\overt B} := R_{A,B,A,B}\circ \Delta^{A}_{\overt}\ominus\Delta^{B}_{\overt},\qquad \eta_{A\overt B}= \eta_{A}\overt\eta_{B}.
\end{equation*}
Moreover, the bicollection $\pmb{\mathbb{C}}_{\overt}$ is a an horizontal monoid
\begin{equation*}
m^{\overt}_{\ominus}\colon \ \pmb{\mathbb{C}}_{\overt} \ominus \pmb{\mathbb{C}}_{\overt} \to \pmb{\mathbb{C}}_{\overt}, \qquad \eta^{\overt}_{\ominus}\colon \ \pmb{\mathbb{C}}_{\ominus}\to \pmb{\mathbb{C}}_{\overt},
\end{equation*}
which are respectively a horizontal algebra and a horizontal unity.
\end{Proposition}

\begin{Definition}\label{def:PROS}We call PROS a monoid in the monoidal category $(\mathrm{Alg}_{\ominus}, \overt, \pmb{\mathbb{C}}_{\overt})$. That is an horizontal monoid $\big(C,m^{C}_{\ominus}, \eta^C_{\ominus}\big)$, endowed with a couple of bicollections morphisms
\[m^{C}_{\overt}\colon \ C\overt C\rightarrow C, \qquad \eta^C_{\overt}\colon \ \pmb{\mathbb{C}}_{\overt}\to C, \]
defining a vertical monoidal structure on $C$. In addition, these morphisms are horizontal morphisms.
\end{Definition}
We recall that the same structure takes also the name of double monoid in the literature, see, e.g.,~\cite{aguiar2010monoidal}.

\subsection*{Acknowledgements}
The authors thank Kurusch Ebrahimi-Fard for many enlightening discussions. CB is supported by the DFG Research Unit FOR2402 and NG is funded by DAAD kurzstipendium. We thank the anonymous referees for their detailed reports.

\pdfbookmark[1]{References}{ref}
\LastPageEnding

\end{document}